\documentclass[11pt]{article}
\input epsf.tex
\usepackage{latexsym}
\usepackage{amsmath,amsthm,amsfonts,amssymb}
\usepackage{epsfig}
\newtheorem{pro}{Proposition}[section]
 \newtheorem{thm}[pro]{Theorem}
 \newtheorem{lem}[pro]{Lemma}

 \newtheorem{cnj}[pro]{Conjecture} 
 \newtheorem{cor}[pro]{Corollary}

\newtheorem{question}[pro]{Question}
\newtheorem{remark}[pro]{Remark}

\def\A{{\cal A}}
 
 \def\C{{\mathbb C}}

 \def\H{{\mathbb H}}
 
\def\g{{\it {\bf g}}}
 \def\P{{\cal P}}

 \def\M{{\cal M}}
\def\m{{\bar M}}
 \def\F{{\mathbb F}}
 \def\T{{\cal T}}

 \def\R{{\mathbb R}}
 \def\chix{{\raise.5ex\hbox{$\chi$}}}
 
\def\Z{{\mathbb Z}}

\def\hF{{\cal F}}

\def\hG{{\cal G}}

\def\hG{{\cal G}}
\def\hH{{\cal H}}
\def\hS{{\cal S}}

\def\hK{{\cal K}}

\def\gF{{\tilde F}}
\def\gH{{\tilde H}}
\def\gG{{\tilde G}}

\def\gM{{\tilde M}}

\def\gK{{\tilde K}}

\def\wF{{F}}
\def\wH{{H}}
\def\wG{{G}}

\def\wK{{K}}

\def\Tree{{Tree}}
\def\Free{{F}}

\def\X{{\tilde X}}

\def\arccosh{{\textnormal{arccosh}}}

\begin{document}
\title{Weak Forms of the Ehrenpreis Conjecture and the Surface Subgroup Conjecture}
\author{Lewis Bowen\footnote{Research supported in part by a Max Zorn Postdoctoral Fellowship}}
\maketitle
\begin{abstract}
We prove the following:
\begin{enumerate}
\item 
Let $\epsilon>0$ and let $S_1,S_2$ be two closed hyperbolic surfaces. Then there exists locally-isometric covers ${\tilde S_i}$ of
$S_i$ (for $i=1,2$) such that there is a $(1+\epsilon)$ bi-Lipschitz
homeomorphism between ${\tilde S_1}$ and ${\tilde S_2}$ and both covers ${\tilde S_i}$ ($i=1,2$) have bounded injectivity radius. 
\item 
Let $\M$ be a closed hyperbolic $3$-manifold. Then there exists a map $j: S \to \M$ where $S$ is a surface of bounded injectivity radius and $j$ is $\pi_1$-injective local isometry onto its image.
\end{enumerate}
\end{abstract}
\noindent
{\bf MSC}: 20E07, 57M10, 57M50\\
\noindent
{\bf Keywords: } Ehrenpreis conjecture, surface subgroup, hyperbolic groups.










\section{Introduction}

\subsection{The Ehrenpreis Conjecture}

\begin{cnj}(The Ehrenpreis Conjecture) 
Let $\epsilon>0$ and let $S_1,S_2$ be two closed Riemann surfaces of the same genus. Then there exists finite-sheeted conformal covers ${\tilde S_i}$ of $S_i$ (for $i=1,2$) such that there is a $(1+\epsilon)$-quasiconformal homeomorphism between ${\tilde S_1}$ and ${\tilde S_2}$.
\end{cnj}

The Ehrenpreis conjecture was introduced in \cite{Ehr} where it was
proven in the case that $S_1$ and $S_2$ are tori. In the appendix to \cite{Gen} it is shown that in the remaining cases the conjecture is equivalent to the following.

\begin{cnj}(The hyperbolic Ehrenpreis conjecture)
Let $\epsilon >0$ and let $S_1, S_2$ be two closed hyperbolic surfaces. Then there exists finite-sheeted locally isometric covers  ${\tilde S_i}$ of $S_i$ (for $i=1,2$) such that there is a $(1+\epsilon)$ bi-Lipschitz homeomorphism between ${\tilde S_1}$ and ${\tilde S_2}$.
\end{cnj}

As a corollary to our first main theorem we will prove the above conjecture with ``finite-sheeted''
 replaced by ``bounded injectivity radius''. We need some terminology. A pair of pants $H$ is a surface homeomorphic to the $2-$sphere minus three open disks. It is boundary-ordered if the boundary components are ordered, in which case we may refer to the first boundary component of $H$ by $\partial_1 H$, the second by $\partial_2 H$ and so on. A labeled pants decomposition of a surface $S$ is a collection $\P$ of boundary-ordered pants embedded in $S$ whose interiors are pairwise-disjoint and whose union is all of $S$. We also require that ``labels'' match: if $\gamma$ is a simple closed curve in $S$ and $\gamma$ is in the boundary of $H_1,H_2 \in \P$ then we require $\gamma= \partial_i H_1=\partial_i H_2$ for some $i\in\{1,2,3\}$. Let $\P^* = \{\partial_i H :\, i=1,2,3$ and $ H\in\P\}$. We require every curve in $\P^*$ to be oriented (but no restrictions are put on the orientations).

\begin{thm}\label{thm:ehrenpreis case}(Main theorem: Ehrenpreis case)
Let $S_1,S_2$ be two closed hyperbolic surfaces and let $\epsilon>0$
be given. Suppose that $S_1$ and $S_2$ are incommensurable. Then there exists an $L_0 >0$ such that for all $L> L_0$
there exists locally isometric covers $\pi_i:{\tilde
  S_i} \to S_i$ (for $i=1,2$) such that for $i=1,2$, ${\tilde S_i}$
has a labeled pants decomposition $\P_i$ and there is a
homeomorphism $h:{\tilde S_1}\to {\tilde S_2}$ such that
$h(\P_1)=\P_2$, $h(\P^*_1)=\P^*_2$ and for every curve $\gamma \in \P_1^*$,
\begin{eqnarray*}
|length(\gamma)-L| &\le & \epsilon,\\
|length(h(\gamma))-L| &\le &\epsilon,\\
|twist(\gamma)| &\le & {\hat T}\exp(-L/4),\\
|twist(h(\gamma))| &\le & {\hat T}\exp(-L/4) \textrm{ and }\\
|twist(\gamma)-twist(h(\gamma))| &\le & \epsilon\exp(-L/4).
\end{eqnarray*}
Here ${\hat T}$ is a constant that depends only on $S_1,S_2$ and $\epsilon$.
\end{thm}
For the definition of $twist(\gamma)$, see section \ref{section:pants decompositions}.

\begin{remark} The covers ${\tilde S_i}$ constructed to prove this theorem are analytically infinite, genus 0 and without boundary. However, they have bounded injectivity radius.
\end{remark}

\begin{remark} The significance of the number $\exp(-L/4)$ lies in the fact that if $P$ is a hyperbolic pair of pants with geodesic boundary in which all boundary components have length $L$ then the distance between any two distinct boundary components is $2\exp(-L/4) + O(\exp(-3L/4))$.
\end{remark}


We also prove:
\begin{thm}\label{thm:bilipschitz}
Let $\epsilon_0 > 0$. Then there exists an $\epsilon, L_1 >0$ such that if two surfaces ${\tilde S_i}$ ($i=1,2$) satisfy the conclusion of
theorem \ref{thm:ehrenpreis case} with $L > L_1$ then there is a $(1+\epsilon_0)$ bi-Lipschitz homeomorphism $h:{\tilde S_1} \to {\tilde S_2}$. 
\end{thm}

The above two theorems imply the following weak form of the Ehrenpreis conjecture:
\begin{thm}
Let $\epsilon >0$ and let $S_1, S_2$ be two closed hyperbolic surfaces. Then there exists locally isometric covers  ${\tilde S_i}$ of $S_i$ (for $i=1,2$) such that
\begin{itemize}
\item there is a $(1+\epsilon)$ bi-Lipschitz homeomorphism between ${\tilde S_1}$ and ${\tilde S_2}$ and
\item ${\tilde S_1}$, ${\tilde S_2}$ have bounded injectivity radius.
\end{itemize}
\end{thm}

\begin{question}\label{question:ehrenpreis case} 
 Can the conclusions to theorem \ref{thm:ehrenpreis case} be strengthened so that ${\tilde S_i}$ is closed?
\end{question}
We will show (in part \ref{part:tree tilings}) that given $L,\epsilon, S_1,S_2$ the above question is equivalent to a linear programming problem. We intend to study this problem in more detail in a future paper.


\subsection{The Surface Subgroup Conjecture}

\begin{cnj}(The Surface Subgroup Conjecture)
Let $\M$ be a closed hyperbolic $3$-manifold. Then there exists a $\pi_1$-injective map $j:S\to \M$ from a closed surface $S$ of genus at least 2 into $\M$.
\end{cnj}

Perhaps the original motivation for the surface subgroup conjecture is its relationship to the Virtual Haken Conjecture. The latter states that every irreducible closed $3$-manifold $\M$ with infinite fundamental group has a finite sheeted cover which contains an embedded incompressible closed surface. It was first stated in \cite{Wal}. The surface subgroup conjecture is an immediate consequence of the Virtual Haken conjecture and one may hope that it is a stepping stone towards the Virtual Haken conjecture.






We will prove the following main theorem.
\begin{thm}\label{thm:surface subgroup case}(Main theorem: surface subgroup case)
Let $\M$ be a closed hyperbolic $3$-manifold. Suppose that $\M$ does not contain a closed totally geodesic immersed surface. Let $\epsilon >0$ be
given. Then there exists an $L_0$ such that for all $L > L_0$ the
following holds. There exists a map $j:S \to \M$ from a surface
$S$ that has a labeled pants decomposition $\P$ such that for every curve $\gamma \in \P^*$
\begin{eqnarray*}
|length_j(\gamma)-L| &\le& \epsilon, \\
|twist(\gamma)| &\le& {\hat T}\exp(-L/4) \textrm{ and }\\
|\Im(twist_j(\gamma))| &\le& \epsilon\exp(-L/4).
\end{eqnarray*}
Here ${\hat T}$ is a constant depending only on $\M$ and $\epsilon$. $length(\gamma)$ and $twist(\gamma)$ denote the complex length and complex twist parameter of $\gamma$ with respect to $j$ (see section \ref{section:pants decompositions}).  
\end{thm}

We will also prove:
\begin{thm}\label{thm:incompressible}
There exist positive numbers $\epsilon_0, L_1>0$ such that for any $0<\epsilon  <
\epsilon_0$ and $L>L_1$ if $j:S \to \M$ is a map from a surface into a hyperbolic $3$-manifold $\M$ satisfying the conclusion of theorem \ref{thm:surface subgroup case} then $j$ is $\pi_1$-injective.
\end{thm}

\begin{question}\label{question:surface subgroup case} Can the surface $S$ in theorem \ref{thm:surface subgroup case} be chosen to be closed?
\end{question}

An affirmative answer would imply the surface subgroup conjecture. We will show (in part \ref{part:tree tilings}) that given $L,\epsilon, \M$ the above question is equivalent to a linear programming problem.


\subsection{A word on the proof and organization of the paper}

The proof of the main theorems rely on what we call the isometry construction. It is a method for perturbing a given isometry into a given discrete group using the horocyclic flow. It is sketched in section \ref{section:sketch}. In part \ref{part:isometry construction}, bounds on the translation distance and position of the perturbed isometry are proven. These bounds are then used to prove the main theorems in section \ref{section:hexagon H}. In part \ref{part:tree tilings}, it shown that any specific instance of questions \ref{question:ehrenpreis case} and \ref{question:surface subgroup case} is equivalent to a certain linear programming problem. We are, as of yet, unable to show that the linear programming has a solution; though we can show (unpublished) that the corresponding system of linear equations has a solution space of relatively small codimension.

The proofs of theorems \ref{thm:bilipschitz} and
\ref{thm:incompressible} are handled in parts \ref{part:bilipschitz} and \ref{part:incompressible}
respectively. 

{\bf Acknowledgements}: We would like to thank Joel Hass for many
encouraging conversations.

\tableofcontents

\part{Background and Notation}\label{part:background}

\section{Estimate Notation}

In this paper, we will have several variables $T, \epsilon, \delta$, etc. However the variable $L$ will be treated in a special way. If $f$
is a function of $L$ and $x$ is a quantity that may depend on several
variables (including $L$) then the
notation $x=O(f(L))$ means there exist positive constants $k, L_0$ that do not
depend on $L$ but may depend on other variables such that for all $L > L_0$
\begin{displaymath}
|x| \le kf(L).
\end{displaymath}

If $x= y + O(f(L))$ and $z = w + O(g(L))$ and if $f(L)/y \to 0$ as $L
\to \infty$ and $g(L)/w \to 0$ as $L \to \infty$ then
\begin{displaymath}
\frac{x}{z} - \frac{y}{w} = O(f(L)/w + yg(L)/w^2).
\end{displaymath}

We use the following to express the above:

\begin{eqnarray*}
\frac{x}{z} = \frac{y + O(f(L))}{w + O(g(L))} = \frac{y}{w} +
O(f(L)/w + yg(L)/w^2)
\end{eqnarray*}

We will write $f(L) \sim g(L)$ to mean
\begin{displaymath}
\lim_{L \to \infty} \, \frac{f(L)}{g(L)} = 1.
\end{displaymath}

We will write $f(L) \approx g(L)$ to mean that there exist positive constants
$k_1,k_2,L_0$ such that for all $L> L_0$
\begin{displaymath}
k_1f(L) \le g(L) \le k_2f(L).
\end{displaymath}

\section{Groups of isometries of Hyperbolic Space}

Throughout this paper (unless explicitly stated otherwise), we will use the upperhalf space model of
$3$-dimensional hyperbolic space which we denote by $\H^3$. In this
model,
\begin{displaymath}
\H^3 = \{(z,t): z \in \C, t >0\}
\end{displaymath}
is equipped with the metric $ds^2 = (|dz|^2 + dt^2)/t^2$. See, for example \cite{Rat} or \cite{Fen} for more details.

The hyperbolic plane $\H^2 = \{(x,t): x \in \R , t>0\} \subset \H^3$
is isometrically embedded in $\H^3$. We identify the group of
orientation-preserving isometries of $\H^3$ with
$PSL_2(\C)$. $PSL_2(\C)$ acts on the complex plane $\C \times \{0\}
\subset \C \times [0,\infty)$ by fractional
linear transformations as follows:
\begin{displaymath}
\left[\begin{array}{cc}
a & b \\
c & d \end{array}\right]z = \frac{az+b}{cz+d}.
\end{displaymath}
The action of $PSL_2(\C)$ on $\H^3$ is defined by extending the above
action to $\H^3$ via the rule that for every $A \in PSL_2(\C)$, the
action of $A$ on $\C \times [0,\infty)$ takes
  semi-circles orthogonal to the boundary to semi-circles orthogonal to
  the boundary. The stabilizer of a point is equal to $SU(2)$, so we may identify $\H^3$ with $PSL_2(\C)/SU(2)$. Also we identify the positively oriented frame bundle of $\H^3$
  with $PSL_2(\C)$.

$\H^2$ is stabilized by $PSL_2(\R)$ which we
identify as the orientation-preserving isometry group of the plane
$\H^2$. $SO(2) < PSL_2(\R)$ is the stabilizer of a point in $\H^2$. So we identify $\H^2$ with $PSL_2(\R)/SO(2)$. In this way, we may also identify the unit tangent bundle of $\H^2$ with $PSL_2(\R)$.

Although we do most of our calculations in the upperhalf space model, the figures are often drawn in the Poincare model (see \cite{Rat} for a description of this model).

A discrete group $G$ of $PSL_2(\C)$ (or $PSL_2(\R)$) is a subgroup whose topology is discrete as a subspace of $PSL_2(\C)$ (or $PSL_2(\R)$). $G$ is cocompact in $PSL_2(\C)$ if the quotient space $PSL_2(\C)/G$ is compact.


\subsection{Displacements}

What follows is covered in more detail in \cite{Fen}. If $g$ is an orientation-preserving hyperbolic isometry of $\H^3$ then there is a unique geodesic called the axis of $g$ (and denoted here by $Axis(g)$) that is preserved under $g$. If $u,v $ are the endpoints of $Axis(g)$ on the boundary at infinity then an orientation of $Axis(g)$ is specified by an ordering of $\{u,v\}$. We associate an element $\mu(g,(u,v)) = \mu(g) = \mu \in \C \mod 2\pi i$ with $g$ and an orientation of its axis in the following way. There exists a unique orientation preserving isometry $A$ such that $Au=0$ and $Av = \infty$ (in the upper half space model). Then $AgA^{-1}$ has the form

\begin{displaymath}
g=\left[ \begin{array}{cc}
z & 0 \\
0 & z^{-1}
\end{array} \right]
\end{displaymath}
where $z \in \C$ is different from $0$ and $1$. We define $\mu$ by the equation
\begin{displaymath}
z = exp(\mu/2).
\end{displaymath}
$\mu$ is called the {\bf displacement} of $g$ (relative to the orientation on $Axis(g)$). Note that $\mu(g,(u,v))=-\mu(g,(v,u))$. In the sequel, we may write $\mu(g)$ if the orientation is understood. If $Axis(g)$ is oriented from the repelling fixed point of $g$ to its attracting fixed point, then we call $\mu(g)$ the {\bf complex translation length of $g$} and denote it by $tr.length(g)$. 

Note

\begin{equation}\label{eqn: trace formula}
\cosh(\mu(g)/2) = \frac{z + z^{-1}}{2} = trace(g)/2.
\end{equation}

\subsection{Fixed Points}\label{subsection:fixed points}

The fixed points of

\begin{displaymath}
\left[ \begin{array}{cc}
a & b \\
c & d
\end{array} \right] \in PSL_2(\C)
\end{displaymath}
acting by fractional linear transformations on $\C$ are 
\begin{equation}\label{eqn:fixed points}
\frac{ a - d \pm \sqrt{(a+d)^2 - 4} }{ 2c}
\end{equation}

\section{Trigonometry}

\subsection{Cross Ratio}\label{subsection:cross ratio}

The cross ratio $R$ of $(a,b,c,d) \in \C^4$ is defined by

\begin{equation}\label{eqn:cross ratio}
R(a,b,c,d) = \frac{ (a-c)(b-d) }{(a-d)(b-c) }.
\end{equation}

The cross ratio is invariant under the action of $PSL_2(\C)$ by fractional linear transformations on $\C$. Note that
\begin{displaymath}
R(b,a,c,d) = \frac{1}{R(a,b,c,d)}
\end{displaymath}

and $R(a,b,c,d) = R(c,d,a,b)$.

\subsection{Double Crosses}\label{subsection:doublecross}
The material in this subsection is detailed more thoroughly in \cite{Fen}. Suppose $u, u'$ are the endpoints of a geodesic $\gamma_1$ oriented from $u$ to $u'$ and $v,v'$  are the endpoints of a geodesic $\gamma_2$ oriented from $v$ to $v'$. Suppose also that $\gamma_3$ is a geodesic perpendicular to both $\gamma_1$ and $\gamma_2$. Let $w,w'$ be the endpoints of $\gamma_3$ which we assume is oriented from $w$ to $w'$. The triple $(\gamma_1,\gamma_2; \gamma_3)$ is called a {\bf double cross}. We define the {\bf width} $\mu(\gamma_1,\gamma_2;\gamma_3) = \mu \in \C/<2\pi i>$ of the double cross by the equation
\begin{displaymath}
exp(\mu) = R(u,v,w',w) = -R(u,v',w',w) = -R(u',v,w',w) = R(u',v',w',w).
\end{displaymath}
So
\begin{displaymath}\begin{array}{ll}
\exp(\mu(\gamma_1,\gamma_2;\gamma_3)) &= R(u,v,w',w)\\
                                &= 1/R(v,u,w',w)\\
                                &=\exp(-\mu(\gamma_2,\gamma_1;\gamma_3 )).
\end{array}
\end{displaymath}
Hence $\mu(\gamma_1,\gamma_2;\gamma_3) =
-\mu(\gamma_2,\gamma_1;\gamma_3)$. $\mu$ also satisfies the equation

\begin{equation}\label{eqn: displacement}
R(u,u',v,v') = \tanh^2(\mu/2).
\end{equation}

The latter equation determines $\mu$ only up to a sign. If we denote $\gamma_i$ with the opposite orientation by $-\gamma_i$, then we have

\begin{displaymath}
\mu(\gamma_1,\gamma_2;-\gamma_3) = -\mu(\gamma_1,\gamma_2;\gamma_3)
\end{displaymath}

and

\begin{displaymath}
\mu(-\gamma_1,\gamma_2;\gamma_3) = \mu(\gamma_1,\gamma_2;\gamma_3) + i\pi.
\end{displaymath}

 The real part of $\mu$ is the signed distance between $\gamma_1$ and $\gamma_2$. The imaginary part measures the amount of turning between $\gamma_1$ and $\gamma_2$. To be precise, $\mu$ is the displacement of the isometry $g$ with oriented axis $Axis(g)=(w,w')$ such that $g \gamma_1=\gamma_2$. 

Suppose as above that
\begin{displaymath}
R=R(u,u',v,v') = \tanh^2(\mu/2).
\end{displaymath}
Then
\begin{eqnarray*}
\frac{1+R}{1-R} &=&\cosh(\mu)\\
\frac{\sqrt{R}}{1-R} &=&(1/2)\sinh(\mu).
\end{eqnarray*}

\subsection{Right Angled Hexagons and Pentagons}\label{right angled polygons}

If $(S_1,S_2,S_3,S_4,S_5,S_6)$ is an {\it ordered} $6$-tuplet of oriented geodesics such that $S_i$ is orthogonal to $S_{i+1}$ and $S_i$ is not equal to $S_{i+2}$ for any $i$ modulo $6$, then it is called a {\bf right angled hexagon}. By orthogonal, we will mean that $S_i$ and $S_{i+1}$ intersect in $\H^3$ at a right-angle (this constrasts a little with \cite{Fen} where the word ``normal'' is used to allow the possibility that $S_i$ and $S_{i+1}$ share an endpoint at infinity). In the terminology of \cite{Fen} all the side-lines that we allow are proper.

Similarly, if $(S_1,S_2,S_3,S_4,S_5)$ is an {\it ordered} $5$-tuplet of oriented geodesics such that $S_i$ is orthogonal to $S_{i+1}$ and $S_i$ is not equal to $S_{i+2}$ for any $i$ modulo $5$, then it is called a {\bf right angled pentagon}. The following lemmas are classical. They appear in \cite{Fen}.

\begin{lem}\label{lem:right angled hexagons}
Let $(S_1,S_2,S_3,S_4,S_5,S_6)$ be a right-angled hexagon. Let $\sigma_i = \mu(S_{i-1}, S_{i+1}; S_i)$ denote the width of the double cross $(S_{i-1},S_{i+1};S_i)$ (see subsection \ref{subsection:doublecross} for defintions). Then the following relations hold.

\begin{enumerate}
\item The law of sines: $\frac{\sinh(\sigma_1)}{\sinh(\sigma_4)} = \frac{\sinh(\sigma_3)}{\sinh(\sigma_6)} = \frac{\sinh(\sigma_5)}{\sinh(\sigma_2)}$.

\item The law of cosines: $\cosh(\sigma_i) = \cosh(\sigma_{i-2}) \cosh(\sigma_{i+2}) + \sinh( \sigma_{i-2}) \sinh( \sigma_{i+2} ) \cosh(\sigma_{i+3})$ for all $i$ with indices considered modulo $6$.
\end{enumerate}
\end{lem}

\begin{lem}\label{lem:right angled pentagons}
If $(S_1,...,S_5)$ is a right-angled pentagon and $\sigma_n= \mu(S_{n-1},S_{n+1}; S_n)$ then
\begin{enumerate}
\item $\cosh(\sigma_n) = -\sinh(\sigma_{n-2})\sinh(\sigma_{n+2})$ for all $n$ mod $5$, and
\item $\cosh(\sigma_n) = -\coth(\sigma_{n-1})\coth(\sigma_{n+1})$ for all $n$ mod $5$.
\end{enumerate}
\end{lem}

Note that in the above, the (real) distance between $S_i$ and $S_{i+2}$ may be zero in which case $\sigma_{i+1}$ is purely imaginary.

We say that $\hS=(S_1,...,S_6)$ is {\bf standardly oriented} if the following holds. For any $i$, if $S_{i-1}$ and $S_{i+1}$ do not intersect then $S_i$ is oriented from its intersection with $S_{i-1}$ to its intersection with $S_{i+1}$. Otherwise let $u_j$ be a unit tangent vector at the point of intersection in the direction of $S_j$ for $j=i-1,i,i+1$. Then we require that $(u_{i-1},u_{i+1},u_i)$ is a positively oriented basis for the tangent space at $x \in \H^3$.  

We will, at times, also use the term ``right-angled hexagon'' to denote a 6-sided polygon (in $\H^n$) such that every pair of adjacent sides meets at a right angle. To any right-angled hexagon $(S_1,...,S_6)$ as above there exists a canonical polygon with vertices $v_1,..,v_6$ where $v_i$ is the intersection of $S_i$ with $S_{i+1}$ for all $i$ mod $6$. We may abuse notation at times by confusing $(S_1,...,S_6)$ with this polygon.

\section{Nearly Symmetric Right-Angled Hexagons}\label{section:L/2 hexagons}

Let $\rho_1,\rho_3,\rho_5 \in \C$ be three numbers that do not depend on the variable $L$. Let $\hG=(\gG_1,..,\gG_6)$ be the standardly oriented right-angled hexagon with $\wG_j=L/2 + \rho_j/2+ i\pi$ for $j=1,3,5$. Here $\wG_j = \mu(\gG_{j-1},\gG_{j+1};\gG_j)$. Assume that $|\rho_j| \le \epsilon$ for all $j=1,3,5$ and that $L$ is large. We will call any hexagon $\hG$ satisfying these properties an $(L,\epsilon)$ nearly-symmetric hexagon. The purpose of this section is to estimate various quantities related to $\hG$. 

\begin{lem}\label{lem:L/2 hexagon}
For $k=2,4,6$
\begin{eqnarray*}
\cosh(\wG_k) &=& -1 - 2 \exp(-L/2 + \rho_{k+3}/2- \rho_{k+1}/2 - \rho_{k-1}/2) + O(\exp(-L))\\
\wG_k &=& 2\exp(-L/4 + \rho_{k+3}/4 - \rho_{k+1}/4 - \rho_{k-1}/4) + i\pi + O(\exp(-3L/4)).
\end{eqnarray*}
with indices $\mod 6$.

\end{lem} 

\begin{proof}

The law of cosines implies that
\begin{displaymath}
\cosh(\wG_1)=\cosh(\wG_3)\cosh(\wG_5)+\sinh(\wG_3)\sinh(\wG_5)\cosh(\wG_4).
\end{displaymath}
So,
\begin{eqnarray*}
\cosh(\wG_4) &=& \frac{\cosh(\wG_1)-\cosh(\wG_3)\cosh(\wG_5)}{\sinh(\wG_3)\sinh(\wG_5)}\\
             &=& -\coth(\wG_3)\coth(\wG_5) + \frac{\cosh(\wG_1)}{\sinh(\wG_3)\sinh(\wG_5)}.
\end{eqnarray*}

But,
\begin{eqnarray*}
\coth(\wG_3) &=& \frac{\exp(L/2 + \rho_3/2) + \exp(-L/2 - \rho_3/2)}{\exp(L/2 + \rho_3/2) - \exp(-L/2 - \rho_3/2)}\\
             &=& 1 +2\frac{\exp(-L/2 - \rho_3/2)}{\exp(L/2 + \rho_3/2) - \exp(-L/2 - \rho_3/2)}\\
             &=&1+ O(\exp(-L)).
\end{eqnarray*}
Similarly, $\coth(\wG_5)=1 + O(\exp(-L))$. Hence,
\begin{eqnarray*}
\cosh(\wG_4) &=& -1 + O(\exp(-L)) + \frac{-(1/2)\exp(L/2 + \rho_1/2)}{(1/4)\exp(L + \rho_3/2 + \rho_5/2) + O(1)  }\\
             &=& -1  + \frac{-(1/2)\exp(L/2 + \rho_1/2)}{(1/4)\exp(L + \rho_3/2 + \rho_5/2) }  + O(\exp(-L))\\
             &=& -1 - 2 \exp(-L/2 + \rho_1/2- \rho_3/2 - \rho_5/2) + O(\exp(-L)).
\end{eqnarray*}

This implies that
\begin{displaymath}
\wG_4 = 2\exp(-L/4 + \rho_1/4 - \rho_3/4 - \rho_5/4) + i\pi + O(\exp(-3L/4)).
\end{displaymath}

The other statements follow in a similar manner.

\end{proof}

\begin{cor}\label{cor:M}
Suppose that in the above lemma $\rho_1=\rho_3=\rho_5=0$. Let $M(L)$ be the real part of $\wG_2 = \wG_4=\wG_6$. Then
\begin{displaymath}
M(L) = 2\exp(-L/4) + O(\exp(-3L/4)).
\end{displaymath}
\end{cor}
{\bf Remark}: If $P$ is a hyperbolic 3-holed sphere with geodesic boundary components all of length $L$ then the distance between any two distinct components is $M(L)=2\exp(-L/4)+O(\exp(-3L/4))$. This can be seen by considering that $P$ canonically decomposes into the union of two isometric right-angled hexagons by cutting $P$ along the three shortest arcs between distinct boundary components.

The proof of the next lemma is similar to that of the one above so we omit it.
\begin{lem}\label{lem:trlength}
Let $\rho_1, \rho_3 \in \C$ such that $|\rho_i|<\epsilon$. Let $\hG=(\gG_1,..,\gG_6)$ be the standardly oriented right-angled hexagon with $\gG_1=L/2 + \rho_1/2 + i\pi$, $\gG_3 = L/2 + \rho_3/2 + i\pi$, $\gG_2=2\exp(-L/4) + i\pi + O(\exp(-L/2))$. Then $\gG_5=L/2 + \rho_1/2 + \rho_3/2 + i\pi + O(\exp(-L/4))$.
\end{lem}

\subsection{Altitudes}\label{subsection:altitudes}

An {\bf altitude} of a right-angled hexagon $\hH$ is a geodesic that is perpendicular to two opposite sides of the hexagon $\hH$. If $\hH$ is a convex planar hexagon it is known (\cite{Bus}) that the three altitudes intersect in a single point and thus decompose $\hH$ into six trirectangles (convex $4$-gons with three right angles).

Let $\hG$ be the hexagon defined above. Let $\hK=(\gK_1,...,\gK_5)$ be the standardly oriented right-angled
pentagon defined by $\gK_k=\gG_k$ for $k=1,2,3,4$ and $\gK_5$ is the
common perpendicular of $\gG_1$ and $\gG_4$ (so it is the altitude between $\gG_1$ and $\gG_4$). If we let $\wK_k=\mu(\gK_{k-1},\gK_{k+1};\gK_k)$ (for all $k$ mod 5) then $\wK_k=\wG_k$ for $k=2,3$. We obtain estimates for the widths of $\hK$ in the next lemma.

\begin{lem}\label{lem:altitudes}
The widths of the pentagon $\hK$ satisfy the following estimates.
\begin{eqnarray*}
\wK_5 &=& L/4 + \log(2) + \rho_5/4 + \rho_3/4 - \rho_1/4 + i\pi + O(\exp(-L/2)).\\
\wK_1 &=& L/4 + \rho_5/4-\rho_3/4-\rho_1/4 + i\pi + O(\exp(-L/2)).\\
\wK_4 &=& \exp(-L/4 -\rho_5/4 - \rho_3/4 +\rho_1/4) + i\pi + O(\exp(-3L/4)).
\end{eqnarray*}

\end{lem}
\begin{proof}
By the right-angled pentagon identities lemma \ref{lem:right angled pentagons} we have
\begin{eqnarray*}
\cosh(\wK_5) &=& -\sinh(\wK_2)\sinh(\wK_3)\\
             &=& -\sinh(\wG_2)\sinh(\wG_3)\\
             &=& -\big(-2\exp(-L/4 + \rho_{5}/4 - \rho_{3}/4 - \rho_{1}/4) + O(\exp(-3L/4))\big)\\
&&\times \big(-(1/2)\exp(L/2 + \rho_3/2) + O(\exp(-L/2))\big)\\
             &=&-\exp(L/4 +\rho_5/4 + \rho_3/4 - \rho_1/4) + O(\exp(-L/4)).
\end{eqnarray*}
Thus
\begin{displaymath}
\wK_5 = L/4 + \log(2) + \rho_5/4 + \rho_3/4 - \rho_1/4 + i\pi + O(\exp(-L/2)).
\end{displaymath}
The estimates for $\wK_5$ follow. Note this implies $\coth^2(\wK_5)=1 + 1/\sinh^2(\wK_5) = 1 +
\exp(-L/2 -\rho_5/2 - \rho_3/2 +\rho_1/2) + O(\exp(-L))$. So $\coth(\wK_5)=1 +(1/2)\exp(-L/2 -\rho_5/2 - \rho_3/2 +\rho_1/2) + O(\exp(-L)) $. Since $\wK_2=\wG_2$
we have $\coth(\wK_2)=(1/2)\exp(L/4+\rho_5/4-\rho_3/4-\rho_1/4) + O(\exp(-L/4))$.

The pentagon identities lemma \ref{lem:right angled pentagons} implies
\begin{eqnarray*}
\cosh(\wK_1) &=& -\coth(\wK_2)\coth(\wK_5)\\
             &=& -(1/2)\exp(L/4 + \rho_5/4-\rho_3/4-\rho_1/4) + O(\exp(-L/4)).
\end{eqnarray*}
The estimate for $\wK_1$ follows. Note that $\coth^2(\wK_3)=1 +
1/\sinh^2(\wK_3)=1 + 4\exp(-L-\rho_3)+O(\exp(-2L))$. So $\coth(\wK_3)=1+2\exp(-L-\rho_3)+O(\exp(-2L))$. The pentagon identities lemma \ref{lem:right angled pentagons} implies
\begin{eqnarray*}
\cosh(\wK_4) &=& -\coth(\wK_3)\coth(\wK_5)\\
             &=& -    [1+2\exp(-L-\rho_3)+O(\exp(-2L))]\\
&&\times             [1+(1/2)\exp(-L/2 -\rho_5/2 - \rho_3/2 +\rho_1/2) +
             O(\exp(-L))]\\
             &=& -1 - (1/2)\exp(-L/2 -\rho_5/2 - \rho_3/2 +\rho_1/2) +
             O(\exp(-L)).
\end{eqnarray*}
The estimates for $\wK_4$ follow.

\end{proof}

\begin{lem}\label{lem:midpoint distances}
Suppose the hexagon $\hG$ is defined as above and $\rho_k=0$ for $k=1,3,5$. Let $p_k$ be the intersection point $\gG_k\cap\gG_{k+1}$ for $k$ mod 6. Let $m_k$ be the midpoint of $\overline{p_{k-1}p_k}$. Then
\begin{displaymath}
\cosh(d(m_1,m_3)) = 3/2 + O(\exp(-L/2)).
\end{displaymath}
\end{lem}

\begin{proof}
Consider the planar 4-gon with vertices $m_1,p_1,p_2,m_3$. It has right angles at $p_1$ and $p_2$. We use the formulas for the convex quadrangles with two right angles (\cite{Fen} page 88) to obtain
\begin{eqnarray*}
\cosh(\overline{m_1m_3}) &=& -\sinh(\overline{m_1p_1})\sinh(\overline{p_2m_3})+\cosh(\overline{m_1p_1})\cosh(\overline{p_2m_3})\cosh(\overline{p_1p_2})\\
                         &=&-\sinh^2(L/4)+\cosh^2(L/4)\cosh(\overline{p_1p_2})\\
                         &=&1 + \cosh^2(L/4)[\cosh(\overline{p_1p_2})-1]\\
                         &=&1 + (1/4)\exp(L/2)(2\exp(-L/2)) + O(\exp(-L/2))\\
                         &=&3/2 + O(\exp(-L/2)).
\end{eqnarray*}
\end{proof}

\section{Labeled Pants Decompositions}\label{section:pants decompositions}

In this section we define $length_j(\gamma)$ and $twist_j(\gamma)$ where $\gamma \in \P^*$, $\P$ is a labeled pants decomposition of a hyperbolic surface $S$, and $j:S \to \M$ is a map into $\M$ (either a hyperbolic $3$-manifold or the Cartesian product of two hyperbolic surfaces). 

Suppose $\gamma \in S$ is such that $\gamma = \partial_k H_1 = \partial_k H_2$ for some $k \in \{1,2,3\}$ and $H_1,H_2 \in \P$. Assume that $H_1$ is on the left of $\gamma$ and $H_2$ is on the right side of $\gamma$. For $i=1,2$ let $m_i$ denote the shortest path in $H_i$ between $\partial_{k+1} H_i$ and $\partial_k H_i$ (indices mod 3). 

\begin{figure}[htb] 
\begin{center}
 \ \psfig{file=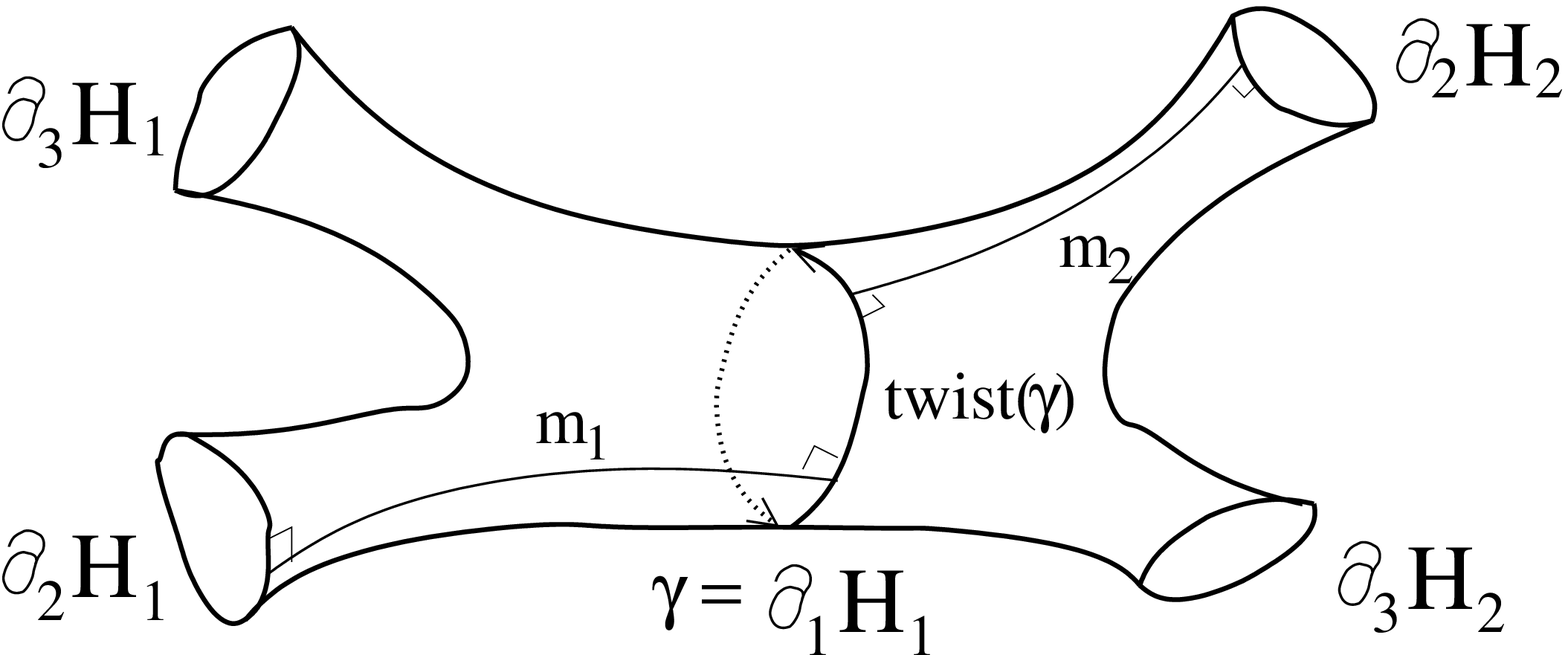,height=2in,width=3.5in}
 \caption{The twist parameter of $\gamma$.}
 \label{fig:twistreal}
 \end{center}
 \end{figure}

Define $twist_0(\gamma)$ equal to the signed distance from $m_1$ to $m_2$ along $\gamma$.  See figure \ref{fig:twistreal}. Let $length(\gamma)$ denote the length of $\gamma$ with respect to the hyperbolic metric on $S$. Let $twist(\gamma)\in\R$ be such that $twist_0(\gamma)\equiv twist(\gamma) \mod length(\gamma)$ and $twist(\gamma)$ has the smallest possible absolute value. 

The definitions of $length_j(\gamma)$ and $twist_j(\gamma)$ are generalizations of the above. If $j:S \to \M$ is continuous and $\M=S_1\times S_2$ is the product of two hyperbolic surfaces, then for $i=1,2$, let $j_i:S \to S_i$ equal $j$ followed by projection. For any curve $\gamma \subset S$, let $length_{j_i}(\gamma)$ be the length of the geodesic homotopic to $j_i(\gamma)$ or zero if $j_i(\gamma)$ is null-homotopic. Let $length_j(\gamma)=(length_{j_1}(\gamma),length_{j_2}(\gamma))\in\R^2$. 


Suppose $\gamma =\partial_k H=\partial_k H'$ for $H,H' \in \P^*$. If $j_i$ restricted to $H \cup H'$ is $\pi_1$-injective, then after homotopy we may assume that it is a local isometry. The hyperbolic structure on $S_i$ pulls-back to a hyperbolic metric on $H \cup H'$. Define $twist_{j_i}(\gamma)$ to be the twist parameter of $\gamma$ with respect to this metric. If $j_i$ restricted to $H \cup H'$ is not $\pi_1$-injective then we do not define $twist_{j_i}(\gamma)$. Let $twist_j(\gamma)=(twist_{j_1}(\gamma),twist_{j_2}(\gamma))\in\R^2$ when this makes sense.

\begin{figure}[htb]
 \begin{center}
 \ \psfig{file=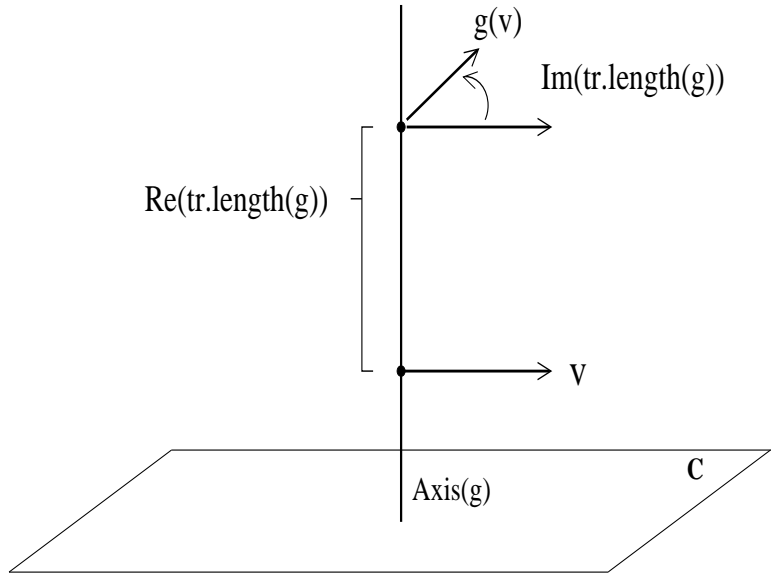,height=3in,width=4in}
 \caption{The complex translation length of $g$.}
 \label{fig:complex length}
 \end{center}
 \end{figure}

Let $\H^n$ denote $n$-dimensional hyperbolic space, $Isom^+(\H^n)$ denote the group of orientation preserving isometries of $\H^n$ and $d(x,y)$ be the distance between points $x,y \in \H^n$. If $g \in PSL_2(\C) = Isom^+(\H^3)$ is a hyperbolic (or loxodromic) isometry then its complex translation length $tr.length(g) \in \C$ is the complex number whose real part is the smallest number $r$ such that there is a $z \in \H^3$ such that $d(z,gz)=r$ and whose imaginary part measures the amount of rotation caused by $g$. To be precise, if $Axis(g)$ denotes the axis of $g$, $z \in Axis(g)$ and $v$ is a unit vector based at $z$ perpendicular to $Axis(g)$ then $Im(tr.length(g))$ is the angle from $\pi(v)$ to $gv$ where $\pi(v)$ equals $v$ parallel transported along $Axis(g)$ to lie in the tangent space of $g(z)$. See figure \ref{fig:complex length}.

Suppose $j_*:\pi_1(S) \to PSL_2(\C)$ is a representation. If $\gamma$ is a curve in $S$ then there is a unique conjugacy class $[\gamma]\subset \pi_1(S)$ representing it. Since conjugate elements of $PSL_2(\C)$ have the same translation length we can define the complex length of $\gamma$ (with respect to $j$) to be the complex translation length of any element in $j_*([\gamma])$. We denote it by $length_j(\gamma)$ or $length(\gamma)$ when $j$ is understood.

Let $\gamma \in \P^*$, $H_1,H_2 \in \P$ and $\gamma=\partial_k H_1=\partial_k H_2$. Let $S' = H \cup H'$. Assume that the representation $j_*$ restricts to a discrete faithful representation of $\pi_1(S')$ and that the image contains no parabolics. If this is not the case, then we do not define $twist_j(\gamma)$. Now let $\Gamma < PSL_2(\C)$ be the image $j_*(\pi_1(S'))$. By general homotopy theory there exists a map $j: S' \to \H^3/\Gamma$ that induces the representation $j_*$.  Assume that $H_1$ is on the left of $\gamma$ and $H_2$ is on the right. After homotoping $j$ if necessary we may assume that $j$ maps each curve in $\P^*$ onto a geodesic. We may also assume that there exists an oriented path $m_1 \subset S$ from $\partial_{k+1} H_1$ to $\gamma$ such that the length of $j(m_1)$ is as small as possible (over all such paths in $S$ and over all maps $j'$ homotopic to $j$ such that $j'$ maps each curve in $\P^*$ to a geodesic). Similarly we may also assume that there exists an oriented path $m_2 \subset S$ from $\gamma$ to $\partial_{k+1} H_2$ such that the length of $j(m_2)$ is as small as possible (over all such paths in $S$ and over all maps $j'$ homotopic to $j$ such that $j'$ maps each curve in $\P^*$ to a geodesic).

\begin{figure}[htb] 
\begin{center}
 \ \psfig{file=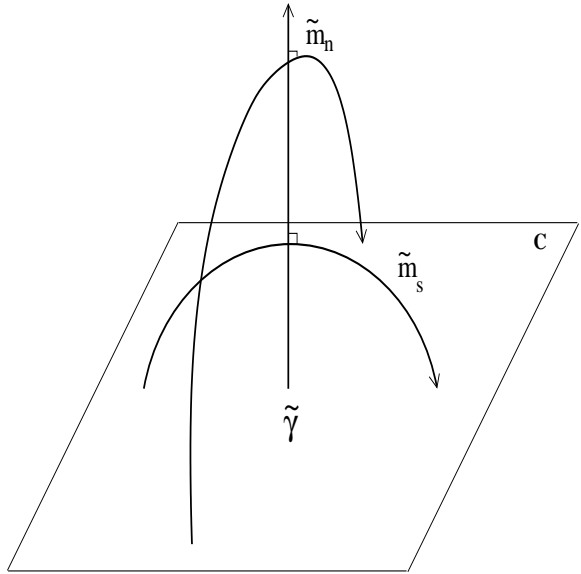,height=3in,width=3in}
 \caption{The complex twist parameter of $\gamma$. Here ${\tilde m_i}$ is the geodesic containing the lift of $j(m_i)$ ($i=1,2$).}
 \label{fig:twistcomplex}
 \end{center}
 \end{figure}

It follows that $j(m_1)$ and $j(m_2)$ are geodesic segments perpendicular to the images of the respective boundary curves and the middle curve. Now we lift the image of the middle curve $j(\gamma)$ and $j(m_1)$ and $j(m_2)$ up to the universal cover $\H^3$ as shown in figure \ref{fig:twistcomplex}. We assume the lifting is done so that the union of the lifts is connected and the distance between the lift of $j(m_1)$ and $j(m_2)$ along the lift of $\gamma$ is as small as possible.

Let ${\tilde m_i}$ be the geodesic containing the lift of $j(m_i)$. Let $g \in PSL_2(\C)$ be the isometry whose axis is the lift of the
image of $\gamma$ and such that $g({\tilde m}_1)={\tilde m_2}$
(where ${\tilde m_i}$ is the oriented geodesic containing the lift of
$j(m_i)$ for $i=1,2$). Let $Axis(g)$ be oriented from its repelling point to its attracting point.

We define $twist_j(\gamma)=\pm tr.length(g)$ where the sign is positive if the orientation on $Axis(g)$ agrees with the orientation induced by $\gamma$ and is negative otherwise. This generalizes the previous definition of $twist(\gamma)$ (when $S$ was a totally geodesic surface). For example, if the imaginary part of $twist(\gamma)$ is small then the surface is only ``lightly bent'' at $\gamma$.

Whether $\M$ is a product of surfaces or a $3$-manifold the definition of $length_j(\cdot)$ and $twist_j(\cdot)$ depends only on the homotopy class of $j$. Therefore, if $j_*:\pi_1(S) \to \pi_1(\M)$ is a homomorphism, then we may let $length_j(\cdot)$ and $twist_j(\cdot)$ be the length and twist parameter with respect to $j: S \to \M$ where $j$ is any map inducing $j_*$.

\section{The Horocyclic Flow}\label{section:horocyclic}

 Let
\begin{displaymath}
N_t = \left[ \begin{array}{cc}
1 & 0 \\
t & 1
\end{array} \right].
\end{displaymath}
and let $N = \{N_t : t\in\R\} < PSL_2(\R) <PSL_2(\C)$. Let ${\mathbb F}$ denote either $\R$ or $\C$. If $\Gamma < PSL_2({\mathbb F})$ is discrete group, the {\bf horocyclic flow} on the frame bundle $\Gamma \backslash PSL_2({\mathbb F})$ is the right action of $N$ on $\Gamma \backslash PSL_2({\mathbb F})$. 

Let $PSL^2_2(\R)$ denote $PSL_2(\R) \times PSL_2(\R)$. If $\Gamma_1, \Gamma_2 <PSL_2(\R)$ are discrete groups, then the {\bf diagonal horocyclic flow} on $(\Gamma_1 \times \Gamma_2)\backslash PSL^2_2(\R)$ is the action of the group ${\hat N} = \{(N_t,N_t) | t\in \R \}$ on $(\Gamma_1 \times \Gamma_2)\backslash PSL^2_2(\R)$.

We say that $\Gamma_1, \Gamma_2 < PSL_2(\R)$ are {\bf commensurable} if there exists an element $g \in PSL_2(\R)$ such that $g\Gamma_1 g^{-1} \cap \Gamma_2$ has finite index in both $g\Gamma_1 g^{-1}$ and $\Gamma_2$. In such a case, if $S_i=\H^2/\Gamma_i$ is a closed surface for $i=1,2$ then there exists a closed surface ${\tilde S} = \H^2/\Gamma$ and local isometries $\pi_i:{\tilde S} \to S_i$. In particular, the Ehrenpreis conjecture for $S_1$ and $S_2$ is trivial.

It seems likely that the following results are well-known. Except for the first statement below, we did not find them in the literature.

\begin{thm}\label{thm:horocyclic flow 2d}
Let $\H_2/\Gamma_1, \H_2/\Gamma_2$ be closed hyperbolic surfaces. Let $X=(\Gamma_1 \times \Gamma_2)\backslash (PSL_2(\R) \times PSL_2(\R))$. Then the following hold. 
\begin{enumerate}
\item Every orbit of the horocyclic flow in $\Gamma_i \backslash PSL_2(\R)$ is dense in $\Gamma_i \backslash PSL_2(\R)$ for $i=1,2$.

\item Every orbit of the diagonal horocyclic flow in $X$ is dense in $X$ {\it unless} $\Gamma_1$ and $\Gamma_2$ are commensurable.
\end{enumerate}

\end{thm}

\begin{proof}
The first statement was proven by Hedlund \cite{Hed}. Let $\Gamma= \Gamma_1 \times \Gamma_2$. Let $g \in \Gamma$. Then the closure of the ${\hat N}$-orbit of $\Gamma g$ equals
\begin{displaymath}
\overline{\Gamma g {\hat N}} = \overline{\Gamma \backslash g{\hat N} g^{-1}} g \subset \Gamma \backslash PSL^2_2(\R).
\end{displaymath}

 By Ratner's theorems on unipotent flows \cite{Ra5}, there exists a closed subgroup $P$ of $PSL_2(\R) \times PSL_2(\R)$ such that $g{\hat N}g^{-1} < P$ and 
\begin{displaymath}
\overline { \Gamma  g{\hat N}g^{-1}} = \Gamma P.
\end{displaymath}
Let $P_0$ be the component of $P$ containing the identity. By the classification of Lie subgroups of $PSL_2(\R)\times PSL_2(\R)$ and since $g{\hat N} g^{-1}$ is properly contained in $P_0$, $P_0$ must be conjugate to one of the following.
\begin{enumerate}
\item $B \times B$ (where $B < PSL_2(\R)$ is the set of upper triangular matrices).
\item $\{(X,X) | \, X \in PSL_2(\R)\}$.
\item $PSL_2(\R) \times PSL_2(\R)$.
\end{enumerate}
By Hedlund's result, for $k=1,2$, $\overline{\Gamma g{\hat N}g^{-1}}$ projects onto $\Gamma_k \backslash PSL_2(\R)$ under the canonical projection map. Therefore the first possibility cannot occur. If the second possibility occurs then let $\pi:PSL_2(\R) \to \H^2 = PSL_2(\R)/SO(2)$ by the quotient map. Then $(\pi\times\pi)(P)/\Gamma \subset (\H^2\times\H^2)/\Gamma$ is a closed hyperbolic surface and for $k=1,2$, the projection maps from $(\H^2\times\H^2)/\Gamma$ to $\H^2/\Gamma_k$ are local isometries. This implies that $\Gamma_1$ and $\Gamma_2$ are commensurable.

The third possibility is equivalent to the statement that the ${\hat N}$-orbit of $g$ in $\Gamma \backslash PSL^2_2(\R)$ is dense.

\end{proof}

The proof of the next theorem is similar to the previous theorem.

\begin{thm}\label{thm:horocyclic flow 3d}
Let $\M = \H^3/\Gamma$ be a closed hyperbolic $3$-manifold where $\Gamma < PSL_2(\C)$ is a discrete cocompact group. Then either $\M$ contains a totally geodesic immersed closed surface or every orbit of the horocyclic flow in $\Gamma \backslash PSL_2(\C)$ is dense. 
\end{thm}

\part{The Isometry Construction}\label{part:isometry construction}

\section{Sketch}\label{section:sketch}

Here we sketch the isometry construction in the 2-dimensional case. Suppose that $\Gamma < PSL_2(\R)$ is a discrete cocompact group. Let $a,b$ be distinct points in $\H^2$. Let $v_b$ be the unit vector based at $b$ that points away from $a$. Let $v_a$ be the unit vector based at $a$ that points towards $b$. Then there is a unique isometry $\gamma \in PSL_2(\R)$ that maps $a$ to $b$ and $v_a$ to $v_b$. 

\begin{figure}[htb]
 \begin{center}
 \ \psfig{file=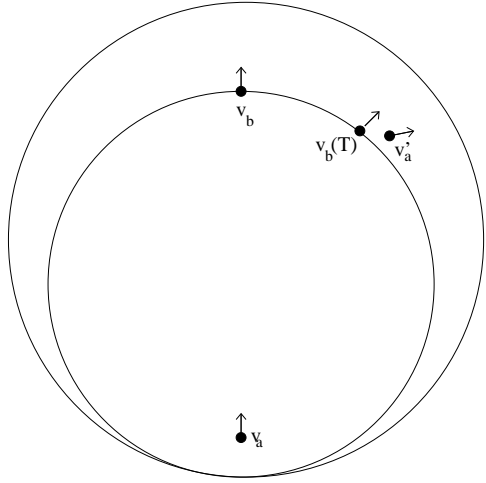,height=2.5in,width=2.5in}
 \caption{The vectors $v_a$, $v_b$, $v_b(T)$ and $g(v_a)$ in the Poincare model.}
 \label{fig:intuitive}
 \end{center}
 \end{figure}

Suppose we would like $\gamma$ to be an
element of $\Gamma$ but it is not. Then we perturb $\gamma$ so
that it is an element of $\Gamma$. To carry this out, let
$\epsilon>0$. Make the ray from $b$ through $a$. It limits on a
point $c$ on the circle at infinity. Let $h$ be the horocycle centered
at $c$ that passes through $b$. Now move the point $b$ along the
horocycle $h$ and carry the vector $v_b$ along with it. Let $b(t)$ and
$v_b(t)$ denote the point $b$ and the vector $v_b$ after time $t$. As
we move $v_b$ along the horocycle we look in its
$\epsilon$-neighborhood (with respect to some metric on the unit
tangent bundle of $\H^2$). As soon as we see a $\Gamma$ translate of
$v_a$ we stop. Let $v'_a$ be the translate of $v_a$ that we first
encounter. By definition then there is an isometry $g \in \Gamma$ such
that $g(v_a)=v'_a$. This isometry $g$ will be our new isometry, a
perturbed copy of $\gamma$. This is what we call the isometry
construction (theorem \ref{thm:isometry-3d}). See figure
\ref{fig:intuitive}. 


Since every orbit of the horocycle flow on $\Gamma \backslash
PSL_2(\R)$ is dense in $\Gamma \backslash PSL_2(\R)$ (theorem
\ref{thm:horocyclic flow 2d}) the time $T$ at which we stop moving the
point $b$ is bounded by a function of $\epsilon$ and $\Gamma$. In
particular, the bound does not depend on the points $a$ and $b$. 


Using this fact and some explicit calculations we will make estimates
regarding the translation length of $g$ and the location of its
axis. For instance, if the distance $d(a,b)=L$ between $a$ and $b$ is
very large, then the translation length of $g$ will be about
$\epsilon$-close to $L$. Just how large $L$ needs to be depends only on $\epsilon$ and $\Gamma$. Other estimates will prove useful in bounding the geometry of not just a single isometry $g \in \Gamma$ but interesting
subgroups of $\Gamma$, mainly $2$-generator free subgroups whose convex hull quotients are pairs of pants.

\section{The Isometry Construction}\label{section:isometry construction}


Given $L>0$, $T, \delta, \theta \in \R$ and $\nu \in \C$ let $\g =
\g(L, T, \nu, \delta, \theta) \in SL_2(\C)$ be defined by the following.

\begin{displaymath}\begin{array}{l}
\g=\left[ \begin{array}{cc}
\exp(L/2) & 0 \\
0 & \exp(-L/2) 
\end{array} \right]
\left[ \begin{array}{cc}
1 & 0 \\
T & 1
\end{array} \right]
\left[ \begin{array}{cc}
e^\nu & 0 \\
0 & e^{-\nu} 
\end{array} \right]
\left[ \begin{array}{cc}
\cos(\delta) & \sin(\delta) \\
-\sin(\delta) & \cos(\delta)
\end{array} \right]
\left[ \begin{array}{cc}
\exp(i\theta) & 0 \\
0 & \exp(-i\theta)
\end{array} \right].\\

\end{array}
\end{displaymath}

So,

\begin{equation}\label{eqn:g}
\g=\left[ \begin{array}{ccc}
\exp(L/2 + \nu + i\theta) \cos(\delta) & & \exp(L/2 + \nu - i\theta)\sin(\delta) \\
\exp(-L/2)[T e^{\nu+i\theta} \cos(\delta) - e^{-\nu + i\theta}\sin(\delta)] & & \exp(-L/2)[Te^{\nu-i\theta} \sin(\delta) + e^{-\nu-i\theta}\cos(\delta)] 
\end{array} \right].
\end{equation}

\begin{pro}\label{pro:translation lengths}
The translation length of $g= \g(L,T,\nu,\delta,\theta)$ satisfies
\begin{displaymath}
tr.length(g) =  L + 2\nu + 2i\theta + 2\log(\cos(\delta)) + O(\exp(-L)).
\end{displaymath}
Moreover, the constant implicit in the $O(\cdot)$ notation depends only on an upper bound for $|T|$ and an upper bound for $|\nu|$.
\end{pro}

\begin{proof}

From equation \ref{eqn:g} we see that
\begin{eqnarray*}
trace(g) &=& \exp(L/2 + \nu + i\theta)\cos(\delta) + \exp(-L/2 + \nu -
         i\theta)T\sin(\delta) + \exp(-L/2 - \nu - i\theta)\cos(\delta)\\
         &=& \exp(L/2 + \nu + i\theta)\cos(\delta) + O(\exp(-L/2)).
\end{eqnarray*}
Let $\mu$ be the displacement of $g$ when $Axis(g)$ is oriented from
its repelling fixed point to its attracting fixed point. Then,
\begin{eqnarray*}
\cosh(\mu(g)/2) &=& trace(g)/2\\
                &=& (1/2)\exp(L/2 + \nu + i\theta)\cos(\delta) + O(\exp(-L/2)) .
\end{eqnarray*}

This implies that
\begin{eqnarray*}
\mu(g)/2 &=& L/2 + \nu + i\theta + \log(\cos(\delta)) + O(\exp(-L)).
\end{eqnarray*}

Equivalently,
\begin{eqnarray*}
tr.length(g)&=&\mu(g) = L + 2\nu + 2i\theta + 2\log(\cos(\delta)) + O(\exp(-L)).
\end{eqnarray*}
The last statement is easy to check.

\end{proof}


\begin{thm}\label{thm:isometry-3d}(The Isometry Construction, 2d and 3d case)
Let $\F$ denote either $\R$ or $\C$. Let $\Gamma < PSL_2(\F)$ be a cocompact discrete torsion-free
orientation-preserving group. If $\F=\C$ then assume that $\H^3/\Gamma$ does not contain an immersed totally geodesic closed surface. Let
$\epsilon >0$. Let $I_r, I_i \subset [-\epsilon,\epsilon]$ be closed sets with nonempty interior.

Then there exists positive numbers $ {\hat T}={\hat T}(\Gamma,\epsilon,I_r,I_i)$ and $L_0=L_0(\Gamma,\epsilon,I_r,I_i)$ such that for every orientation-preserving isometry $A \in PSL_2(\F)$ and for every $L > L_0$ there exists parameters $\delta, \theta \in \R$ and $T, \nu \in \C$ satisfying all of the following.
\begin{itemize}
\item $g:=A^{-1}\g(L,T,\nu,\delta, \theta)A \in \Gamma$.
\item $|\nu|, |\delta|, |\theta|, |\Im(T)| < \epsilon$.
\item $|\Re(T)| < {\hat T}$.
\item $tr.length(g) \in L + I_r + iI_i$. 
\item If $\F=\R$ then $\Im(T)=\Im(\nu)=\theta=0$.
\end{itemize}
In the above, we have abused notation by identifying $\g$ with its projection to $PSL_2(\F)$.
\end{thm}

{\bf Remark}: To see how this is related to the sketch given in section \ref{section:sketch}, let $a$ be the point $(0,1)$ (in the upper half space model) and $b$ be the point $(0,\exp(L))$. Then $g$ is the isometry $\g(L,T,\nu,\delta,0)$.

\begin{proof}

For any $\nu \in C$ and $\delta \in \R$ define matrices $C_\nu,
R_\delta \in SL_2(\C)$ by
\begin{displaymath}\begin{array}{l}
C_\nu = \left[ \begin{array}{cc}
e^\nu & 0 \\
0 & e^{-\nu} 
\end{array} \right]\\
\\
R_\delta = \left[ \begin{array}{cc}
\cos(\delta) & \sin(\delta) \\
-\sin(\delta) & \cos(\delta)
\end{array} \right].
\end{array}
\end{displaymath}

Let $I'_r$ be a closed subset contained in the interior of $I_r$ such that $I'_r$ has nonempty
interior. Let $I'_i$ be a closed subset contained in the interior of
$I_i$ such that $I'_i$ has nonempty interior. Let $N_t \in SL_2(\C)$
be defined as in section \ref{section:horocyclic}. Let
$B=B(\epsilon,I'_r, I'_i) \subset PSL_2(\F)$ be the set of all
matrices of the form $N_t C_\nu R_\delta C_{i\theta}$ where
\begin{itemize}
\item $\delta, \theta \in \R$
\item $t, \nu, i\theta \in {\mathbb F}$,
\item $|t|, |\delta|, |\nu|, |\theta| <\epsilon$,
\item $\Re(2\nu+2\log(\cos(\delta))) \in I'_r$ and
\item $\Im(2\nu) + 2\theta \in I'_i$.
\end{itemize}
Note $B$ has nonempty interior in $PSL_2(\F)$. Theorem \ref{thm:horocyclic flow 2d} and the hypothesis that either $\F=\R$ or $\H^3/\Gamma$ does
not contain an immersed totally geodesic closed surface imply that
for every element $Z \in PSL_2(\F)$ the orbit 
\begin{displaymath}
\{\Gamma ZN_t | \, t \in \R \} \subset \Gamma \backslash PSL_2(\F)
\end{displaymath}
is dense in $\Gamma \backslash PSL_2(\F)$. Since $\Gamma \backslash PSL_2(\F)$ is compact this
implies that there exists a $T' > 0$ such that for all $Z_1,Z_2 \in
PSL_2(\F)$, there exists a $T$ with $|T| < T'$ such that $\Gamma Z_1
N_T \in \Gamma Z_2 B^{-1}$. We apply this fact to the isometries $Z_1=A^{-1}X$ and $Z_2 = A^{-1}$ where $X=C_{L/2}$.


So there exists a $T_0$ with $|T_0| <T'$ and a $g
\in \Gamma $ such that $A^{-1}X N_{T_0} = g A^{-1} (N_t C_\nu
R_\delta C_{i\theta})^{-1} $ for some $\delta, \theta \in \R$, $t, \nu \in \F$ with $N_t
C_\nu R_\delta C_{i\theta} \in B$. So,
\begin{displaymath}
g  = A^{-1} X N_{T_0} N_t C_\nu R_\delta C_{i\theta} A =
A^{-1}\g(L,T_0+t,\nu,\delta, \theta )A.
\end{displaymath}
By proposition \ref{pro:translation lengths}, $tr.length(g) = L+
2\nu + 2i\theta + 2\cos(\delta) +O(\exp(-L))$. Since $\Re(2\nu+2\cos(\delta)) \in
I'_r \subset int(I_r)$ and $\Im(2\nu) + 2\theta \in I'_i \subset
int(I_i)$ there exists an $L_0$ (that depends only on $T'$, $\epsilon$ and $\Gamma$) such that if $L>L_0$ then
\begin{displaymath}
tr.length(g) \in L + I_r + I_i.
\end{displaymath}
To finish, let ${\hat T}=T'+\epsilon$.
\end{proof}


The proof of the theorem below is similar to the one above.

\begin{thm}\label{thm:isometry-ehrenpreis}(The Isometry Construction, product of two surfaces case)
Let $\Gamma_1, \Gamma_2 < PSL_2(\R)$ be two cocompact discrete
torsion-free orientation-preserving groups. Assume that $\Gamma_1$ and $\Gamma_2$ are incommensurable. Let $\epsilon>0$. Let $I_1, I_2$ be two closed subsets of
$[-\epsilon,\epsilon]$ such that both $I_1$ and $I_2$ have
nonempty interior. Then there exists positive numbers $ {\hat T} = {\hat T}(\Gamma_1,\Gamma_2,\epsilon,I_1,I_2)$, $L_0=L_0(\Gamma_1,\Gamma_2,\epsilon,I_1,I_2) >0$
such that the following holds.

For every pair of orientation-preserving isometries $(A_1,A_2) \in PSL_2(\R) \times PSL_2(\R)$, for every $L>L_0$ and for $k=1,2$ there exists parameters $T_k,\nu_k, \delta_k \in \R$ such that all of the following statements hold.

\begin{itemize}
\item $g_k:=A^{-1}_k\g(L,T_k,\nu_k,\delta_k, 0)A_k \in \Gamma_k$.
\item $|T_1|, |T_2| < {\hat T}$.
\item $|\delta_k|, |\nu_k|, |T_1 - T_2|< \epsilon$.
\item $tr.length(g_k) \in L + I_k$.
\end{itemize}

\end{thm}

\section{The Hexagon $\hH$}\label{section:hexagon H}

{\bf Estimate Notation}: Throughout this section when we write $x=O(f(L))$ the constants implicit in the $O(\cdot)$ notation will depend only on upper bounds for $|T|$, $|\nu|$ and $|\delta|$.

For this section we fix quantities $L, {\hat T}, \delta,\theta, \alpha, \epsilon \in \R$, $T, \nu, \X, \gM \in \C$ satisfying the bounds
\begin{itemize}
\item $|T|<{\hat T}$,
\item $|\delta|, |\theta|, |\alpha|, |\nu|, |\Im(T)|| < \epsilon$,
\item $\X= \exp(L/2)+O(\exp(L/4))$ and
\item $\gM = e^\alpha M + O(\exp(-L/2))=2\exp(-L/4+\alpha)+O(\exp(-3L/2))$ where $M=M(L)$ is defined in corollary \ref{cor:M}.
\end{itemize}
Let $g=\g(L,T,\nu,\delta,\theta)$. Let $\hH=\hH(g,\X, \gM) = (\gH_1,...,\gH_6)$ be the standardly oriented right-angled hexagon satisfying the following (where $\wH_k = \mu(\gH_{k-1},\gH_{k+1};\gH_k)$ for all $k$ mod 6, see subsection \ref{subsection:doublecross}).

\begin{figure}[h]
 \begin{center}
 \ \psfig{file=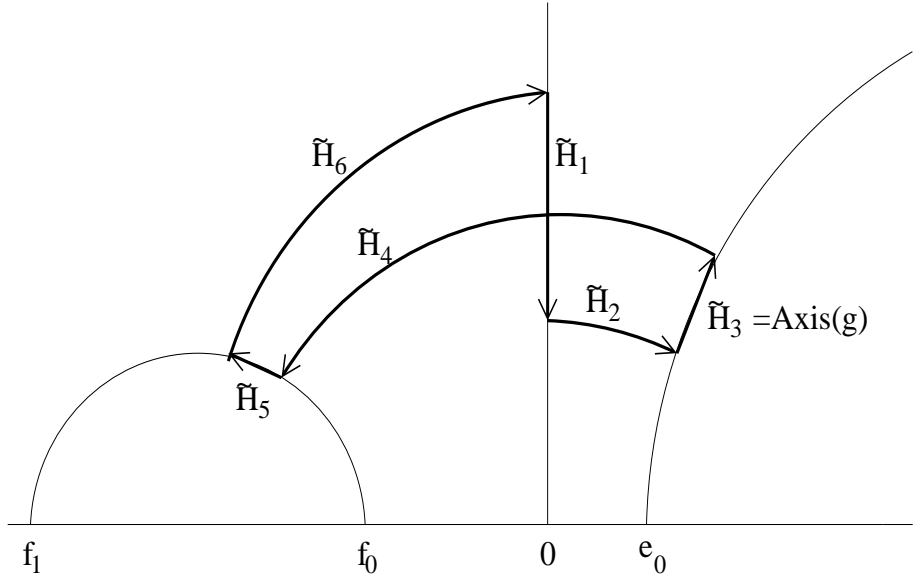, height=3in, width=4.75in}
 \caption{The hexagon $\hH$ (in the upperhalf space model).}
 \label{fig:hexagonH}
 \end{center}
 \end{figure}

\begin{enumerate}
\item $\gH_1$ is equal to the geodesic with endpoints $\{0,\infty\}$ (as unoriented geodesics).
\item $\gH_3$ is equal to the axis of $g$ (as unoriented geodesics). 
\item $\gH_6$ has endpoints $\pm \X$. 
\item $\wH_6=\gM + i\pi$.
\end{enumerate}
Hexagon $\hH$ is depicted in figure \ref{fig:hexagonH}.

\begin{thm}(Hexagon $\hH$ estimates)\label{thm:H estimates}
Assume that $|\tan(\delta)| \le 2\epsilon$ and $|e^{-2\nu}| \le
2$. Then:
\begin{eqnarray*}
\wH_2 &=& i\pi + O(\exp(-L/2)),\\
\wH_4 &=& \gM + O(\exp(-L/2)),\\
\wH_5 &=& (\sqrt{2}/2)\coth(\gM)\exp(-L/2)(T +\tau  ) +i\pi  + O(\exp(-L/2))
\end{eqnarray*}
where $\tau \in \C$ is such that $|\tau| \le 6\epsilon$.
\end{thm}

We now complete the proofs of the main theorems given the above result (which is proven in the next 4 subsections). Suppose $\M$ is a closed hyperbolic $3$-manifold not containing any totally geodesic closed immersed surfaces. Let $\Gamma < PSL_2(\C)$ be a discrete group such that $\M$ is isometric to $\H^3/\Gamma$. Suppose $S$ is a surface with a labeled pants decomposition $\P$. Suppose $H \in \P$ is such that $\partial_1 H \subset \partial S$. Suppose $j_*: \pi_1(S) \to \Gamma$ is a discrete representation such that the restriction of $j_*$ to $\pi_1(H) < \pi_1(S)$ is faithful and without parabolics. 
 Let $S' \supset S$ be the (topological) surface extending $S$ such that $S'$ has a labeled pants decomposition $\P'$ extending $\P$ so that $\P'-\P=\{H'\}$ and $\partial_1 H'=\partial_1 H$. We will use the above theorem to extend $j_*$ to a homomorphism $j'_*: \pi_1(S') \to \Gamma$ satisfying certain geometric bounds.

\begin{lem}
Let $\epsilon>0$ be given.
 Then there exists an $L_1$ (depending only on $\Gamma$ and $\epsilon$) such that if $L>L_1$ and
\begin{eqnarray*}
\Big| length_j\big(\partial_1 H\big)-L\Big| &\le& \epsilon\\
\end{eqnarray*}
then there exists a discrete representation $j'_*:\pi_1(S') \to \Gamma$ extending $j_*$ such that
 \begin{eqnarray*}
\Big|length_{j'}\big(\partial_k H' \big)-L\Big| &\le& \epsilon \textnormal{ for $k=1,2,3$},\\
\Big|\Re \big(twist_{j'}(\partial_1 H)\big)\Big| &\le& {\hat T}\exp(-L/4),\\
\Big|\Im \big(twist_{j'}(\partial_1 H)\big)\Big| &\le& \epsilon \exp(-L/4)
\end{eqnarray*}
where ${\hat T}$ is a number depending only on $\M$ and $\epsilon$.
\end{lem}

\begin{proof}
Let $I_-=[-\epsilon,-\epsilon/2]$ and $I_+=[\epsilon/2, \epsilon]$. Choose $L_1$ to be larger than
\begin{eqnarray*}
\max \, L_0(\Gamma,\epsilon,I_{\sigma_1},I_{\sigma_2})
\end{eqnarray*}
where $L_0(\cdot)$ is given by theorem \ref{thm:isometry-3d} and the maximum is overall all $\sigma_1,\sigma_2 \in \{+,-\}$. If necessary, choose $L_1$ larger so that for all $L>L_1$ the error estimates in theorem \ref{thm:H estimates} and lemmas \ref{lem:L/2 hexagon} and \ref{lem:trlength} are at most $\epsilon/8$. 

By general homotopy theory, there exists a map $j: S \to \M=\H^3/\Gamma$ inducing $j_*$. After homotoping $j$ we may assume that for $k=1,2,3$, $j$ maps $\partial_k H$ to a geodesic. We may also assume that there exists a path $m$ from $\partial_2 H$ to $\partial_1 H$ such that the length of $j(m)$ is as small as possible over all such path and over all maps $j'$ homotopic to $j$ that map the boundary of $H$ to geodesics. Thus $j(m)$ is a geodesic segment perpendicular to $j(\partial_1 H)$. 

Let ${\tilde m}$ and $\gamma$ be lifts of $j(m)$ and $j(\partial_1 H)$ respectively (so that ${\tilde m} \cup \gamma$ is a lift of $j(m) \cup j(\partial_1 H)$). Orient ${\tilde m}$ towards $\gamma$. Orient $\gamma$ to be consistent with the given orientation on $\partial_1 H$.  Let $\eta$ be the oriented geodesic containing ${\tilde m}$. Let $\Pi$ be the geodesic plane containing $\eta$ and $\gamma$. Orient $\Pi$ so that $(v_1,v_2)$ forms a positively-oriented bases for $\Pi$ at $\gamma \cap \eta$ where $v_1$ points in the direction of $\eta$ and $v_2$ points in the direction of $\gamma$. Let $p$ be a point on $\eta$ such that $d(p,\eta\cap\gamma)=M$ (where $M=M(L)$ is as in corollary \ref{cor:M}) and $p$ comes after $\gamma \cap \eta$ with respect to the orientation on $\eta$. See figure \ref{fig:hexagons}.

Let $A \in PSL_2(\C)$ be the orientation-preserving isometry that maps $\Pi$ to $\H^2$ (that is the plane bounded by the real line in the upperhalf-space model with the standard orientation), $\eta$ to the oriented geodesic from $-\exp(L/2)$ to $\exp(L/2)$, and $p$ to $(0,\exp(L/2))$. 

Let $\rho_1 \in \C$ be defined by $length(\partial H_1)=L+\rho_1$. Let $I_r=[\epsilon/2,\epsilon]$ or $[-\epsilon,-\epsilon/2]$ depending on whether $\Re(\rho_1)$ is negative or positive. Similarly, let $I_i=[\epsilon/2,\epsilon]$ or $[-\epsilon,-\epsilon/2]$ depending on whether $\Im(\rho_1)$ is negative or positive.

By the isometry construction theorem \ref{thm:isometry-3d}, there exists parameters $T,\nu,\delta,\theta$ such that
\begin{eqnarray*}
|T| &\le& {\hat T},\\
|\nu|, |\delta|, |\theta|, |\Im(T)| &\le & \epsilon,\\
g_3:=A^{-1}\g(L,T,\nu,\delta,\theta)A &\in& \Gamma\\
tr. length(g_3) &\in & L + I_r + iI_i.
\end{eqnarray*}

Let $\hH=(\gH_1,...,\gH_6)$ be the standardly oriented right-angled hexagon defined by

\begin{itemize}
\item $\gH_5$ is equal to $\gamma$ (as unoriented geodesics).
\item $\gH_3$ is equal to the axis of $g_3$ (as unoriented geodesics). 
\item $\gH_6$ is equal to $\eta$ (as unoriented geodesics). 
\item $\gH_1 \cap \gH_6 =p$ so $\wH_6=M + i\pi$ where $\wH_i=\mu(\gH_{i-1},\gH_{i+1};\gH_i)$ for all $i$ mod $6$. 
\end{itemize}

Let $\hG=(\gG_1,...,\gG_6)$ and $\hF=(\gF_1,...,\gF_6)$ be the standardly oriented right-angled hexagons defined by

\begin{itemize}
\item $\gG_1=\gF_1=\gamma$ but $\gG_1$ has the same orientation as $\gamma$ whereas $\gF_1$ has the opposite orientation,
\item $\gG_3=\gF_2=\eta$ (with orientation)
\item $\gG_3=\gF_2=Axis(g_3)$ (as unoriented geodesics)
\item $\wG_1=\wF_1=length(\partial_1 H)$,
\item $\wG_3=\wF_3=tr.length(g_3)$.
\end{itemize}

\begin{figure}[htb]
 \begin{center}
 \ \psfig{file=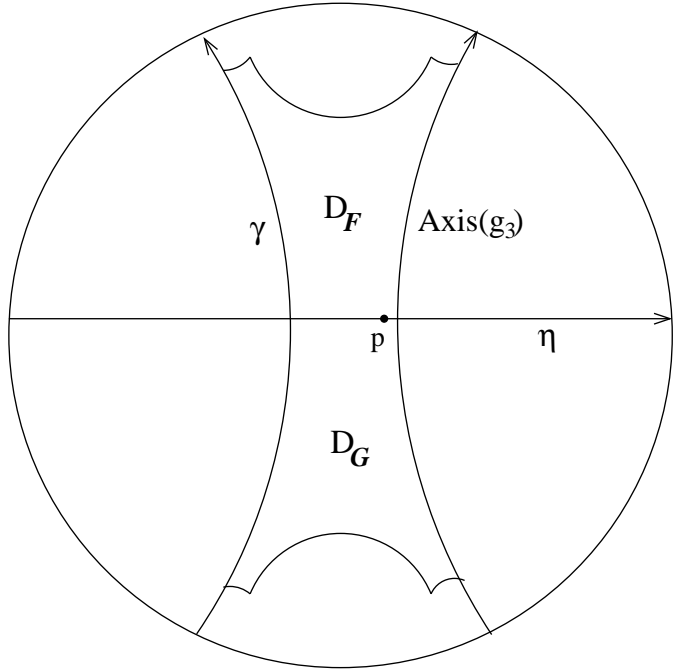,height=3.5in,width=3.5in}
 \caption{$\gamma$, $\eta$ and the disks $D_G$ and $D_F$ in the Poincare ball model.}
 \label{fig:hexagons}
 \end{center}
 \end{figure}

Here $\wG_i=\mu(\gG_{i-1},\gG_{i+1};\gG_i)$ and $\wF_i=\mu(\gF_{i-1},\gF_{i+1};\gF_i)$ for all $i$ mod $6$. Let $\rho_3, \rho_5 \in \C$ be defined by $tr.length(g_3)=L + \rho_3$, $\wG_5=L/2 + \rho_5/2 + i\pi$. The hypothesis and construction imply $|\rho_1|,|\rho_3| < \epsilon$. By definition, $\wG_2=\wH_4 + i\pi$. By the hexagon $\hH$ estimates theorem above, this implies $\wG_2 = M + O(\exp(-L/2))$. Lemma \ref{lem:trlength} implies $\rho_5 = \rho_1 + \rho_2 + O(\exp(-L/4))$. By the hypotheses on $L_1$, the error term is bounded by $\epsilon/8$. This implies that $|\rho_5| < \epsilon$. 

Let $D_G, D_F \subset \H^3$ be two disks with boundaries $\partial \hG$ and $\partial \hF$. Here $\partial \hG$ is the piecewise geodesic cycle with vertices $v_i = \gG_i \cap \gG_{i+1}$ for $i$ mod $6$ (and similarly for $\partial \hF$). 

Let $g_1 \in \Gamma$ be the hyperbolic element with axis $\gamma$ and translation length equal to $length(\partial_1 H)$. Clearly $(D_G \cup D_F)/<g_1,g_3>$ is a pair of pants $H'$ where $<g_1,g_3>$ denotes the group generated by $g_1$ and $g_3$. We order the boundary components so that $\partial_1 H'$ is the image $\gamma$, $\partial_2 H'$ is the image of $axis(g_3)$ and $\partial_3 H'$ is the image of $\gG_5 \cup \gF_5$. The inclusion map $<g_1,g_3> < \Gamma$ induces a map $j'':H' \to \M$. $\partial_1 H'$ is identified with $\partial_1 H$ as both are identified with $\gamma/g_1$. So we may define $S'=S \cup_{\partial_1 H'=\partial_1 H} H'$ and $j': S' \to \M$ is the map extending both $j$ and $j''$. 

The length of $\partial_3 H'$ (with respect to $j'$) is $2\wG_5=L + \rho_5$. The twist parameter at $\partial_1 H$ is, by definition, equal to $\wH_5 - i\pi$. By the previous theorem this real part bounded by $10 {\hat T} \exp(-L/4)$ and imaginary part bounded by $100\epsilon \exp(-L/4)$. Since $\epsilon$ is arbitrary this concludes the lemma.
\end{proof}

\begin{proof}(of theorem \ref{thm:surface subgroup case})

Let $\M$ be a closed hyperbolic 3-manifold. Suppose that there does not exist any closed totally geodesic surface immersed in $\M$. We identify $\M$ with $\H^3/\Gamma$ for some cocompact discrete group $\Gamma < PSL_2(\C)$. Let $\epsilon>0$. Let $L_0=L_1$ be given by the previous lemma. Let $L > L_0$.

Using the proof of the above lemma, it can be shown that there exists a map $j:H \to \M$ from a boundary-ordered pair of pants into $\M$ such that
\begin{eqnarray*}
\Big| length_j\big(\partial_k H\big)-L\Big| &\le& \epsilon\\
\end{eqnarray*}
for $k=1,2,3$. Applying the lemma successively, we obtain a map $j':S \to \M$ from a surface $S$ into $M$ satisfying the conclusions of the theorem.

\end{proof}

The proof of theorem \ref{thm:ehrenpreis case} involves only notational changes to the above proof so we omit it.

\subsection{Fixed Points of $\g$}


Let $\{e_0,e_1\} \subset \C$ be the fixed points of $\g$. Let
\begin{displaymath}\begin{array}{ll}
N_1 &= \exp(L/2 + \nu + i\theta)\cos(\delta) - \exp(-L/2 + \nu - i\theta)T\sin(\delta) - \exp(-L/2 - \nu-i\theta)\cos(\delta).\\
N_2 &= \Big[\Big(\exp(L/2 + \nu + i\theta)\cos(\delta) + \exp(-L/2
  -i\theta)\big(e^\nu T\sin(\delta) + e^{- \nu}\cos(\delta)\big)\Big)^2 - 4\Big]^{1/2}.\\
D   &=  2\exp(-L/2)[Te^{\nu+i\theta}\cos(\delta) - e^{-\nu + i\theta}\sin(\delta)].
\end{array}
\end{displaymath}

From equation \ref{eqn:fixed points} (subsection
\ref{subsection:fixed points}) and equation \ref{eqn:g} (section
\ref{section:isometry construction}) we obtain 
\begin{equation}\label{eqn:fix(g)}
\{e_0,e_1\} = \left\{ \frac{N_1 - N_2}{D}, \frac{N_1 + N_2}{D} \right\}.
\end{equation}

After relabeling if necessary we may assume $e_0=(N_1-N_2)/D$ and $e_1=(N_1+N_2)/D$. When $L$ is large and the other parameters are
small, $e_0$ is close to zero and $e_1$ is ``close'' to $\infty$. The
next proposition will be useful in estimating the widths of the
hexagon $\hH$ and the pentagon $\hK$.

\begin{pro}\label{pro:fixed points}

The following estimates and identities hold:
\begin{eqnarray*}
N_1 & \approx& \exp(L/2).\\
N_2 & \approx& \exp(L/2).\\
N_1 - N_2 &=& O(\exp(-L/2)).\\
N_1+N_2    &\approx& \exp(L/2).\\
N_1^2 - N_2^2 &=& O(1).\\
D          &=& O(\exp(-L/2)).\\
\frac{N_1^2-N_2^2}{D} &=&-2\sin(\delta)\exp(L/2 + \nu -i\theta) = O(\exp(L/2)).\\
\frac{ \frac{N_1^2-N_2^2}{D} + D\exp(L)}{\exp(L/2)N_2} &=& O(\exp(-L/2)).
\end{eqnarray*}

\end{pro}

\begin{proof}
The estimates for $N_1, N_2, N_1+N_2$ and $D$ are immediate. Recall that for $x$ close to zero, $\sqrt{1-x}=1 -(1/2)x + O(x^2)$. So if $x$ is very large, $\sqrt{x^2 -4} = x\sqrt{1 - 4/x^2} = x - 2/x + O(1/x^3)$. Hence,\begin{eqnarray*}
N_2 &=& \exp(L/2 + \nu+i\theta)\cos(\delta) + O(\exp(-L/2)).
\end{eqnarray*}
Thus $N_1 - N_2 = O(\exp(-L/2))$ as required. We compute $N_1^2-N_2^2$ as follows:
\begin{eqnarray*}
N_1^2-N_2^2 &=& \Big(\exp(L/2 + \nu + i\theta)\cos(\delta) - \exp(-L/2 -i\theta)\big(e^\nu T\sin(\delta) + e^{- \nu}\cos(\delta)\big) \Big)^2\\
            &&- \Big(\exp(L/2 + \nu+i\theta)\cos(\delta) + \exp(-L/2 -i\theta)\big(e^\nu T\sin(\delta) + e^{- \nu}\cos(\delta)\big) \Big)^2 + 4\\
            &=& -4 \exp(\nu)\cos(\delta)\Big( \exp(\nu)T\sin(\delta)  + \exp(- \nu)\cos(\delta)\Big) + 4\\
            &=& 4\sin(\delta)e^\nu\Big(e^{-\nu}\sin(\delta) - e^\nu T\cos(\delta) \Big)= O(1).
\end{eqnarray*}
Therefore,
\begin{eqnarray*}
\frac{N_1^2-N_2^2}{D} &=& \frac{4\sin(\delta)e^\nu [e^{-\nu}\sin(\delta) - e^\nu T \cos(\delta)] }{2\exp(-L/2+i\theta)[Te^\nu\cos(\delta) - e^{-\nu}\sin(\delta)]}\\
                      &=& -2\sin(\delta)\exp(L/2 + \nu-i\theta) = O(\exp(L/2)).
\end{eqnarray*}

Since $D\exp(L)=O(\exp(L/2))$ and $\exp(L/2)N_2 \approx \exp(L)$ this implies that
\begin{eqnarray*}\label{eqn:weird2}
\frac{ \frac{N_1^2-N_2^2}{D} + D\exp(L)}{\exp(L/2)N_2} = O(\exp(-L/2))
\end{eqnarray*}
as required.
\end{proof}

\subsection{The Width $\wH_2$}

In this subsection we prove:

\begin{pro}\label{pro:H2}
The width $\wH_2$ satisfies
\begin{displaymath}\begin{array}{ll}
\cosh(\wH_2) &= -N_1/N_2.\\
\wH_2 &= i\pi + O(\exp(-L/2)).
\end{array}
\end{displaymath} 

\end{pro}

\begin{proof}
Recall the definition of the cross ratio $R$ (subsection \ref{subsection:cross ratio}). Let
\begin{displaymath}
R:=R(\infty,0,e_0,e_1) = \frac{e_1}{e_0} = \frac{N_1+N_2}{N_1-N_2}.
\end{displaymath}
Recall that $\wH_2 = \mu(\gH_1,\gH_3;\gH_2)$. By subsection \ref{subsection:doublecross}, $\tanh^2(\mu(\gH_1,\gH_3;\gH_2)/2)=R$. So
\begin{eqnarray*}
\cosh(\wH_2) = \frac{1+R}{1-R} = \frac{-N_1}{N_2}.
\end{eqnarray*}
This proves the first statement. Note
\begin{displaymath}
N_1/N_2 = 1 + \frac{N_1 - N_2}{N_2} = 1 + O(\exp(-L))
\end{displaymath}
since $N_1 - N_2 = O(\exp(-L/2))$ and $N_2 \approx \exp(L/2)$ by the previous proposition. So $\wH_2 = i\pi + O(\exp(-L/2))$ as required.

\end{proof}

\subsection{The Width $\wH_4$}

\begin{pro}
Let $\{f_0,f_1\}$ be the endpoints of $\gH_5$ with $|f_0| < |f_1|$. Let $\m=\tanh(\gM/2)$. Then
\begin{displaymath}\begin{array}{ll}
f_0 &= -\m \X\\
f_1 &= -\X/\m.
\end{array}
\end{displaymath}

\end{pro}

\begin{proof}

By definition $\gH_5$ is orthogonal to the geodesic with endpoints $\{-\X,\X \}$. Thus 
\begin{displaymath}
R(f_0,f_1,\X,-\X)=\tanh^2(i\pi/4)=-1.
\end{displaymath}
So,
\begin{displaymath}
\frac{(f_0-\X) (f_1 + \X ) }{ (f_0 + \X) (f_1 - \X) } = -1.
\end{displaymath}  
Equivalently,
\begin{displaymath}
f_0f_1 +(f_0- f_1)\X - \X^2 = -f_0f_1 + (f_0 - f_1)\X + \X^2.
\end{displaymath}
So, $f_0f_1 = \X^2$. The width $\wH_6=\mu(\gH_5,\gH_1;\gH_6)$ equals $\gM$. By subsection \ref{subsection:doublecross} $R(f_0,f_1,0,\infty)=\tanh^2(\wH_6/2)$. So,
\begin{displaymath}
\frac{(f_0-0)(f_1-\infty)}{(f_0-\infty)(f_1-0)} = \frac{f_0}{f_1} = \tanh^2(\gM/2).
\end{displaymath}
Thus
\begin{eqnarray*}
f_0^2 &=& f_0f_1 \frac{f_0}{f_1} = \tanh^2(\gM/2)\X^2 = (-\m\X)^2 \textrm{  and}\\
f_1^2 &=& f_0f_1\frac{f_1}{f_0} = \coth^2(\gM/2)\X^2 = (-\X/\m)^2.
\end{eqnarray*}

The choice of sign is justified by figure \ref{fig:hexagonH}.

\end{proof}






\begin{pro}\label{pro:H4}
The width $\wH_4$ satisfies:
\begin{eqnarray*}
\cosh(\wH_4) &=& (N_1/N_2)\cosh(\gM)+ (1/2)\sinh(\gM)Z;\\
\sinh(\wH_4) &=& \sinh(\gM)+ (1/2)\cosh(\gM)Z - (1/4)\cosh(\gM)\coth(\gM)Q + O(\exp(-L));\\
\wH_4 &=& \gM + O(\exp(-L/2)).
\end{eqnarray*}

where 
\begin{eqnarray*}
Z &=& \frac{D\X^2 + (N_1^2 - N_2^2)/D }{\X N_2 } = O(\exp(-L/2)), \\
Q &=& \frac{D^2 \X^4 + (N_1^2 - N_2^2)^2/D^2 }{\X^2 N_2^2} = O(\exp(-L )).
\end{eqnarray*}

\end{pro}

\begin{proof}

We compute $R(f_0,f_1,e_0,e_1)$ as follows.

\begin{eqnarray*}
R(f_0,f_1,e_0,e_1) &=& \frac{(-\m \X -\frac{N_1-N_2}{D}) (-\X/\m -\frac{N_1+N_2}{D}) }{ (-\m \X -\frac{N_1+N_2}{D}) (-\X/\m -\frac{N_1-N_2}{D}) }\\
                   &=&\frac{(-D\m \X - (N_1 - N_2))(-D\X - \m(N_1 + N_2)) }{ (-D \m \X - (N_1 +N_2)  )( -D\X - \m(N_1 - N_2) ) }\\
                   &=&\frac{ D^2\m \X^2 + D\X(N_1 - N_2) + D \m^2 \X(N_1 + N_2) + \m(N_1^2 - N_2^2) }{ D^2\m \X^2 + D\X(N_1 + N_2) + D \m^2 \X(N_1 - N_2) + \m(N_1^2 - N_2^2)}\\
                   &=&\frac{ D\m \X^2 + \X(N_1 - N_2) +  \m^2 \X(N_1 + N_2) + \m(N_1^2 - N_2^2)/D }{ D\m \X^2 + \X(N_1 + N_2) +  \m^2 \X(N_1 - N_2) + \m(N_1^2 - N_2^2)/D} = (a/b)
\end{eqnarray*}
where $a$ is the numerator in the line above and $b$ is the denominator. Since $R = \tanh^2(\wH_4/2)$ we obtain the following.

\begin{eqnarray*}
\cosh(\wH_4) &=&\frac{1 + R}{1-R} = \frac{1 + (a/b)}{1- (a/b)} = \frac{b+a}{b-a}\\
             &=&\frac{2D\m \X^2 + 2(1+\m^2)\X N_1 + 2\m(N_1^2 - N_2^2)/D}{ 2(1-\m^2)\X N_2}\\
             &=&\frac{(1+\m^2)\X N_1}{ (1-\m^2)\X N_2} + \frac{D\m \X^2  + \m(N_1^2 - N_2^2)/D}{ (1-\m^2)\X N_2}\\
             &=&\cosh(\gM)(N_1/N_2) + \frac{\m}{1-\m^2} \left(\frac{D\X^2 + (N_1^2 - N_2^2)/D}{\X N_2}\right)\\
             &=&\cosh(\gM)(N_1/N_2) + (1/2)\sinh(\gM)\left(\frac{D\X^2 + (N_1^2 - N_2^2)/D}{\X N_2}\right)\\
&=& \cosh(\gM)(N_1/N_2) + (1/2)\sinh(\gM)Z.
\end{eqnarray*}
This proves the first statement. Let $X=\exp(L/2)$. Then

\begin{eqnarray*}
Z &=& \frac{D\X^2 + (N_1^2 - N_2^2)/D}{\X N_2}\\
&=& \frac{D(X^2 + O(\exp(3L/4))) + (N_1^2 - N_2^2)/D}{(X+O(\exp(L/4))) N_2}\\
&=& \frac{DX^2 + (N_1^2 - N_2^2)/D }{X N_2 } + O(\exp(-3L/4))\\
&=& O(\exp(-L/2)).
\end{eqnarray*}
The last estimate comes from proposition \ref{pro:fixed points}. Now we estimate $\sinh(\wH_4)$.

\begin{eqnarray*}
&&\sinh(\wH_4) = \left( \cosh^2(\wH_4)-1 \right)^{1/2}\\
             &=& \left( \left( \cosh(\gM)(N_1/N_2) + (1/2)\sinh(\gM)Z \right)^2 - 1\right) ^{1/2}\\
             &=&\left( \cosh^2(\gM)(N_1/N_2)^2 +\cosh(\gM)\sinh(\gM)Z(N_1/N_2)+ (1/4)\sinh^2(\gM)Z^2  - 1\right) ^{1/2}\\
 &=&\left( \sinh^2(\gM) + \cosh^2(\gM)[(N_1/N_2)^2 -1]+\cosh(\gM)\sinh(\gM)Z(N_1/N_2)+ (1/4)\sinh^2(\gM)Z^2  \right) ^{1/2}\\
             &=&\sinh(\gM)\left( 1 + \coth^2(\gM)[(N_1/N_2)^2 -1] +\coth(\gM)Z(N_1/N_2)+ (1/4)Z^2  \right)^{1/2}\\
             &=&\sinh(\gM)\left( 1 + \coth^2(\gM)[(N_1/N_2)^2 -1] +\coth(\gM)Z + O(\exp(-L))  \right)^{1/2}.
\end{eqnarray*}

We check the order of magnitude of the terms above. Since $N_1^2 - N_2^2 = O(1)$ and $N_2^2 \approx \exp(-L)$ the term $[(N_1/N_2)^2 -1]= [N_1^2-N_2^2]/N_2^2=O(\exp(-L))$. Since $\gM \approx \exp(-L/4)$, $\coth(\gM) \approx \exp(L/4)$. So 
\begin{displaymath}
 \coth^2(\gM)[(N_1/N_2)^2 -1] =O(\exp(-L/2)).
\end{displaymath}

Since $Z=O(\exp(-L/2))$ we have that $\coth(\gM)Z =  O(\exp(-L/4))$. So the expression above is equal to $\sinh(\gM)\sqrt{1+x}$ where $x = O(\exp(-L/4))$. Recall that $\sqrt{1 + x} = 1 + (1/2)x -(1/4)x^2 + O(x^3)$. So \begin{eqnarray*}
&&\sinh(\wH_4)\\ 
&=&\sinh(\gM)\big( 1 + (1/2)\coth^2(\gM)[(N_1/N_2)^2 -1] +(1/2)\coth(\gM)Z+ O(\exp(-L))\\
 && -(1/4)\big[\coth^2(\gM)[(N_1/N_2)^2 -1] +\coth(\gM)Z+ O(\exp(-L))\big]^2 + O(\exp(-3L/4)) \big)\\
 &=&\sinh(\gM)\big( 1 + (1/2)\coth^2(\gM)[(N_1/N_2)^2 -1] +(1/2)\coth(\gM)Z+ O(\exp(-L))\\
 && -(1/4)[\coth(\gM)Z+ O(\exp(-L/2))]^2 + O(\exp(-3L/4)) \big)\\
&=&\sinh(\gM) + (1/2)\cosh(\gM)Z + (1/2)\cosh(\gM)\coth(\gM)[(N_1/N_2)^2 -1] \\
&& - (1/4)\cosh(\gM)\coth(\gM)Z^2 + O(\exp(-L))\\
&=&\sinh(\gM) + (1/2)\cosh(\gM)Z\\
&&+ (1/2)\cosh(\gM)\coth(\gM)\left([(N_1/N_2)^2 -1] - (1/2)Z^2\right) + O(\exp(-L)).
\end{eqnarray*}


We estimate the last coefficient:
\begin{eqnarray*}
[(N_1/N_2)^2 -1] - (1/2)Z^2 &=&\frac{N_1^2 - N_2^2}{N_2^2} - (1/2)\left(\frac{D\X^2 + (N_1^2 - N_2^2)/D}{\X N_2}\right)^2\\
&=&\frac{N_1^2 - N_2^2}{N_2^2} - (1/2)\frac{D^2\X^4 + (N_1^2 - N_2^2)^2/D^2 + 2\X^2(N_1^2-N_2^2)}{\X^2 N_2^2}\\
&=& - (1/2)\frac{D^2\X^4 + (N_1^2 - N_2^2)^2/D^2 }{\X^2 N_2^2}= -(1/2)Q.
\end{eqnarray*}

$D^2 \X^4 = O(\exp(L))$ since $D=O(\exp(-L/2))$ and $\X \approx \exp(L/2)$. The previous proposition implies $(N_1^2 - N_2^2)^2/D^2=O(\exp(L))$. So the numerator is on the order of $\exp(L)$. Since $N_2 \approx \exp(L)$ the denominator $\X^2 N_2^2 \approx \exp(2L)$. Thus $Q = O(\exp(-L))$. Thus we have:
\begin{displaymath}
\sinh(\wH_4) = \sinh(\gM)+ (1/2)\cosh(\gM)Z - (1/4)\cosh(\gM)\coth(\gM)Q + O(\exp(-L))
\end{displaymath}
as required. Next we compute $\sinh(\wH_4 - \gM)$ as follows.
\begin{eqnarray*}
\sinh(\wH_4 -\gM) &=& \sinh(\wH_4)\cosh(\gM) - \cosh(\wH_4)\sinh(\gM)\\
&=&\cosh(\gM)\sinh(\gM)+ (1/2)\cosh^2(\gM)Z - (1/4)\cosh^2(\gM)\coth(\gM)Q\\
&& + O(\exp(-L)) -\cosh(\gM)\sinh(\gM) - (1/2)\sinh^2(\gM)Z + O(\exp(-5L/4))\\
&=&(1/2)Z - (1/4)\cosh^2(\gM)\coth(\gM)Q + O(\exp(-L))\\
&=&O(\exp(-L/2)).
\end{eqnarray*}
\end{proof}

\subsection{The Width $\wH_5$}

\begin{pro}\label{pro:H5}
Assume that $|\tan(\delta)| \le 2\epsilon$ and $|e^{-2\nu}| \le
2$. Then we have the following estimates:
\begin{displaymath}
\wH_5 =(\sqrt{2}/2)\coth(\gM)\exp(-L/2)(T +\tau  ) +i\pi  + O(\exp(-L/2))
\end{displaymath}
where $\tau \in \C$ is such that $|\tau| \le 6\epsilon$.
\end{pro}

\begin{proof}
Recall that $\wH_6=\gM + i\pi$. By the law of cosines we obtain
\begin{eqnarray*}
\cosh(\wH_5) &=& \frac{\cosh(\wH_2) - \cosh(\wH_6)\cosh(\wH_4)} {\sinh(\wH_6)\sinh(\wH_4)}\\
             &=& \frac{-N_1/N_2 +   \cosh(\gM)[(N_1/N_2)\cosh(\gM)+ (1/2)\sinh(\gM)Z]}{-\sinh(\gM)[\sinh(\gM)+ (1/2)\cosh(\gM)Z - (1/4)\cosh(\gM)\coth(\gM)Q + O(\exp(-L)) ]}\\
           &=& \frac{(N_1/N_2)\sinh^2(\gM)+ (1/2)\cosh(\gM)\sinh(\gM)Z }{-\sinh(\gM)[\sinh(\gM)+ (1/2)\cosh(\gM)Z- (1/4)\cosh(\gM)\coth(\gM)Q + O(\exp(-L)) ]}\\
&=&  - \frac{ \sinh(\gM)(N_1/N_2 ) +(1/2)\cosh(\gM)Z  }{\sinh(\gM)+ (1/2)\cosh(\gM)Z - (1/4)\cosh(\gM)\coth(\gM)Q} + O(\exp(-3L/4))\\ 
&=&  -1 - \frac{\sinh(\gM)(N_1/N_2 - 1)+ (1/4)\cosh(\gM)\coth(\gM)Q }{\sinh(\gM)+ (1/2)\cosh(\gM)Z - (1/4)\cosh(\gM)\coth(\gM)Q} + O(\exp(-3L/4))\\ 
&=&  -1 - \frac{(1/4)\cosh(\gM)\coth(\gM)Q }{\sinh(\gM)} + O(\exp(-3L/4))\\ 
&=&  -1 - (1/4)\coth^2(\gM)Q + O(\exp(-3L/4)).
\end{eqnarray*}
The issue now is to compute $\sqrt{Q}$. Recall that
\begin{eqnarray*}
Q &=& \frac{D^2 \X^4 + (N_1^2 - N_2^2)^2/D^2 }{\X^2 N_2^2}\\
(N_1^2 - N_2^2)^2/D^2 &=& [-2\sin(\delta)\exp(L/2 + \nu - i\theta)]^2\\
D &=&2\exp(-L/2 + i\theta)[Te^\nu\cos(\delta) - e^{-\nu}\sin(\delta)]\\
N_2^2 &=& \Big(\exp(L/2 + \nu + i\theta)\cos(\delta) + \exp(-L/2
  -i\theta)\big(e^\nu T\sin(\delta) + e^{- \nu}\cos(\delta)\big)\Big)^2 - 4.
\end{eqnarray*}
Next we estimate $Q$. Let $X=\exp(L/2)$. 
\begin{eqnarray*}
Q &=& \frac{D^2 \X^4 + (N_1^2 - N_2^2)^2/D^2 }{\X^2 N_2^2}\\
  &=& \frac{D^2 (X +O(\exp(L/4)))^4 + (N_1^2 - N_2^2)^2/D^2 }{(X + O(\exp(L/4)))^2 N_2^2}\\
&=& \frac{D^2 X^4 + (N_1^2 - N_2^2)^2/D^2 }{X^2 N_2^2 } + O(\exp(-5L/4)).
\end{eqnarray*}

So,

\begin{eqnarray*}
Q &=& X^{-2}\frac{4\exp(L+2i\theta)[Te^\nu\cos(\delta) - e^{-\nu}\sin(\delta)]^2 + [-2\sin(\delta)\exp(L/2 + \nu -i\theta)]^2}{\big[\exp(L/2 + \nu + i\theta)\cos(\delta) + \exp(-L/2 -i\theta)\big(e^\nu T\sin(\delta) + e^{- \nu}\cos(\delta)\big)\big]^2 - 4}\\
&& + O(\exp(-5L/4))\\
&=&\frac{4e^{2i\theta}[Te^\nu\cos(\delta) - e^{-\nu}\sin(\delta)]^2 + 4\sin^2(\delta)e^{2\nu-2i\theta} }{\big[\exp(L/2 + \nu + i\theta)\cos(\delta) + \exp(-L/2 -i\theta)\big(e^\nu T\sin(\delta) + e^{- \nu}\cos(\delta)\big) \big]^2 - 4} + O(\exp(-5L/4)).
\end{eqnarray*}

We only need to know $Q$ up to $O(\exp(-L))$. So we simplify the denominator as follows.
\begin{displaymath}
Q = \frac{4e^{2i\theta}[Te^\nu\cos(\delta) - e^{-\nu}\sin(\delta)]^2 + 4\sin^2(\delta)e^{2\nu-2i\theta} }{\exp(L + 2\nu+2i\theta)\cos^2(\delta)}+ O(\exp(-5L/4)).
\end{displaymath}
The numerator equals
\begin{eqnarray*}
&&4e^{2i\theta}[Te^\nu\cos(\delta) - e^{-\nu}\sin(\delta)]^2 + 4\sin^2(\delta)e^{2\nu-2i\theta}\\
&=&4T^2 e^{2\nu+2i\theta}\cos^2(\delta) - 8Te^{2i\theta}\cos(\delta)\sin(\delta) + 4\sin^2(\delta)(e^{2\nu-2i\theta}+e^{-2\nu+2i\theta}).
\end{eqnarray*}
So
\begin{eqnarray*}
Q &=& 4\exp(-L)\Big(T^2 - 2T\tan(\delta)e^{-2\nu} + \tan^2(\delta)(e^{-4i\theta}+e^{-4\nu})\Big) + O(\exp(-5L/4))\\
  &=& 4\exp(-L)\Big( (T - \tan(\delta)e^{-2\nu})^2 + \tan^2(\delta)e^{-4i\theta} \Big) + O(\exp(-5L/4)).
\end{eqnarray*}

We need to estimate $\sqrt{Q}$. Notice that there is a choice of a square root for
\begin{displaymath}
(T-\tan(\delta)e^{-2\nu})^2+\tan^2(\delta)e^{-4i\theta}
\end{displaymath}
such that
\begin{displaymath}
\Big| T- \sqrt{ (T - \tan(\delta)e^{-2\nu})^2 + \tan^2(\delta)e^{-4i\theta} }\Big| \le  \Big|\tan(\delta)(|e^{-2\nu}|+1)\Big| \le 6\epsilon.
\end{displaymath}
Above we used the hypotheses $|\tan(\delta)| \le 2\epsilon$ and $|e^{-2\nu}| \le 2$. Hence there exists a number $\tau \in \C$ with $|\tau| \le 6\epsilon$ such that
\begin{displaymath}
(T + \tau)^2 = (T - \tan(\delta)e^{-2\nu})^2 + \tan^2(\delta)e^{-4i\theta}.
\end{displaymath}
So we obtain a square root of $Q$ as follows:
\begin{displaymath}
\sqrt{Q} = \exp(-L/2)(T + \tau) +O(\exp(-3L/4)).
\end{displaymath}

Recall that $\cosh(x)= 1 + x^2/2 + O(x^4)$. Hence 

\begin{eqnarray*}
\wH_5 &=&(\sqrt{2}/2)\coth(\gM)\sqrt{Q} + i\pi + O(\exp(-L/2))\\
      &=&(\sqrt{2}/2)\coth(\gM)\exp(-L/2)(T +\tau  ) +i\pi  + O(\exp(-L/2)).\\
\end{eqnarray*}
The choice of square root is justified by figure \ref{fig:hexagonH} (which is drawn in the case that $T>0$).
\end{proof}


\part{Tree Tilings}\label{part:tree tilings}

In this part, we develop a formalism to describe tilings of $Tree$, the Cayley graph of $F=\Z/2\Z * \Z/2\Z *\Z/2\Z$. We use this to convert the problem of showing the existence of a closed surface with a labeled pants decomposition of the kind required by questions \ref{question:ehrenpreis case} and \ref{question:surface subgroup case} into the problem of showing the existence of a periodic tiling of $Tree$. Then we convert that problem into a linear programming problem.

\section{Definitions}

Let $F=\Z/2\Z *\Z/2\Z*\Z/2\Z=<a,b,c |, a^2=b^2=c^2=1>$. Let $Tree$ denote the labeled graph with vertex set $F$ such that for every $w\in\F$ there exists an edge labeled $a$ from $w$ to $wa$, an edge labeled $b$ from $w$ to $wb$ and an edge labeled $c$ from $w$ to $wc$. These are all of the edges.

We say that $G$ is a {\bf tileset graph} for $Tree$ if $G$ is a finite graph such that each edge is labeled $a$, $b$ or $c$. A {\bf tiling} of $Tree$ by a tileset graph $G$ is a map $\phi: Tree \to G$ that sends vertices to vertices, edges to edges, preserves incidence and labels.

 $\Free$ acts on the set of vertices of $\Tree$ by group
 multiplication on the left. This action extends to the edges in the
 obvious way so that labels and directions are preserved. $F$ also acts on the set of tilings $\phi:F \to G$ by ``moving the tiles around''. To be precise:
\begin{displaymath}
(g\phi)(f) = \phi(g^{-1}f)
\end{displaymath}
where $f$ is either an edge or a vertex of $\Tree$ and $g \in F$. We
say that a tiling is {\bf periodic} if its stabilizer has
finite index in $F$. Equivalently, its $F$ orbit is finite.

\subsection{$Y$-Graphs: the $3$-manifold case}\label{subsection:y-graphs 3-manifold case}

Let $\M$ be a closed hyperbolic $3$-manifold. We will construct a tileset graph so that tilings correspond to immersions $j:S \to \M$ of surfaces into $\M$ satisfying certain geometric constraints. 

Suppose for $i=1,2$, $j_i:H_i \to \M$ is a map from a boundary-ordered pair of pants to $\M$ and there is a map $\Phi: H_1 \to H_2$ that preserves the boundary order such that $\Phi \circ j_2$ is homotopic to $j_1$. Then we say that $j_1$ is homotopic to $j_2$. 

Let $V=V(\M,L,\epsilon)$ denote the set of homotopy classes of maps $v=(j:H \to \M)$ satisfying:
\begin{itemize}
\item $H$ is a boundary-ordered pair of pants,
\item $|length_j(\partial_k H)-L| < \epsilon$ for $k=1,2,3$.
\end{itemize}
$V$ is a finite set because there is only a finite number of geodesics in $\M$ of length no greater than $L+\epsilon$.

 Let $d \in \{a,b,c\}$ and let ${\bar d}$ equal 1 if $d=a$, equal 2 if $d=b$ and equal $3$ if $d=c$. Define a graph $Y=Y(\M,L,\epsilon)$ with vertex $V$ as follows. Roughly speaking, there exists an $d$-labeled edge in $Y$ from $(j_1: H_1 \to \M)$ to $(j_2: H_2 \to \M)$ iff we can glue $H_1$ to $H_2$ along $\partial_{\bar d}(H_1)$ and $\partial_{\bar d}(H_2)$ to obtain a map $b: H_1 \cup_{\partial_d H_1} H_2 \to \M$ extending $j_1$ and $j_2$. Precisely, there is a $d$-labeled edge between $v$ and $v'$ iff there exists a map $e=(b_e: B_e \to \M)$ such that:
\begin{itemize}
\item $B_e$ is a four-holed sphere,
\item $B_e$ has a labeled pants decomposition $\P=\{H,H'\}$,
\item $v=(b_e|_H: H \to \M)$ and $v'=(b_e|_{H'}: H' \to \M)$,
\item the unique simple closed curve in $\P^*$ in the interior of $B$ equals $\partial_{\bar d}(H)=\partial_{\bar d}(H')$,
\item $\big| \Im\big(twist(\partial_{\bar d} H)\big)\big| \le \epsilon \exp(-L/4)$.
\end{itemize}

\begin{lem}\label{lem:single surface periodic tiling implies closed surface}
Suppose there exists a periodic tiling $\phi:Tree \to Y=Y(\M,L,\epsilon)$. Then the conclusions to theorem \ref{thm:surface subgroup case} can be strenghened so that the surface $S$ is closed.
\end{lem}

The proof is similar to a standard construction in graphs of groups theory.

\begin{proof}

For each vertex $v \in V$ choose a representative $j:H \to \M$ of $v$ such that $j$ is locally $1-1$ on the boundary and maps each boundary curve to a geodesic. For each edge $e$ in $Y$, choose a representative $b_e:B_e \to \M$ of $e$ so that the following holds. If the endpoints of $e$, $v_1,v_2$ have chosen representatives $j_1:H_1 \to \M$, $j_2:H_2 \to \M$ then there is a map $i_e:H_1 \cup H_2 \to B_e$ from the disjoint union of $H_1$ and $H_2$ whose restriction to either $H_1$ or $H_2$ is inclusion and such that $b_e \circ i_e$ is equal to $j_i$ when restricted to $H_i$ (for $i=1,2$).

For $f \in F=\Z/2\Z*\Z/2\Z*\Z/2\Z$, let $j_f: H_f \to \M$ be a copy of the chosen representative of $\phi(f) \in V$. Define an equivalence relation $\sim$ on the disjoint union ${\hat S} = \bigcup_{f \in F} \, H_f$ as follows. Suppose $x \in H_f$ and $y \in H_{fd}$ (for some $d\in\{a,b,c\}$). Let $e=\phi(\{f,fd\})$ where $\{f,fd\}$ denotes the edge between $f$ and $fd$ in $Tree$. Let $x \sim y$ if $i_e(x)=i_e(y)$ where $i_e:H_f \cup H_{fd} \to B_e$ is as above. Let $\sim$ be the smallest equivalence relation on $\bigcup_{f \in F} \,  H_f$ satisfying the above.

Let $S_\phi$ be the surface defined by
\begin{displaymath}
S_\phi = \big( \bigcup_{f \in F} \,  H_f \big)/\sim
\end{displaymath}

The images of the boundary-ordered pants $H_f$ comprise a labeled pants decomposition $\P_\phi$ for $S_\phi$.
Define $j_\phi:S_\phi \to \M$ by $j_\phi([x])=j_f(x)$ where $[x]$ denotes the equivalence class of $x \in H_f$. $j_\phi$ is well-defined by construction. By construction, for every curve $\gamma \in \P_\phi^*$,
\begin{eqnarray*}
|length_{j_\phi}(\gamma)-L| &\le& \epsilon, \\
|\Im(twist_{j_\phi}(\gamma))| &\le& \epsilon\exp(-L/4).
\end{eqnarray*}

 Let $\Gamma$ be the stabilizer of $\phi$ in $F$. Since $\phi$ is periodic, $\Gamma$ has finite index in $F$. $\Gamma$ acts on the disjoint union $\bigcup_{f \in F} \,  H_f$ in the obvious way. This action descends to an action on $S_\phi$ since $\Gamma$ preserves equivalence classes. This action preserves the labeled pants decomposition $\P$. So the quotient space $S=S_\phi/\Gamma$ is a surface with a labeled pants decomposition $\P=\P_\phi/\Gamma$. The map $j_\phi: S_\phi \to \M$ is preserved under the action of $\Gamma$ (i.e. $j_\phi(\gamma x)=j_\phi(x)$ for all $x\in S_\phi$ and $\gamma \in \Gamma$). Hence it descends to a map $j: S \to \M$. By construction all curves $\gamma$ in $\P^*$ satisfy the above bounds. $S$ is compact (and thus closed since it has no boundary) since $\Gamma$ has finite index in $F$. 
 
\end{proof}

\subsection{$Y$-Graphs: the product of two surfaces case}\label{subsection:y-graphs}

This section is similar to the previous one; for a given pair of closed hyperbolic surfaces $S_1,S_2$ and positive numbers $L,\epsilon$ we construct a tileset graph $Y(S_1 \times S_2, L, \epsilon)$ so that periodic tilings of $\Tree$ by $Y$ correspond to pairs of finite covers $\pi:{\tilde S_i} \to S_i$ satisfying the geometric constraints of theorem \ref{thm:ehrenpreis case}.

Let $\M=S_1\times S_2$. Suppose for $i=1,2$, $j_i:H_i \to \M$ is a map from a boundary-ordered pair of pants to $\M$ and there is a map $\Phi: H_1 \to H_2$ that preserves the boundary order such that $\Phi \circ j_2$ is homotopic to $j_1$. Then we say that $j_1$ is homotopic to $j_2$.

Let $V=V(\M,L,\epsilon)$ denote the set of homotopy classes of maps $v=(j:H \to \M)$ satisfying:
\begin{itemize}
\item $H$ is a boundary-ordered pair of pants,
\item $||length_j(\partial_k H)-(L,L)||_\infty < \epsilon$ for $k=1,2,3$.
\end{itemize}
$V$ is a finite set because there is only a finite number of geodesics in $S_1$ or $S_2$ of length no greater than $L+\epsilon$.

 Let $d \in \{a,b,c\}$ and let ${\bar d}$ equal 1 if $d=a$, equal 2 if $d=b$ and equal $3$ if $d=c$. Define a graph $Y=Y(\M,L,\epsilon)$ with vertex $V$ as follows. Roughly speaking, there exists an $d$-labeled edge in $Y$ from $(j_1: H_1 \to \M)$ to $(j_2: H_2 \to \M)$ iff we can glue $H_1$ to $H_2$ along $\partial_{\bar d}(H_1)$ and $\partial_{\bar d}(H_2)$ to obtain a map $b: H_1 \cup_{\partial_d H_1} H_2 \to \M$ extending $j_1$ and $j_2$. Precisely, there is a $d$-labeled edge between $v$ and $v'$ iff there exists a map $e=(b_e: B_e \to \M)$ such that:
\begin{itemize}
\item $B_e$ is a four-holed sphere,
\item $B_e$ has a labeled pants decomposition $\P=\{H,H'\}$,
\item $v=(b_e|_H: H \to \M)$ and $v'=(b_e|_{H'}: H' \to \M)$,
\item the unique simple closed curve in $\P^*$ in the interior of $B$ equals $\partial_{\bar d}(H)=\partial_{\bar d}(H')$,
\item if $\pi_i:\R^2 \to \R$ denotes projection onto the $i$-th factor then
\begin{eqnarray*}
\Big| \pi_1\big( twist_j(\partial_{\bar d} H)\big) -\pi_2\big( twist_j(\partial_{\bar d} H)\big) \Big| \le \epsilon \exp(-L/4).
\end{eqnarray*}
\end{itemize}

\begin{lem}\label{lem:single surface periodic tiling implies closed surface 2}
Suppose there exists a periodic tiling $\phi:Tree \to Y=Y(\Gamma,L,\epsilon)$. Then the conclusions to theorem \ref{thm:ehrenpreis case} can be strenghened so as to require the covers $\pi_i:{\tilde S_i}\to S_i$ to be finite-sheeted.
\end{lem}
The proof is essentially the same as the proof of lemma \ref{lem:single surface periodic tiling implies closed surface} in the previous section so we omit it.

\section{Periodic Tilings}\label{section:periodic tilings}

In this section, we show that the existence of a periodic tiling $\phi: \Tree \to Y$ is equivalent to the existence of a ``positive flow'' on $Y$ in the following sense. Let $V$ and $E$ denote the vertex set and edge set of $Y$. For a vertex $v \in V$ let $A(v), B(v), C(v) \subset E$ denote the set of edges labeled $a$, $b$, $c$ respectively that are incident to $v$. By a {\bf flow} on $G$ we mean a function $f:\{V \cup E\} \to \R$ such that
\begin{eqnarray*}
f(v) = \Sigma_{e \in A(v)} \, f(e) = \Sigma_{e \in B(v)} \, f(e) = \Sigma_{e \in C(v)} \, f(e).
\end{eqnarray*}
A flow is said to be {\bf positive} if all of its entries are non-negative but not all are zero.

\begin{lem}
There exists a positive flow $Z$ on $Y$ iff there exists a periodic tiling $\phi:\Tree \to Y$.
\end{lem}

\begin{proof}

Suppose that $\phi:\Tree \to Y$ is a periodic tiling. Let $\Gamma<F$ denote
the stabilizer of $\phi$. Let $\Tree/\Gamma$ denote the quotient
graph. This is the graph whose
vertex set is equal to the set of cosets $\Gamma \backslash \Free$ such
that there is a edge from $\Gamma g$ to $\Gamma h$ labeled $d$ if
$\Gamma gd = \Gamma h$ (for $d \in \{a,b,c\}$). Since $\Gamma$ has finite index, $\Tree/\Gamma$ is a finite graph.

The tiling $\phi: \Tree \to Y$ projects to a tiling on the quotient
\begin{displaymath}
{\bar \phi}: \Tree/\Gamma \to Y.
\end{displaymath}

Let $Z_\phi \in \R^{V \cup E}$ be the vector defined by
\begin{displaymath}\begin{array}{ll}
Z_\phi(v) &= |{\bar \phi}^{-1}(v)|\\
Z_\phi(e) &= |{\bar \phi}^{-1}(e)|
\end{array}
\end{displaymath}
for any $v \in V$ or $e \in E$.

Claim: $Z_\phi$ is a positive flow on $Y$. Suppose $v \in V$. Then $Z_\phi(v)$, the number of vertices of $\Tree/\Gamma$ which map to $v$ is equal to the number of $a$-labeled edges of $\Tree/\Gamma$ which are adjacent to such vertices. Since each such edge must map into $A(v)$, the number of such edges equals $\Sigma_{e \in A(v)} \, Z_\phi(e)$. Similar statements hold for $b$ and $c$. Thus $Z_\phi$ is a positive flow.

To prove the other direction, assume there exists a positive flow on $Y$. The flow conditions are integral linear conditions. Specifically, there is an integer matrix $A$ such that $Z \in \R^{V \cup E}$ is a flow iff $AZ=0$. So without loss of generality we can assume that $Z$ is a positive integral flow. We construct a periodic tiling $\phi$ by first constructing a tiling ${\bar \phi}:\Tree/\Gamma \to Y$ from a finite quotient of $\Tree$ to $Y$ and then pulling it back to $\Tree$.

For each vertex $v \in V$, let $X_v$ be a finite set such that
\begin{displaymath}
|X_v| = Z(v).
\end{displaymath}
For each vertex $v \in V$ and edge $e$ incident to $v$ let $X_{v,e} \subset X_v$ be such that
\begin{itemize}
\item $|X_{v,e}| = Z(e)$ and
\item $X_{v,e} \cap X_{v,e'} = \emptyset$ whenever $e \ne e'$ and $e$ and $e'$ have the same label.
\end{itemize}
Because $Z$ is a flow, it is possible to find sets $X_{v,e}$ satisfying the above. 

For each edge $e$ of $Y$ with endpoints $v_1, v_2 \in V$, let $\beta_e:X_{v_1,e} \to X_{v_2,e}$ be a bijection. Let $Q$ be the labeled graph with vertex set equal to the disjoint union $\cup_v X_v$ such that there is a $d$-labeled edge from $x_1 \in X_{v_1,e}$ to $x_2 \in X_{v_2,e}$ iff $\beta_e(x_1)=x_2$ and $e$ is labeled $d$ (for $d \in \{a,b,c\}$). 

By the second property above, if $x$ is a vertex of $Q$ then $x$ is contained in exactly three sets of the form $X_{v,e}$; one for each label $\{a,b,c\}$. Therefore, $Q$ is consistently labeled in the sense that for each connected component $Q'$ of $Q$, the universal covering space of $Q'$ is equal to $\Tree$ (labels included). 

Define ${\bar \phi}:Q \to Y$ by ${\bar \phi}(x) = v$ if $x \in X_v$. If $y$ is an edge from $x_1 \in X_{v_1,e}$ to $x_2 \in X_{v_2,e}$ so that $\beta_e(x_1)=x_2$ then define ${\bar \phi}(y)=e$. This makes ${\bar \phi}$ a graph homomorphism that preserves labels. 

Let $Q'$ be a connected component of $Q$. Let $\pi: \Tree \to Q'$ be a universal covering map ($\pi$ is unique up to the choice of $\pi(id)$). ${\bar \phi}$ pulls back under $\pi$ to a tiling $\phi:\Tree \to Y$. If $\Gamma$ denotes the stabilizer of $\phi$ then $\Tree/\Gamma=Q'$ is a finite graph. Hence $\Gamma$ has finite index in $F$ and thus $\phi$ is a periodic tiling.

\end{proof}

The general problem of determining whether or not there exists a positive flow on a given tileset graph is algorithmically decidable; it is a linear programming problem. From the lemma above and lemmas \ref{lem:single surface periodic tiling implies closed surface} and \ref{lem:single surface periodic tiling implies closed surface 2}, it follows that for fixed $\M,L,\epsilon$, questions \ref{question:ehrenpreis case} and \ref{question:surface subgroup case} are algorithmically decidable. We intend to study the graphs $Y(\M,L,\epsilon)$ in more detail in future work.

\part{Bi-Lipschitz maps}\label{part:bilipschitz}

The goal of this part is to prove theorem \ref{thm:bilipschitz}.

\section{Definitions and Main Theorems}

If $X, Y$ are metric spaces and $f: X \to Y$ is a homeomorphism and $k \ge 1$ then
$f$ is said to be {\bf $k$ bi-Lipschitz} if for every $x_1,x_2 \in X$,
\begin{displaymath}
k^{-1}d(f(x_1),f(x_2)) \le d(x_1,x_2) \le k d(f(x_1),f(x_2))
\end{displaymath}
and for every $y_1,y_2 \in Y$
\begin{displaymath}
k^{-1}d(f^{-1}(y_1),f^{-1}(y_2)) \le d(y_1,y_2) \le k
d(f^{-1}(y_1),f^{-1}(y_2)).
\end{displaymath}

We say that a map $f: X \to Y$ is a {\bf similarity} if there exists a $k > 0$ such that for all $x_1,x_2 \in X$,
\begin{displaymath}
d(f(x_1),f(x_2)) = kd(x_1,x_2).
\end{displaymath}

If $P$ is a hyperbolic three-holed sphere with geodesic boundary then
a point $p \in P$ is called a {\bf special point} if
\begin{itemize}
\item $p$ is on a boundary component $c_1$ of $P$ and
\item there exists a different boundary component $c_2$ of $P$ such
  that the shortest path from $c_1$ to $c_2$ starts at $p$.
\end{itemize}

In this part we prove the following main theorems.

\begin{thm}\label{thm:bilipschitz pants}
Let $\epsilon>0$ be given. Then there exists $\epsilon_1, L_1>0$ such
that the following holds. Suppose $L>L_1$ and $P$ is a hyperbolic 3-holed sphere with geodesic
boundary such that the length of every boundary component of $P$ is in
the interval $(L-\epsilon_1,L+\epsilon_1)$. Let $P_L$ be the hyperbolic
3-holed sphere such that every boundary component of $P_L$ has length
$L$. Then there exists an orientation-preserving homeomorphism $F:P \to P_L$ such that
\begin{itemize}
\item $F$ is $(1+\epsilon)$ bi-Lipschitz,
\item $F$ takes special points to special points and
\item the restriction of $F$ to any boundary component is a similarity onto its image.
\end{itemize}

\end{thm}
We prove this theorem at the end of section
\ref{section:trirectangles}. For $s \in \R$ let
\begin{displaymath}
A_s = \left[ \begin{array}{cc}
e^{s/2} & 0 \\
0 & e^{-s/2} 
\end{array} \right] \in PSL_2(\R)
\end{displaymath}
Consider the annulus $\A=\H^2/<A_l>$ where $l>0$ and $<A_l>$ denotes the discrete group isomorphic to $\Z$ generated by $A_l$. Note that for any $s \in \R$ the action of $A_s$ on $\H^2$ descends to an isometric action on the annulus $\A$. Let ${\bar \gamma}$ denote the projection of the axis of $A_l$ to $\A$. Assume that the axis of $A_l$ is oriented from $0$ to $\infty$ and give ${\bar \gamma}$ the induced orientation.

\begin{thm}\label{thm:bilipschitz annuli} 
Let $\epsilon>0$.  Then there exist positive numbers $E=E(\epsilon)>0$ and $L_0(\epsilon)$ such that if $L>L_0(\epsilon)$, $0\le t \le E\exp(-L/4)$ and $w\ge (1/2)\exp(-L/4)$ then there exists a homeomorphism $F: \A \to \A$ such that the following hold.
\begin{itemize}
\item If $x \in \A$ is at least a distance $w$ from the central geodesic ${\bar \gamma}$ and $x$ is on the left side of ${\bar \gamma}$ then $F(x)=x$.
\item If $x \in \A$ is at least a distance $w$ from the central geodesic ${\bar \gamma}$ and $x$ is on the right side of ${\bar \gamma}$ then $F(x)=A_t(x)$.
\item $F$ is $(1+\epsilon)$ bi-Lipschitz. 
\end{itemize}
\end{thm}
We prove this theorem in section \ref{section:annuli}. Given the above theorems we now prove theorem \ref{thm:bilipschitz}.

\begin{proof}(of theorem \ref{thm:bilipschitz})
Let $\epsilon>0$. Let $\epsilon_2,L_2 >0$ be constants satisfying the following.
\begin{itemize}

\item The conclusion of
theorem \ref{thm:bilipschitz pants} is satisfied when
$\epsilon_1=\epsilon_2$ and $L_1=L_2$. 

\item If $L>L_2$ and $P$ is a pair of pants with boundary lengths in
  the interval $(L-\epsilon_2, L+ \epsilon_2)$ and ${\tilde M}$ is the
  minimum distance between two boundary components of $P$ then ${\tilde M} > 2w:=\exp(-L/4)$.

\item $3\epsilon_2 \le E(\epsilon)$ and $L > L_0(\epsilon)$ when
  $E(\cdot)$ and $L_0(\cdot)$ are the functions defined in theorem \ref{thm:bilipschitz annuli}.

\item $L_2 \ge 1$.

\end{itemize}

To see that the second condition above is attainable, recall that by lemma \ref{lem:L/2 hexagon} we have
\begin{displaymath}
{\tilde M} =  2\exp(-L/4 + \rho_1/4 - \rho_3/4 - \rho_5/4)  + O(\exp(-3L/4))
\end{displaymath}
for some $\rho_1, \rho_3, \rho_5 \in (-\epsilon_2,\epsilon_2)$.

Let ${\tilde S_1}, {\tilde S_2}$ be a pair of surfaces satisfying the
conclusion to theorem \ref{thm:ehrenpreis case} where the constant $\epsilon$ there is replaced with $\epsilon_2$ and $L>L_2$. Recall this means there that for $i=1,2$ there exists a labeled pants decomposition $\P_i$ of ${\tilde  S_i}$ and a homeomorphism $h: {\tilde S_1} \to {\tilde S_2}$ such that $h(\P_1)=\P_2$, $h(\P_1^*)=\P_2^*$ and for every curve $\gamma \in \P_1^*$,

\begin{eqnarray*}
|length(\gamma) - L| &\le & \epsilon_2\\
|length(h(\gamma)-L| &\le & \epsilon_2\\
|twist(\gamma) - twist(h(\gamma))| & \le & \epsilon_2\exp(-L/4).
\end{eqnarray*}
For $k=1,2$, let ${\tilde S'}_k$ be the surface obtained from ${\tilde
  S_k}$ by deforming every curve in $\P_k^*$ to have length $L$. To be
  precise, ${\tilde S'}_k$ is determined up to isometry by the
  following requirements. ${\tilde S'}_k$ has a labeled pants decomposition $\P'_k$ and there is a homeomorphism $h_k:
  {\tilde S}_k \to {\tilde S'}_k$ such that $h_k(\P_k)=\P'_k$,
  $h_k(\P^*_k)=\P'^*_k$ and for every curve $\gamma \in \P'^*_k$,
\begin{eqnarray*}
length(\gamma)&=&L\\
\frac{twist(\gamma)}{L} &=& \frac{twist(h^{-1}(\gamma))}{length(h^{-1}(\gamma))}.
\end{eqnarray*}
The maps $h_1$, $h_2$ can be chosen to be $(1+\epsilon)$ bi-Lipschitz by defining each map $h_k$
on the closure of each component of ${\tilde S}'_k - \P'^*_k$ so that
it satisfies the conclusion of theorem \ref{thm:bilipschitz pants}.

The map $\Phi=h_2 \circ  h \circ h_1^{-1}: {\tilde S}'_1 \to {\tilde
  S}'_2$ is such that for all $\gamma \in \P'^*_1$,
\begin{eqnarray*}
|twist(\gamma)-twist(\Phi(\gamma))| &\le& 3\epsilon_2\exp(-L/4).
\end{eqnarray*}
This can be seen using the above formula for $twist(\gamma)/L$. $\Phi$ can be chosen to be
$(1+\epsilon)$ bi-Lipschitz by defining $\Phi$ on the
$\exp(-L/4)$-neighborhood of each curve $\gamma \in \P'^*_1$ so that
it satisfies the conclusion of theorem \ref{thm:bilipschitz
  annuli} with $t=twist(\Phi(\gamma)) - twist(\gamma)$. This is
well-defined because the $\exp(-L/4)$-neighborhoods of the curves of
the curves in $\P'^*$ are pairwise disjoint. $\Phi$ can be chosen to
be an isometry on the complement of those neighborhoods. 

Now the map $h_2^{-1} \circ \Phi \circ h_1: {\tilde S}_1 \to {\tilde
  S}_2$ is $(1+\epsilon)^3$ bi-Lipschitz. Since $\epsilon$ is
  arbitrary, this proves the theorem.

\end{proof}

\section{Pairs of Pants and Trirectangles}\label{section:trirectangles}

Let $E,L>0$ be given. For $k=0,1$ let $Q_k$ be a convex 4-gon in $\H^2$ with vertices $w_k,x_k,y_k,z_k$. Assume that the interiors angles are all equal to $\pi/2$ except for the angle at $z_k$. Assume that $w_k$ is opposite $z_k$. Let
\begin{itemize}
\item $d(w_0,x_0)=M(L)/2=\exp(-L/4) + O(\exp(-3L/4))$ (where $M(L)$ is defined in corollary \ref{cor:M}),
\item $d(w_0,y_0)=L/4$,
\item $d(w_1,x_1)=(1+\tau)\exp(-L/4)$ where $|\tau| < E$,
\item $d(w_1,y_1)= L/4 + \rho$ where $|\rho| < E$.
\end{itemize}
See figure \ref{fig:trirect}.

\begin{thm}\label{thm:trirectangle}
Given $\epsilon > 0$ there exists positive numbers $E=E(\epsilon), L_0=
L_0(\epsilon) > 0$ such that if $|\tau|,|\rho| < E$, $L > L_0$ and $Q_0,Q_1$ are as above then
there exists a homeomorphism $F:Q_1 \to Q_0$ such that 
\begin{itemize}
\item $F(w_1)=w_0$, $F(x_1)=x_0$, $F(y_1)=y_0$, $F(z_1)=z_0$.
\item $F$ restricted to any side of $Q_1$ is a similarity.
\item $F$ is $(1+\epsilon)$ bi-Lipschitz.
\end{itemize}
\end{thm}

\begin{proof}
Let $\gamma_0$ be the set of points $q$ in $Q_0$ such that $d(q,\overline{w_0x_0})=L/4$. Similarly let $\gamma_1$ be the set of points $q$ in $Q_1$ such that $d(q,\overline{w_1x_1})=L/4 + \rho$. Then $\gamma_0$ and $\gamma_1$ are segments of an equidistant curve. See figure \ref{fig:trirect}. For $k=0,1$ let $\Delta_k$ be the curvilinear triangle with sides equal to $\gamma_k$, $\overline{y_kz_k}$ and $\overline{p_kz_k}$ where $p_k$ is the point of intersection of $\gamma_k$ and $\overline{x_kz_k}$. Let $Q'_k$ be the complement $Q_k - \Delta_k$.

\begin{figure}[htb]
 \begin{center}
 \ \psfig{file=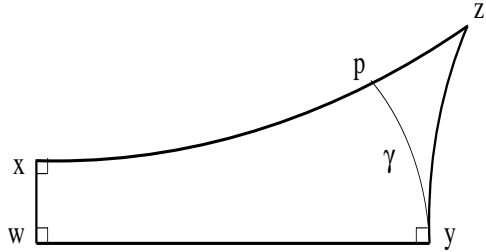,height=1.3in,width=2.5in}
 \caption{The trirectangle $Q_k$}
 \label{fig:trirect}
 \end{center}
 \end{figure}

We first define the map $F$ on $Q'_1$ (with image $Q'_0$) so that $F$ restricted to any
geodesic segment perpendicular to $\overline{w_1x_1}$ is a
similarity. We use rectangular coordinates to compute the bi-Lipschitz
constant of $F$. So let $\R^2_r$ equal $\R^2$ with
the metric
\begin{displaymath}
ds^2 = \cosh^2(y)dx^2 + dy^2.
\end{displaymath}
$\R^2_r$ is isometric to the hyperbolic plane (see \cite{Fen}
page 205). Let $O=(0,0)\in\R^2_r$. Then $(x,0) \in \R^2_r$ is a point
of distance $|x|$ from $O$ and $(x,y)$ is a point at distance $|y|$
from $(x,0)$ such that $\overline{(x,y)(x,0)}$ is perpendicular to $\overline{(0,0)(x,0)}$.

We position the trirectangles $Q_0$ and $Q_1$ as follows. Let $(0,0)=w_1=w_2$, $x_0=(M(L)/2,0)$, $x_1=((1+\tau)\exp(-L/4),0)$, $y_0=(0,L/4)$ and $y_1=(0,L/4+\rho)$. Let
\begin{eqnarray*}
a &=& \frac{d(w_0,x_0)}{d(w_1,x_1)}\\
& =& \frac{\exp(-L/4) + O(\exp(-3L/4))}{(1+\tau)\exp(-L/4)}\\
&=& \frac{1}{1+\tau} + O(\exp(-L/2)).
\end{eqnarray*}
Let
\begin{eqnarray*}
b = \frac{d(w_0,y_0)}{d(w_1,y_1)} = \frac{L/4}{L/4+\rho} = \frac{1}{1+4\rho/L}.
\end{eqnarray*}
We may assume, by taking $E$ small enough and $L$ large enough, that $|a-1|,|b-1|<\epsilon$. Define $F$ on $Q_1'$ by $F(x,y)= (ax,by)$. We claim that this map is $(1+\epsilon)$ biLipschitz if $E$ is small enough and $L$ is large enough. To check this let $Z=Z_{(x,y)}$ be the matrix
\begin{displaymath}
Z = \left[\begin{array}{cc}
\cosh(y) & 0 \\
0 & 1.
\end{array}\right]
\end{displaymath}
Let $||\cdot||_e$ denote the usual Euclidean norm and let $||\cdot||_r$
denote the norm in $\R^2_r$. If $v$ is a vector based at
$(x,y)$ then $||Zv||_e = ||v||_r$. Hence $v$ has Euclidean norm $1$
iff $Z^{-1}v$ has hyperbolic norm $1$. Let $K=K_{(x,y)}=Z_{F(x,y)} DF_{(x,y)} Z_{(x,y)}^{-1}$. So
\begin{displaymath}
K = \left[\begin{array}{cc}
\frac{a\cosh(by)}{\cosh(y)} & 0 \\
&\\
0 & b.
\end{array}\right]
\end{displaymath}

Using the fact that $y \to \cosh(by)/\cosh(y)$ is monotonic it can be shown that for $y\in [0,L/4+\rho]$,
\begin{eqnarray*}
\Big|\frac{\cosh(by)}{\cosh(y)}-1\Big| < \epsilon.
\end{eqnarray*}
(for $E$ small enough and $L$ large enough). So $||K - I||_\infty < \epsilon$ (if $E$ is small enough and $L$ is large enough). This implies that $F$ restricted to $Q_1'$ is a $(1+\epsilon)$ bi-Lipschitz map onto $Q_0'$. Note also that $F$ restricted to $\gamma_1$ is a similarity onto $\gamma_0$ (with respect to the intrinsic metrics on $\gamma_0$ and $\gamma_1$).

Next we define $F$ on the curvilinear triangle $\Delta_1$. We require that $F$ restricted to $\overline{p_1z_1}$ is a similarity and $F$ restricted to $\overline{y_1z_1}$ is a similarity. If $\omega$ is a geodesic segment perpendicular to $\gamma_1$ contained in $\Delta_1$ then we require $F$ restricted to $\omega$ to be a similarity. Since $F$ is already defined on the boundary of $\Delta_1$ (and since $\overline{p_1z_1}$ is perpendicular to $\gamma_1$) this determines $F$ completely. It is clear that $F|_{\Delta_1}$ has continuous derivatives. Since $\Delta_1$ is contained is a circle of radius 100 (for all $E$ small enough and $L$ large enough) it is clear that as $E$ tends to zero and $L$ tends to infinity, the map $F$ restricted to $\Delta_1$ tends to an isometry. Thus by choosing $E$ small enough and $L$ large enough we may assume that $F|_{\Delta_1}$ is $(1+\epsilon$) biLipschitz.

\end{proof}

\begin{proof}(of the Bi-Lipschitz Pants Theorem \ref{thm:bilipschitz pants})

We first decompose $P$ into a union of two right-angled hexagons $H_1$ and $H_2$ in the standard way. To be specific, $H_1$ and $H_2$ are obtained from $P$ by cutting along three distinct geodesic arcs where each arc is the shortest path between two distinct boundary components. The three altitudes (see subsection \ref{subsection:altitudes}) decompose each hexagon into six trirectangles that satisfy the bounds of theorem \ref{thm:trirectangle}. 

In a similar way, we decompose $P_L$ into twelve trirectangles. We
define a map $F: P \to P_L$ so that $F$ restricted to any of the
twelve trirectangles is a $(1+\epsilon)$ biLipschitz map onto a
trirectangle of $P_L$ whose restriction to the boundary is a similarity. This is possible by theorem \ref{thm:trirectangle} above (if $\epsilon_1$ is small enough and $L$ is large enough). Since the special points are contained in the vertices of the trirectangles $F$ maps special points to special points and restricts to a similarity on the boundary components. 

\end{proof}

\section{Annuli}\label{section:annuli}

\begin{proof}(of theorem \ref{thm:bilipschitz annuli})

We will define a function ${\tilde F}$ on the plane $\H^2$ which will descend
to a function $F$ of the annulus $\A$ satisfying the conclusions of
the theorem. If $(x,y) \in \H^2$ (in the upperhalf plane model) define $r=r(x,y)\in (0,\infty)$ and $\theta(x,y) \in (0,\pi)$ by
\begin{eqnarray*}
(x,y)= (r\cos(\theta), r\sin(\theta)).
\end{eqnarray*}
Let $\gamma$ denote the geodesic with endpoints $0$ and $\infty$. Let $\alpha \in (0,\pi/2)$. Then the equation $\theta(x,y)=\alpha$ defines a curve equidistant from $\gamma$. Choose $\alpha$ so that this distance equals $w$ and so that $\alpha \le \pi/2$ (this ensures that the curve is on the right side of $\gamma$).

Define ${\tilde F}={\tilde F}_{\alpha,t}$ by
\begin{displaymath}
{\tilde F}(re^{i\theta}) = f(\theta)re^{i\theta}
\end{displaymath}
where $f=f_{\alpha,t}:(0,\pi) \to \R$ is defined by
\begin{displaymath}
f(\theta)=\left\{\begin{array}{ll}
1 & \pi-\alpha \le \theta < \pi\\
\frac{\theta-\alpha}{\pi-2\alpha} + \Big(1 -
\frac{\theta-\alpha}{\pi-2\alpha}\Big)e^t & \alpha < \theta <
\pi-\alpha\\
e^t & 0 < \theta \le \alpha
\end{array}\right\}
\end{displaymath}
Since ${\tilde F}$ preserves Euclidean lines through the
origin, ${\tilde F}$ commutes with the action of $A_l$. Hence ${\tilde
  F}$ descends to map $F:\A \to \A$ on the annulus. From the
definition, we see that $F$ satisfies the first two properties in the
conclusion of theorem \ref{thm:bilipschitz annuli}. It suffices to show that the bi-Lipschitz constant of ${\tilde F}$ is less than $1+\epsilon$ since it equals the bi-Lipschitz constant of $F$.

Recall that the metric on the upper-half plane is $ds^2 =
(dx^2+dy^2)/y^2$. Also
\begin{displaymath}\begin{array}{ll}
(dx^2+dy^2)/y^2 &= (dr^2 + r^2d\theta^2)/r^2\sin^2(\theta).
\end{array}
\end{displaymath}
Let $\R^2_c$ denote $(0,\infty)\times (0,\pi)$ with the metric $ds^2=(dr^2 + r^2d\theta^2)/r^2\sin^2(\theta)$. Define the matrix $Z_{r,\theta}$ by

\begin{displaymath}
Z_{r,\theta}=\left[\begin{array}{cc}
1/(r\sin(\theta))   & 0 \\
0 & 1/\sin(\theta)
\end{array}\right].
\end{displaymath}

Let $||\cdot||_e$ denote the usual Euclidean norm and let $||\cdot||_c$
denote the norm in $\R^2_c$. If $v$ is a vector based at
$(x,y)$ then $||Zv||_e = ||v||_c$. Hence $v$ has Euclidean norm $1$
iff $Z^{-1}v$ has hyperbolic norm $1$. Let $K=K_{(r,\theta)}$ be the matrix $K=Z_{F(r,\theta)} D{\tilde F}_{(r,\theta)} Z_{(r,\theta)}^{-1}$ where $D{\tilde F}_{(r,\theta)}$ is the matrix representing the differential of ${\tilde F}$ in the coordinates $(r,\theta)$. Then it suffices to show that $||K-I||_\infty < \epsilon$ for $E$ small enough and $L$ large enough.

The derivative of $f$ is given by:
\begin{displaymath}
f'(\theta)=\left\{\begin{array}{ll}
0 & \pi-\alpha \le \theta < \pi\\
\frac{1-e^t}{\pi-2\alpha} & \alpha < \theta <
\pi-\alpha\\
0 & 0 < \theta \le \alpha
\end{array}\right\}
\end{displaymath}
The differential of ${\tilde F}$ with respect to the coordinates $(r,\theta)$ is:
\begin{displaymath}
D{\tilde F}_{(r,\theta)} = \left[\begin{array}{cc}
f(\theta)   & f'(\theta)r \\
0  & 1
\end{array}\right].
\end{displaymath}

The matrix $K$ is:
\begin{eqnarray*}
K &=& \left[\begin{array}{cc}
1 & f'(\theta)/f(\theta) \\
0  & 1 
\end{array}\right]
\end{eqnarray*}

So it suffices to show that $|f'(\theta)/f(\theta)| <\epsilon$ is $E$ is small enough and $L$ is large enough. Since $t\ge 0$ by hypothesis, $1 \le f(\theta) \le e^t$. So
\begin{displaymath}
|f'(\theta)/f(\theta)| \le |f'(\theta)| = \frac{e^t-1}{\pi-2\alpha}.
\end{displaymath}

Now we estimate $\alpha$. Let $z=(\cos(\alpha),\sin(\alpha)) \in \H^2$. By definition of $\alpha$ the distance between $z$ and $(0,1)$ equals $w$. Recall that if $z_1=(x_1,y_1), z_2=(x_2,y_2) \in \H^2$ then the hyperbolic distance between them is given by:
\begin{eqnarray*}
\cosh(d(z_1,z_2)) = 1 + \frac{||z_1-z_2||_e^2}{2y_1y_2}.
\end{eqnarray*}
So
\begin{displaymath}\begin{array}{ll}
\cosh(d((0,1),z)) &= 1 + \frac{||(0,1)-z||_e^2}{2\sin(\alpha)}\\
              &= 1 + \frac{\cos^2(\alpha) +
              (1-\sin(\alpha))^2}{2\sin(\alpha)}\\
              &= \frac{1}{\sin(\alpha)}.
\end{array}
\end{displaymath}
So, $\sin(\alpha) = 1/\cosh(d((0,1),z)) = 1/\cosh(w)$. Recall that $t \le E\exp(-L/4)$ and $w \ge (1/2)\exp(-L/4)$. So
\begin{eqnarray*}
\cos^2(\alpha) &=& 1 - \frac{1}{\cosh^2(w)} = \tanh^2(w)\\
             &\ge& (1/4)\exp(-L/2) + O(\exp(-L)).
\end{eqnarray*}

So $\cos(\alpha) \ge (1/2)\exp(-L/4) + O(\exp(-3L/4))$. This implies
that
\begin{displaymath}
\alpha \le \pi/2 -  (1/2)\exp(-L/4) + O(\exp(-3L/4)).
\end{displaymath}
So,
\begin{eqnarray*}
\frac{e^t-1}{\pi-2\alpha} & \le& \frac{ E \exp(-L/4) +
  O(\exp(-L/2))}{\exp(-L/4) + O(\exp(-3L/4))}\\
                          &=&E + O(\exp(-L/4)).
\end{eqnarray*}

Thus if $E$ is small enough and $L$ is large enough $F$ is $(1+\epsilon)$ bi-Lipschitz as claimed.

\end{proof}

\part{Incompressibility}\label{part:incompressible}

\section{Introduction}

The goal of this part is to prove theorem \ref{thm:incompressible}. The strategy is this: first we define a set of weights on the $1$-skeleton of a natural hexagonal tiling that refines the pants decomposition of $S$. The weights approximate hyperbolic distance. To show that $j:S \to \M$ is $\pi_1$-injective it suffices to show that $j({\hat \gamma})$ is homotopically nontrivial in $\M$ for any homotopically nontrivial curve ${\hat \gamma} \subset S$. We homotope ${\hat \gamma}$ into the graph and require it to have least weight among all curves in the graph homotopic to it. Then we lift ${\hat \gamma}$ to the universal cover of $S$ and then push it forward to $\H^3$ via a lift of $j$. We straighten so that it is a piecewise geodesic curve $\gamma$. With a few exceptions, we associate to each edge of $\gamma$ a plane defined using altitudes of right-angled hexagons. Then we show that successive planes are disjoint. This is used to show that $\gamma$ cannot be a closed curve and therefore $j({\hat \gamma})$ is homotopically nontrivial.

We will need a bit of terminology/notation. We have defined a right-angled hexagon as an ordered list of oriented geodesics (section \ref{right angled polygons}). We will, by abuse of terminology, also call any 6-sided polygon such that any pair of adjacent edges meet in a right angle a right-angled hexagon. If $\hG=(\gG_1,...,\gG_6)$ is a right-angled hexagon (as defined in section \ref{right angled polygons}) then the polygon associated to $\hG$, $poly(\hG)$ is the right-angled hexagon with vertices $v_1,...,v_6$ where $v_i$ is the intersection point $\gG_i\cap\gG_{i+1}$ for all $i$ mod 6. 


If $a,b,c \in \H^3$ then $\angle (a,b,c)$ denotes the interior angle at $b$ of the triangle with vertices $a$, $b$ and $c$. We let $\overline{ab}$ denote the geodesic segment from $a$ to $b$. We will at times abuse notation by confusing $\overline{ab}$ with the distance between $a$ and $b$.

\section{Graphs on the Surface}

Let $\M^3$ be a closed hyperbolic 3-manifold and let $\Gamma
<PSL_2(\C)$ be a discrete group such that $\H^3/\Gamma = \M^3$. Let $\epsilon, L, {\hat T}>0$ and let $j:S \to \M^3$ be a map of a closed surface $S$ into $\M^3$ such that $S$ has a labeled pants decomposition $\P$ satisfying the bounds in the conclusion of theorem \ref{thm:surface subgroup case}. Assume that $2{\hat T}\exp(-L/4) < \epsilon<1/2< L/4$. Later we may choose $L$ larger or $\epsilon$ smaller if necessary.

For convenience we define a hyperbolic structure on the surface $S$ by the following. We require that every curve $\xi \in \P^*$ is a geodesic in $S$ of length $L$. Let $twist_S(\xi)$ denote the twist parameter of $\xi \in \P^*$ with respect to the metric on $S$ and $twist_j(\xi)$ denote the twist parameter of $\xi$ with respect to $j$. Then we require that $twist_S(\xi)=\Re(twist_j(\xi))$.

 Let $\T$ be the tiling of $S$ by right-angled hexagons that respects the pants decomposition $\P$. To be precise, $\T$ is the collection of right-angled hexagons in $S$ with pairwise disjoint interiors whose union is all of $S$ and such that every $H \in \P$ is a union of 2 hexagons in $\T$. Thus every hexagon $\hH \in \T$ has three alternating sides of length $L/2$. We call such sides ``long sides'' of $\hH$. The other three sides are called ``short sides'' of $\hH$.



We associate to each hexagon $\hH \in \T$ a graph $G(\hH)$ as follows. The vertices of the graph are the 6 vertices of $\hH$ and the three midpoints of the long sides of $\hH$. There exists an edge in $G(\hH)$ between every pair of vertices $v_1,v_2$ except if $v_1$ and $v_2$ lie on the same long side and either $v_1$ or $v_2$ is a midpoint of that side.

We let $G_0(\T) \subset S$ be the graph equal to the union
\begin{displaymath}
G_0(\T) = \cup_{\hH \in \T} \, G(\hH).
\end{displaymath}

Let $G(\T)$ be the (smallest) graph containing $G_0(\T)$ such that the following holds. The vertex set of $G(\T)$ equals the vertex set of $G_0(\T)$. Suppose that $\xi \in \P^*$. Suppose that $v_1,v_2 \in \xi$ are vertices of $G(\T)$ and $d(v_1,v_2) \le {\hat T}\exp(-L/4)$. Then there exists an edge $e$ between $v_1$ and $v_2$ in $G(\T)$. These are all the edges in $G(\T)$ that are not in $G_0(\T)$. For example, if all the associated twist parameters of the four-holed sphere decomposition $\P$ equal zero (with respect to the metric on $S$) then $G_0(\T)=G(\T)$.

Let $\pi:{\tilde S} \to S$ be the universal cover of $S$. Let ${\tilde j}:{\tilde S} \to \H^3$ be a lift of $j$. Let ${\tilde \T} = \pi^{-1}(\T)$ be a hexagonal tiling of ${\tilde S}$. Let $G({\tilde \T})$ be the graph associated to ${\tilde \T}$ as above so that $G({\tilde \T})=\pi^{-1}(G(\T))$. Let ${\tilde j}(G({\tilde \T}))=G_0 \subset \H^3$. Let $G$ be the graph whose vertex set coincides with the vertex set of $G_0$ such that every edge of $G$ is a geodesic segment in $\H^3$.

Recall that an $(L,\epsilon)$ nearly-symmetric hexagon $\hG=(\gG_1,...,\gG_6)$ is a standardly oriented right-angled hexagon so that if $\wG_i = \mu(\gG_{i-1},\gG_{i+1};\gG_i)$ (for $i$ mod 6) then there exists numbers $\rho_i \in \C$ (for $i=1,3,5$) such that $|\rho_i|<\epsilon$ and $\wG_i=L/2 + \rho_i/2 + i\pi$.

Recall that an {\bf altitude} of a right-angled hexagon $\hG$ is a geodesic that is perpendicular to two opposite sides of the hexagon $\hG$. By abuse of language, we will also refer to an altitude as the shortest geodesic segment between two opposite sides of $\hG$.

After homotoping $j$ if necessary we may assume the following:

\begin{itemize}
\item For every hexagon $\hH \in {\tilde \T}$ the image ${\tilde j}(\partial \hH) \subset \H^3$ equals $poly(\hG)$ for some $(L,\epsilon)$ nearly-symmetric hexagon $\hG$.
\item If $\hH \in {\tilde \T}$ and ${\tilde j}(\partial \hH)=poly(\hG)$ then the midpoints of the edges of $\hH$ are mapped to the endpoints of the altitudes of $\hG$.
\end{itemize}

\section{Curves on the Surface}
 Let ${\hat \gamma}$ be any arbitrary homotopically nontrivial curve in $S$. To prove theorem \ref{thm:incompressible} it suffices to show that $j({\hat \gamma})$ is homotopically non-trivial.


If $e$ is an edge of the graph $G(\hH)\subset S$ for some hexagon $\hH \in \T$ then let $weight(e)$ be the length of $e$ with respect to the hyperbolic metric on $S$. If $e$ is an edge of $G(\T)$ that is not an edge of $G(\hH)$ for any hexagon $\hH \in \T$ then we define $weight(e)=0$. We define the weight of a path in $G(\T)$ to be the sum of the weights of all the edges contained in the path. We pull back and push forward these weights to obtain weights on the edges of $G({\tilde T})$ and $G$.

After homotoping ${\hat \gamma}$ if necessary we may assume that ${\hat \gamma}$ lies inside the graph $G(\T)$ and has least weight among all curves in the graph homotopic to $\hat \gamma$ (through homotopies in the surface $S$).

Let ${\tilde \gamma} \subset {\tilde S}$ be a lift of ${\hat \gamma}$ and let $\gamma' = {\tilde j}({\tilde \gamma}) \subset \H^3$. Let $\{v_i\}_{i\in\Z}$ be the sequence of vertices of $G$ traversed by $\gamma'$. Let $\gamma$ be the piecewise geodesic path with vertices $\{v_i\}_{i\in\Z}$. To prove theorem \ref{thm:incompressible} it suffices to show that $\gamma$ is not a closed curve since this implies that $\gamma'$ is not a closed curve which implies $j({\hat \gamma})$ is homotopically nontrivial. 

To do this we will associate to each edge $e$ of $\gamma$ a plane transversal to it and then show the set of all such planes is pairwise disjoint (with a few exceptions). We will then use this to show that $\gamma$ is not a closed curve.

If $\hG$ is an $(L,\epsilon)$ nearly symmetric hexagon in $\H^3$ then we define an {\bf altitude plane} of $\hG$ to be a plane in $\H^3$ that contains an altitude of $\hG$ and the short side of $\hG$ which is perpendicular to the altitude. Here a short side of $\hG$ is a side of real length less than $2(1+\epsilon)\exp(-L/4)+O(\exp(-3L/4))$. For example, if $\hG$ is defined as in section \ref{section:L/2 hexagons} then for $k$ even $\gG_k$ are the geodesics containing the short sides of $\hG$. There are only three altitude planes of $\hG$ corresponding to the three short sides of $\hG$.

\begin{figure}[htb]
 \begin{center}
 \ \psfig{file=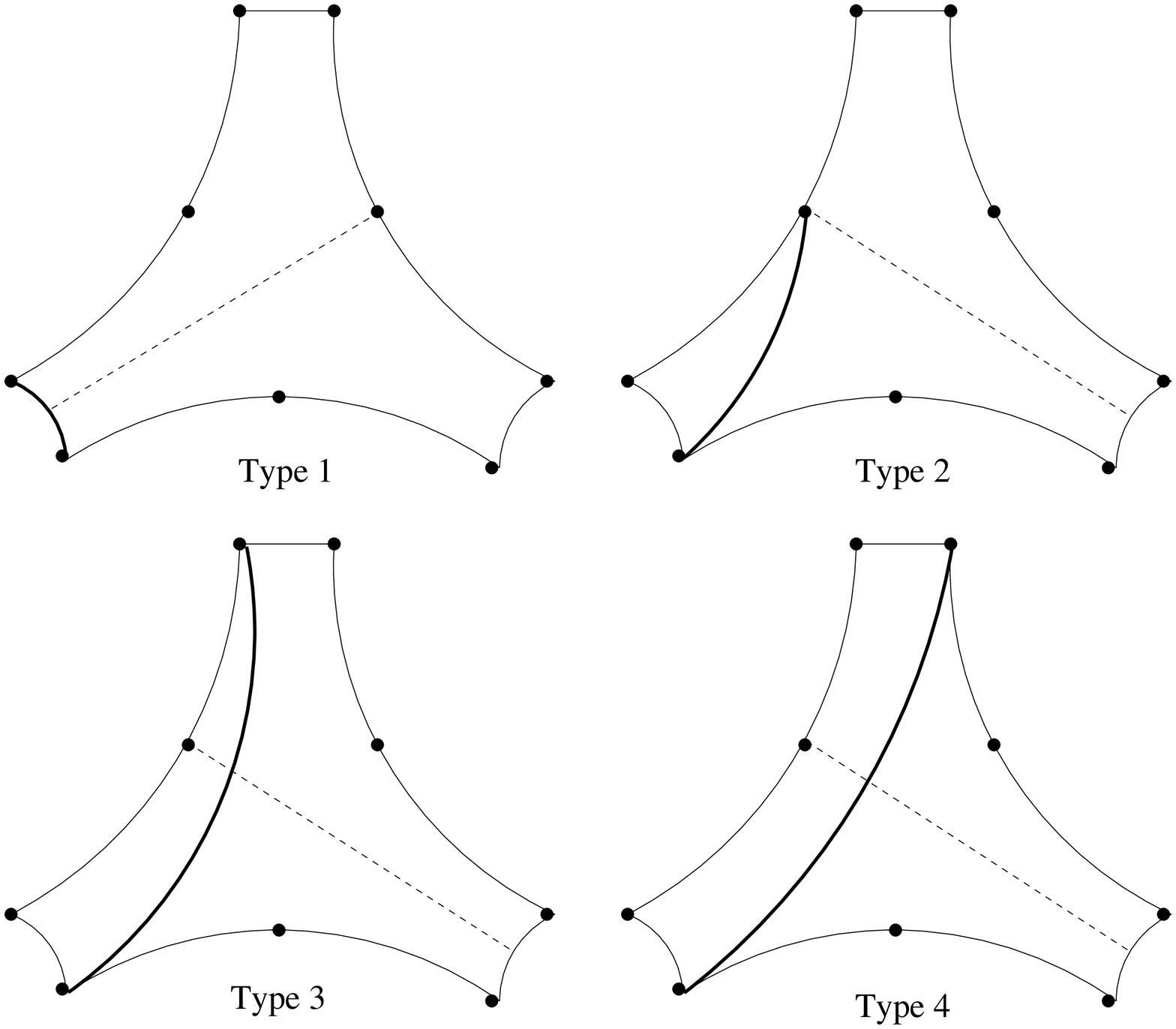,height=4in,width=4in}
 \caption{Types of edges of $\gamma$ (thick lines) and their associated altitude planes (dashed lines)}
 \label{fig:singleedge}
 \end{center}
 \end{figure}

\begin{figure}[htb]
 \begin{center}
 \ \psfig{file=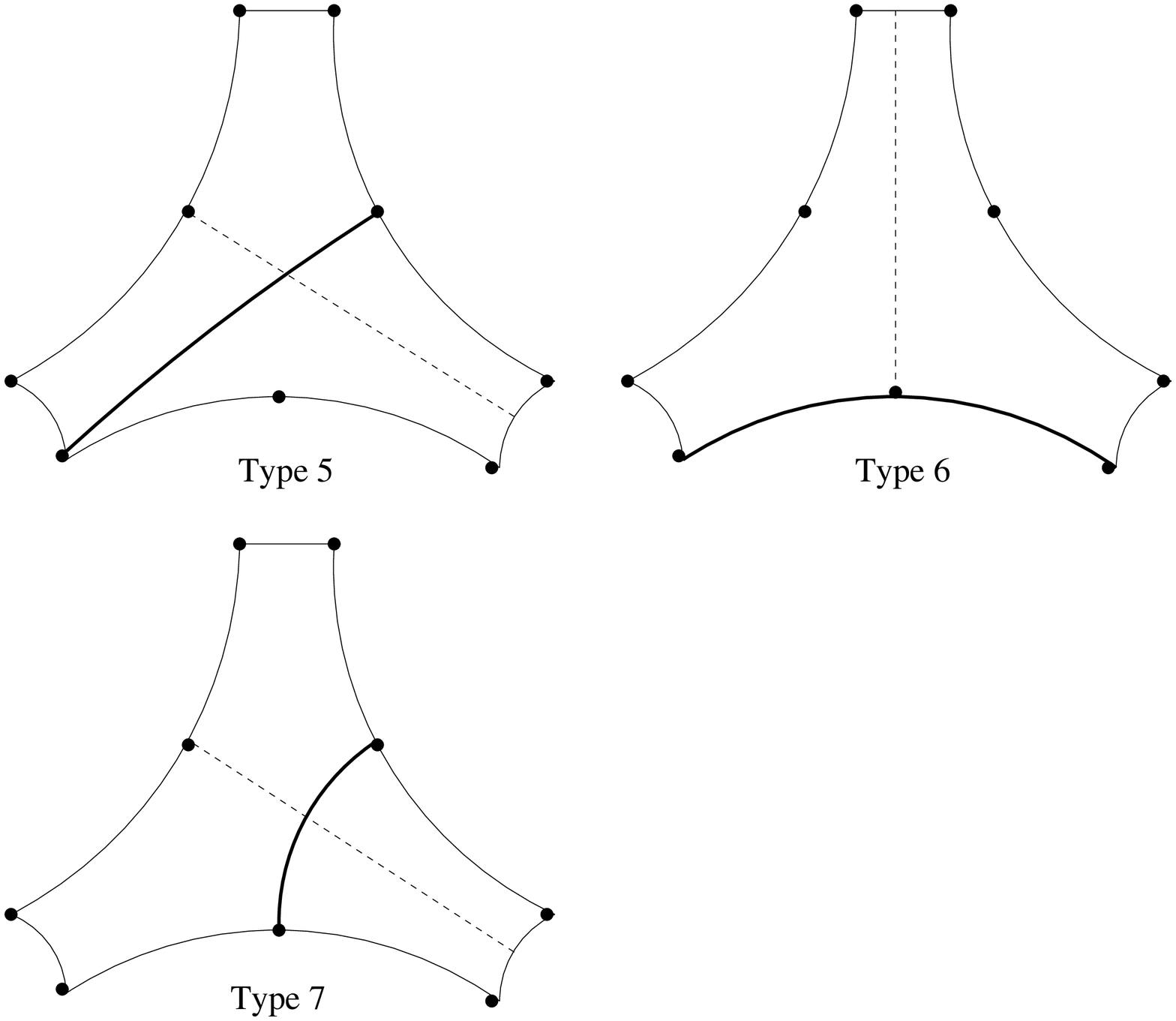,height=4in,width=4in}
 \caption{Types of edges of $\gamma$ (thick lines) and their associated altitude planes (dashed lines)}
 \label{fig:singleedge2}
 \end{center}
 \end{figure}

Suppose $e$ is an edge of $\gamma$ with endpoints $a$ and $b$. Suppose that the weight of $e$ (induced from the weighting of the edges of $G(\T)$) is nonzero. Suppose that $a$ and $b$ are contained in a hexagon $\hH$ of the image (so $\hH={\tilde j}({\tilde \hH})$ for some lift ${\tilde \hH}$ of a hexagon ${\hat \hH}$ in $\T$). Then we associate to $e$ an altitude plane $\Pi'(e)$ of $\hH$ as shown in figures \ref{fig:singleedge} and \ref{fig:singleedge2}. For example, if $e$ is contained in a short edge of the hexagon $poly(\hH)$ then $e$ is said to be Type 1 and the altitude plane $\Pi'(e)$ intersects it transversally. If $e$ connects a vertex of the boundary of $poly(\hH)$ to an altitude's endpoint in such a way that the endpoints of $e$ separate one vertex of $poly(\hH)$ from the others then $e$ is said to be Type 2 and $\Pi'(e)$ passes through the endpoint of $e$ that is also the endpoint of an altitude. There is some ambiguity if $e$ is of Type 6 because in that case $e$ may be contained in two different hexagons $\hH_1, \hH_2$ and the associated planes may be different. If this is the case we arbitrarily choose one of the two hexagons containing $e$ to define $\Pi'(e)$.

Let $\{e_i\}_{i\in\Z}$ denote the sequence of positive-weight edges of $\gamma$. Define $\Pi(e_i)$ for $i\in\Z$ as follows. If $e_i=(a,b)$ and $e_{i+1}=(c,d)$ are both Type 2 edges and $b$ is an endpoint of an altitude (of the hexagon containing $e_i$) then we let $\Pi(e_i)=\emptyset$. Here $b$ is the endpoint of $e_i$ that is closest to $e_{i+1}$. Otherwise we let $\Pi(e_i)=\Pi'(e_i)$. Note that if $\Pi(e_i)=\emptyset$ then both $\Pi(e_{i-1})$ and $\Pi(e_{i+1})$ are nonempty.

\begin{lem}\label{lem:gamma}
There exists positive numbers $\epsilon_0,L_0>0$ such that if $\epsilon<\epsilon_0$, ${\hat T}>0$ is a given constant that may depend on $\epsilon$, $L>L_0$ and $\gamma$ is defined as above then for all $i\in\Z$, $\Pi(e_i)\cap \Pi(e_{i+1})=\emptyset$ and if $\Pi(e_i)=\emptyset$ then $\Pi(e_{i-1})\cap \Pi(e_{i+1})=\emptyset$. 
\end{lem}

We will prove the above lemma after the one below.

\begin{figure}[htb]
 \begin{center}
 \ \psfig{file=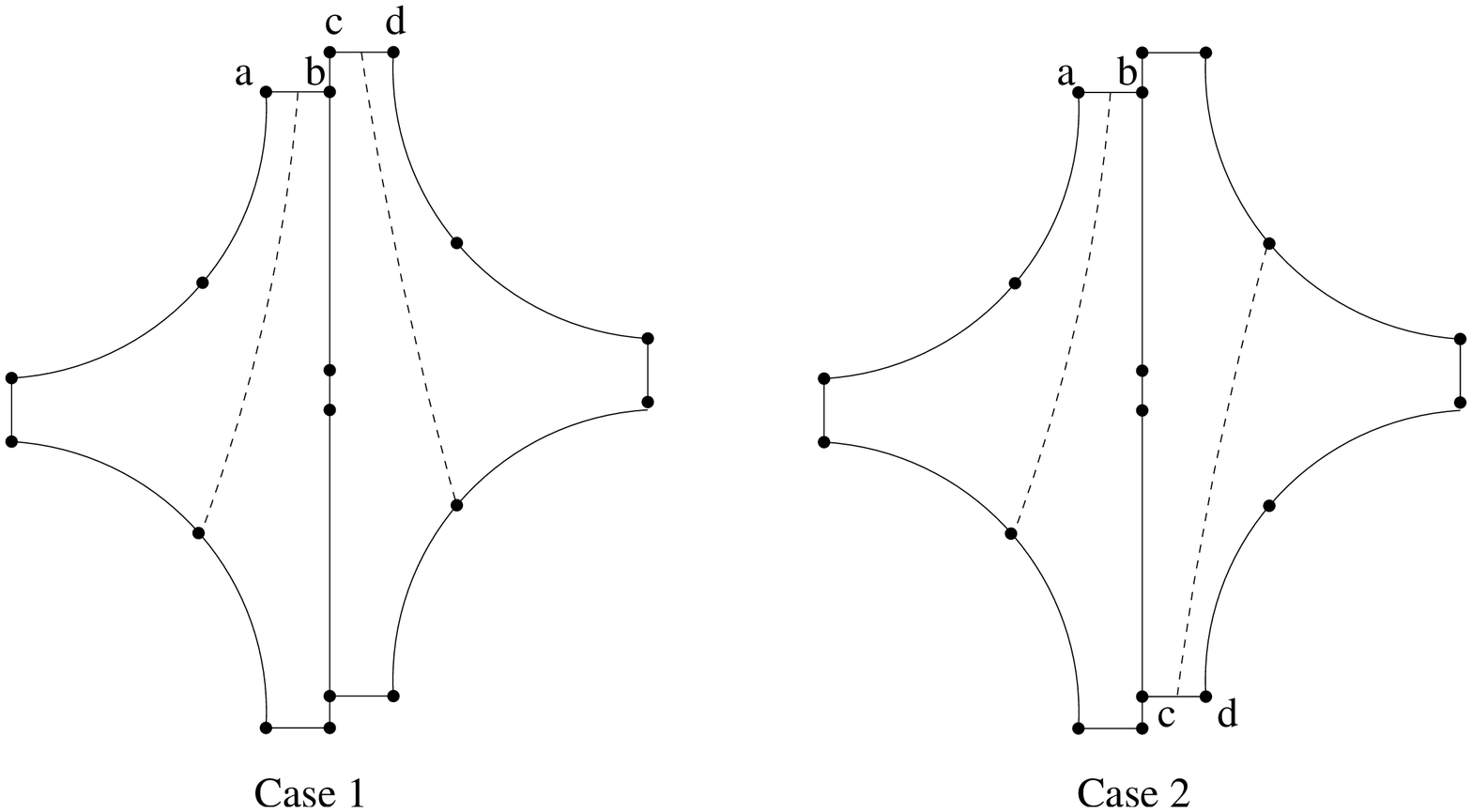,height=3.2in,width=4.3in}
 \caption{Pairs of altitude planes indicated in dashed lines}
 \label{fig:altlong}
 \end{center}
 \end{figure}

\begin{figure}[htb]
 \begin{center}
 \ \psfig{file=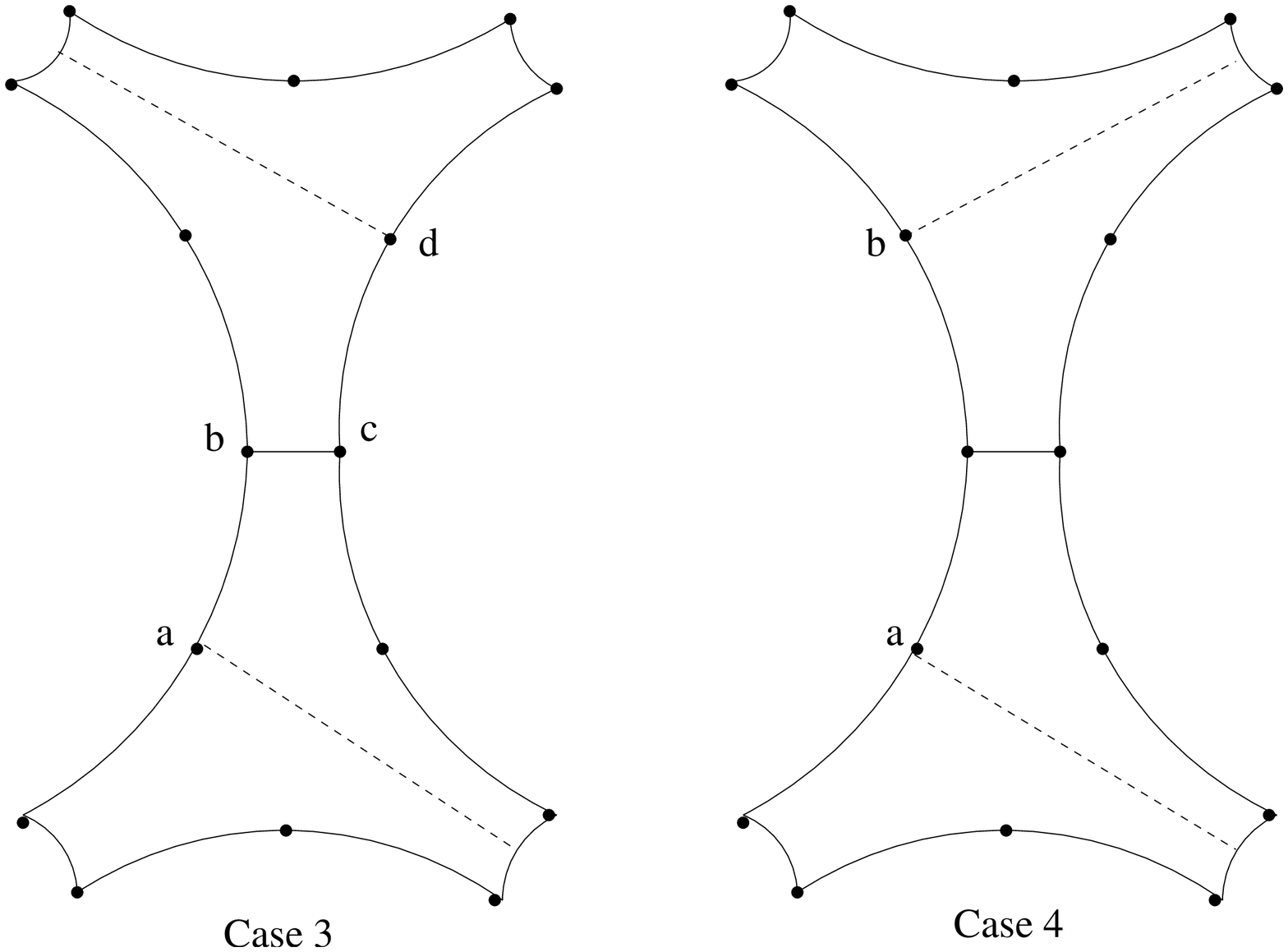,height=3.2in,width=4.3in}
 \caption{Pairs of altitude planes indicated in dashed lines}
 \label{fig:altshort}
 \end{center}
 \end{figure}

\begin{figure}[htb]
 \begin{center}
 \ \psfig{file=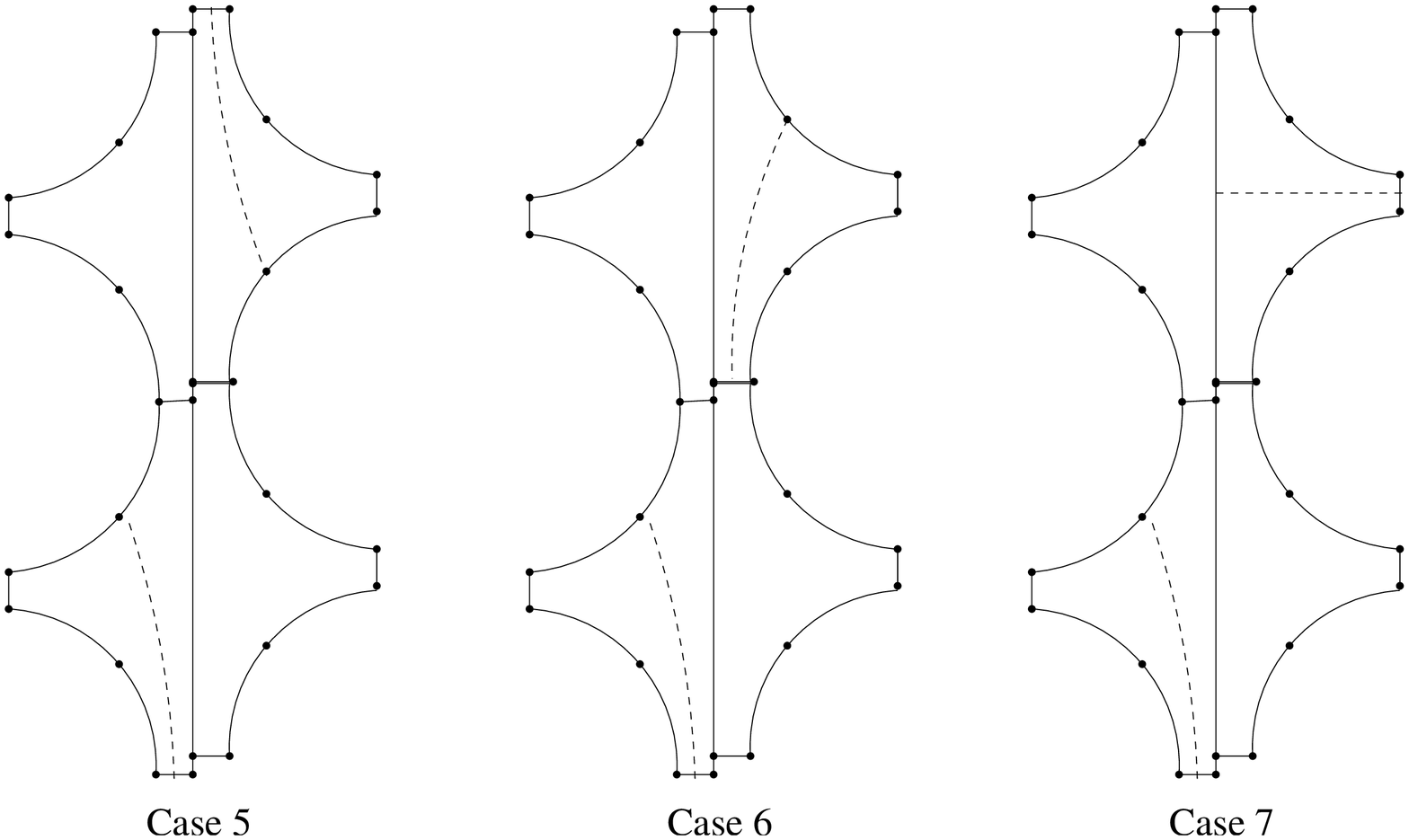,height=3in,width=4in}
 \caption{Pairs of altitude planes indicated in dashed lines}
 \label{fig:altcat}
 \end{center}
 \end{figure}

\begin{figure}[htb]
 \begin{center}
 \ \psfig{file=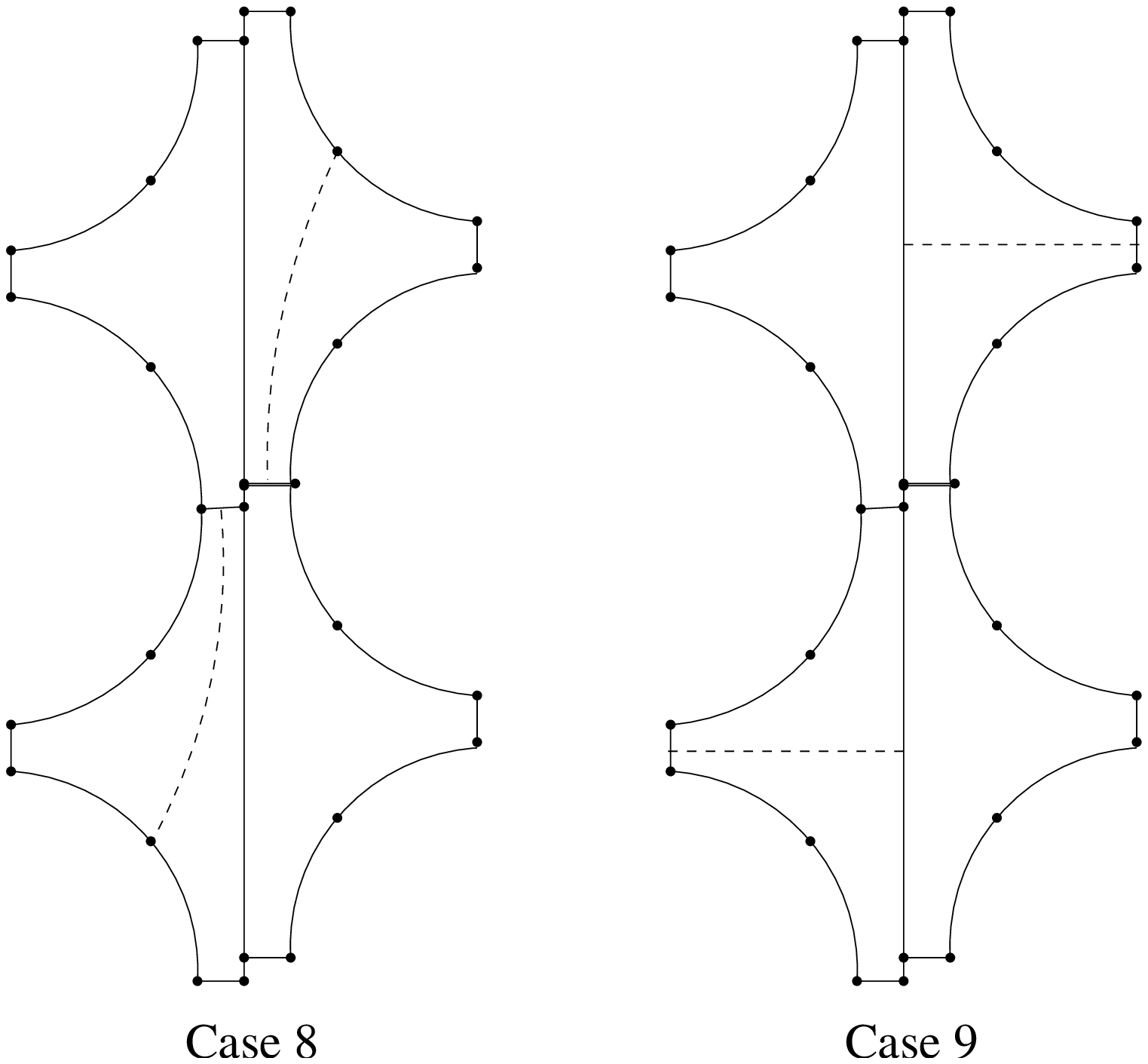,height=3in,width=3in}
 \caption{Pairs of altitude planes indicated in dashed lines}
 \label{fig:altcat2}
 \end{center}
 \end{figure}

\begin{lem}\label{lem:altitude planes}
There exists positive numbers $\epsilon_0,L_0>0$ such that if $\epsilon<\epsilon_0$, ${\hat T}>0$ is a given constant that may depend on $\epsilon$, $L>L_0$ and $\gamma$ is defined as above then each pair of altitude planes represented in figures \ref{fig:altlong}-\ref{fig:altcat2} is disjoint.
\end{lem}

\begin{proof}

{\bf Cases 1 \& 8}: These two cases are the exact same (the planes themselves are the same). They are distinguished only to facilitate proving lemma \ref{lem:gamma}. Let $a,b,c,d$ be the points shown in the Case 1 figure \ref{fig:altlong}. Let $\Pi_1$ be the altitude plane perpendicular to $\overline{ab}$ and let $\Pi_2$ be the altitude plane perpendicular to $\overline{cd}$. Let $\Pi_3$ be the plane perpendicular to $\overline{ab}$ that contains $b$. It is immediate that $\Pi_3$ is disjoint from $\Pi_1$. We will show that $\Pi_3$ is also disjoint from $\Pi_2$. Let $e$ be the intersection point of $\Pi_2$ with $\overline{cd}$. Then by estimates from lemma \ref{lem:altitudes} there exists $\epsilon_0,L_0>0$ such that if $\epsilon<\epsilon_0$ and $L>L_0$ then $d(c,e) \ge (1-\epsilon)\exp(-L/4)$. So we will assume this is the case.

 Let $\Pi_4$ be the plane perpendicular to $\overline{cd}$ that passes through $c$. Then $\Pi_4$ and $\Pi_3$ intersect in an angle equal to the imaginary part of the twist parameter $twist_j(\sigma)$ where $\sigma \in \P^*$ is such that ${\tilde j}({\tilde \sigma})\supset \overline{bc}$ where ${\tilde \sigma}$ is a lift of $\sigma$ to the universal cover ${\tilde S}$ of $S$. Thus the angle between $\Pi_4$ and $\Pi_3$ is at most $\epsilon\exp(-L/4)$.

 Suppose for a contradiction that there is a point $x \in \Pi_2 \cap \Pi_3$. Consider the triangle with vertices $x,c,e$. By definition $\angle(c,e,x) = \pi/2$. $|\angle(x,c,e)-\pi/2|$ is at most the angle between $\Pi_3$ and $\Pi_4$. So by the above estimates $\angle(x,c,e)\ge \pi/2 - \epsilon\exp(-L/4)$. By the law of cosines
\begin{displaymath}
\cosh(\overline{ce}) = \frac{\cos(\pi/2)\cos(\angle(x,c,e)) + \cos(\angle(c,x,e))}{\sin(\pi/2)\sin(\angle(x,c,e))}.
\end{displaymath}
Thus
\begin{eqnarray*}
\cos(\angle(c,x,e)) &=&  \sin(\angle(x,c,e))\cosh(\overline{ce})\\
                    &\ge& \sin(\pi/2-\epsilon\exp(-L/4))\cosh((1-\epsilon)\exp(-L/4)) \\
                    &=& [1-(1/2)\epsilon^2\exp(-L/2)+O(\exp(-L))]\\
&&\times[1+(1/2)(1-\epsilon)^2\exp(-L/2)+O(\exp(-L)]\\
                    &=&1 + (1/2)(1-2\epsilon)\exp(-L/2) + O(\exp(-L)).
\end{eqnarray*}
If $\epsilon < 1/2$ and $L$ is large enough the above implies that $\cos(\angle(c,x,e)) > 1$ which is a contradiction. Hence the two planes do not intersect as claimed. Since $\Pi_3$ separates $\Pi_1$ from $\Pi_2$ this implies that $\Pi_1\cap \Pi_2=\emptyset$. 


{\bf Cases 2, 5 \& 6}: Case 2 and 6 are essentially the same. The proof of Case 5 is very similar to the proof of Case 2 so we will just prove Case 2. Let $a,b,c,d$ be the points indicated in Case 2 figure \ref{fig:altlong}. Let $\Pi_1$ be the altitude plane perpendicular to $\overline{ab}$ and let $\Pi_2$ be the altitude plane perpendicular to $\overline{cd}$. 


Let $\hG=(\gG_1,...,\gG_6)$ be the standardly oriented right-angled
hexagon satisfying
\begin{itemize}
\item $\gG_1 \supset \overline{ab}$; 
\item $\gG_2$ contains the altitude perpendicular to $\overline{ab}$;
\item $\gG_4$ contains the altitude perpendicular to $\overline{cd}$;
\item $\gG_5 \supset \overline{cd}$;
\item $\gG_6 \supset \overline{bc}$. 
\end{itemize}
Let $\wG_i=\mu(\gG_{i-1},\gG_{i+1};\gG_i)$ (for all $i$ mod 6). By lemma \ref{lem:L/2 hexagon} if $\epsilon$ is small enough and $L$ is large enough the following estimates hold.
\begin{itemize}
\item For $k=1,5$, $\wG_k = x_k\exp(-L/4) + i\pi$ where $x_k \in \C$ satisfies $|1-x_k| < \epsilon$;
\item $\wG_6 = L/2 + \rho$ where $\rho \in \C$ satisfies $|\rho|<\epsilon$. 
\end{itemize}
Here we have used the assumption that $\epsilon>2{\hat T}\exp(-L/4)$. 

Suppose for a contradiction that $z$ is a point in $\Pi_1\cap
\Pi_2$. Let $v_k$ be the intersection point $\gG_k \cap \gG_{k+1}$ (for
$k$ mod 6). Then the triangle $zv_2v_3$ has the following properties.
\begin{itemize}
\item $v_2v_3 = |\Re(\wG_3)|$.
\item $|\angle(z,v_2,v_3)-\pi/2| \le |\Im(\wG_2-i\pi)|$.
\item $|\angle(v_2,v_3,z)-\pi/2| \le |\Im(\wG_4-i\pi)|$.
\end{itemize}
For example, the plane $\Pi_3$ perpendicular to $\gG_3$ containing $v_3$ makes an angle $|\Im(\wG_4-i\pi)|$ with the plane $\Pi_2$. This implies the third inequality above; the second inequality is similar.

By the law of cosines,
\begin{eqnarray*}
\cosh(\wG_3) &=& \cosh(\wG_1)\cosh(\wG_5)+\sinh(\wG_1)\sinh(\wG_5)\cosh(\wG_6)\\
             &=& 1 + x_1x_5\exp(-L/2)(1/2)\exp(L/2+\rho) +
O(\exp(-L/2))\\
             &=& 1+ (x_1x_5/2)\exp(\rho) + O(\exp(-L/2)).
\end{eqnarray*}
So, if $L$ is large enough we may write $\cosh(\wG_3)=(3/2)(1+\tau)$ where $\tau \in \C$ is such that $|\tau|<4\epsilon$. By the law of sines,
\begin{eqnarray*}
\sinh^2(\wG_2)&=&\frac{\sinh^2(\wG_5)\sinh^2(\wG_6)}{\sinh^2(\wG_3)}\\
              &=&\frac{x_5^2\exp(-L/2)(1/4)\exp(L+2\rho) + O(1)}{(9/4)(1+\tau)^2-1}\\
&=&\frac{x_5^2\exp(L/2+2\rho) + O(1)}{9(1+\tau)^2-4}.
\end{eqnarray*}
So if $L$ is large enough we may write $\wG_2 = L/4 + (1/2)\ln(4/5) + \tau_2 + i\pi$ where $|\tau_2|<100\epsilon$. Similarly, $\wG_4=L/4+(1/2)\ln(4/5)+\tau_4 + i\pi$ where $|\tau_4|<100\epsilon$. In particular $|\Im(\wG_2-i\pi)|, |\Im(\wG_4-i\pi)| < 100\epsilon$.

So the triangle with vertices $z,v_2,v_3$ satisfies $d(v_2,v_3) \ge \arccosh((3/2)(1-4\epsilon))$ and $|\angle(z,v_2,v_3)-\pi/2| < 100\epsilon$ and $|\angle(v_2,v_3,z)-\pi/2|<100\epsilon$. It is easy to see (or calculate using the law of cosines) that these conditions are contradictory if $\epsilon$ is small enough.

{\bf Cases 3 \& 7}: These cases can be handled similarly, so we will just prove case 3. This case is also very similar to the preceding case. Let $a,b,c,d$ be the points indicated in Case 3 figure \ref{fig:altshort}. Let $\Pi_1$ be the altitude plane containing $a$ and let $\Pi_2$ be the altitude plane containing $d$.

Let $\hG=(\gG_1,...,\gG_6)$ be the standardly oriented right-angled
hexagon satisfying
\begin{itemize}
\item $\gG_1$ contains $a$ and $b$;
\item $\gG_2$ contains the altitude perpendicular to $a$;
\item $\gG_4$ contains the altitude perpendicular to $d$;
\item $\gG_5 \supset \overline{cd}$;
\item $\gG_6$ contains $b$ and $c$.
\end{itemize}
If $\wG_i=\mu(\gG_{i-1},\gG_{i+1};\gG_i)$ (for all $i$ mod 6) then for $\epsilon$ small enough and $L$ large enough we have the following estimates (see lemmas \ref{lem:L/2 hexagon} and \ref{lem:altitudes})
\begin{itemize}
\item For $k=1,5$, $\wG_k = L/4 + \rho_k + i\pi$ where $\rho\in\C$ and $|\rho|<\epsilon$;
\item $\wG_6 = x_6\exp(-L/4)$ where $|x_6-2|<\epsilon$.
\end{itemize}

Suppose for a contradiction that $z$ is a point in $\Pi_1\cap
\Pi_2$. Let $v_k$ be the intersection point $\gG_k \cap \gG_{k+1}$ for
$k$ mod 6). Then the triangle with vertices $z, v_2, v_3$ has the following properties (for $L$ large enough).
\begin{itemize}
\item $v_2v_3 = |\Re(\wG_3)|$.
\item $|\angle(z,v_2,v_3)-\pi/2| \le |\Im(\wG_2-i\pi)|+\epsilon$.
\item $|\angle(v_2,v_3,z)-\pi/2| \le |\Im(\wG_4-i\pi)|+\epsilon$.
\end{itemize}
The $+\epsilon$ term in the above is to take into account the fact that each altitude plane is perpendicular to a short side rather than a long side. For example if $\Pi'_1$ is the plane containing the altitude through $a$ and is perpendicular to the side containing $a$ then the angle between $\Pi'_1$ and $\Pi_1$ is at most epsilon for $L$ large enough by lemma \ref{lem:altitudes} (this angle equals the imaginary part of $\wK_5-i\pi$ if $\hK$ is defined as in the lemma so that $\gK_5$ passes through $a$).

By the law of cosines 
\begin{eqnarray*}
\cosh(\wG_3) &=& \cosh(\wG_1)\cosh(\wG_5)+\sinh(\wG_1)\sinh(\wG_5)\cosh(\wG_6)\\
             &=& \exp(L/2 + \rho_1 + \rho_5) + O(1).
\end{eqnarray*}
So $\wG_3=L/2 + \rho_1 + \rho_5 + \log(2) + O(\exp(-L/2))$. By the law of sines
\begin{eqnarray*}
\sinh(\wG_2)&=&\frac{\sinh(\wG_5)\sinh(\wG_6)}{\sinh(\wG_3)}\\
              &=&\frac{-(1/2)\exp(L/4+\rho_5)x_6\exp(-L/4) + O(\exp(-L/2))  }{\exp(L/2+\rho_1+\rho_5)+O(1)}\\
& =&\frac{-(1/2)x_6  }{\exp(L/2+\rho_1)} + O(\exp(-L)) = O(\exp(-L/2)).
\end{eqnarray*}
So $\wG_2 = O(\exp(-L/2)) + i\pi$. Similarly, $\wG_4 = O(\exp(-L/2)) + i\pi$. So the triangle with vertices $z,v_2,v_3$ satisfies $d(v_2,v_3) \ge L/2 + \log(2) - 3\epsilon$ and $|\angle(z,v_2,v_3)-\pi/2| < 2\epsilon$ and $|\angle(v_2,v_3,z)-\pi/2|<2\epsilon$ if $L$ is large enough. It is easy to see (or calculate using the law of cosines) that these conditions are contradictory if $\epsilon$ is small enough and $L$ is large enough.

{\bf Cases 4 \& 9}: These cases are very similar so we will just do case 4.  Let $a,b,$ be the points indicated in Case 4 figure \ref{fig:altshort}. Let $\Pi_1$ be the altitude plane containing $a$ and let $\Pi_2$ be the altitude plane containing $b$. Suppose for a contradiction that there is a point $z \in \Pi_1 \cap \Pi_2$. Consider the triangle with vertices $a,b,z$. Note that $d(a,b) \ge L/2-\epsilon$ and $|\angle(z,a,b)-\pi/2|<\epsilon$ and $|\angle(a,b,z)-\pi/2|<\epsilon$ if $L$ is large enough by lemma \ref{lem:altitudes}. To elaborate, for example, $|\angle(a,b,z)-\pi/2|$ is at most $|\alpha|$ where $\alpha$ is the angle between $\Pi_2$ and $\Pi'_2$ where $\Pi'_2$ is the plane containing the altitude through $b$ and perpendicular to the side containing $b$. This angle equals $\wK_5$ if the pentagon $\hK$ is defined as in lemma \ref{lem:altitudes} so that $\gK_5$ passes through $b$.

The properties of the triangle $a,b,z$ given above are contradictory if $L$ is large enough and $\epsilon$ is less than $\pi/2$. Thus $\Pi_1$ and $\Pi_2$ are disjoint.

\end{proof}

\begin{figure}[htb]
 \begin{center}
 \ \psfig{file=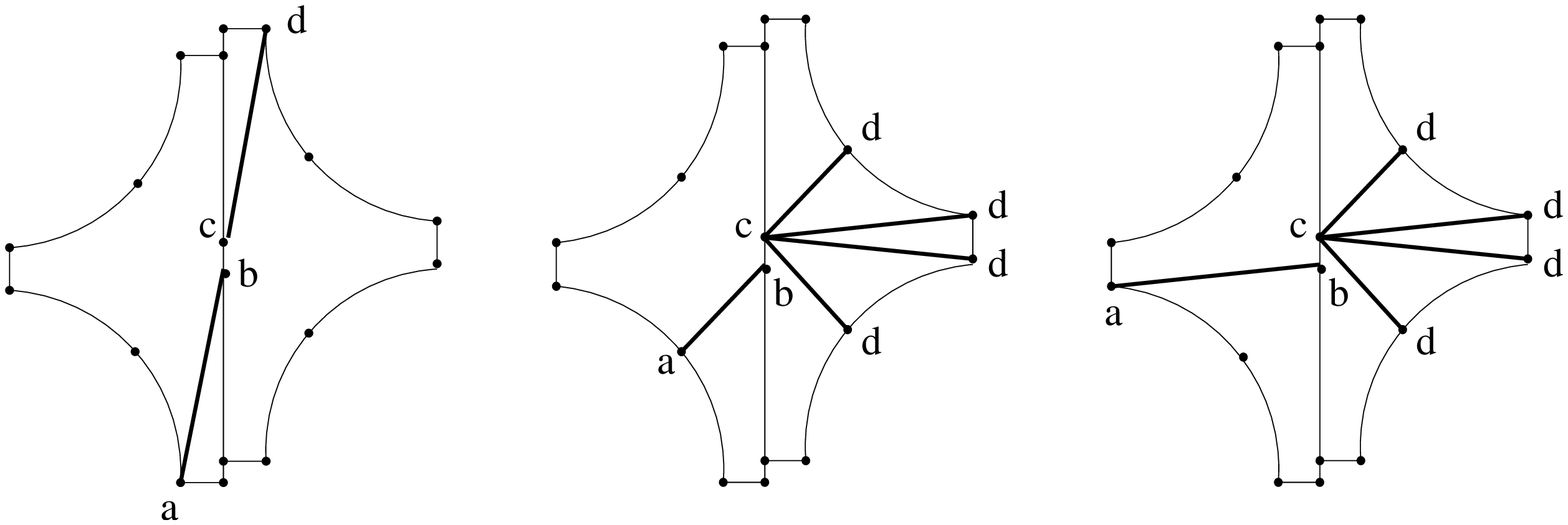,height=2in,width=4in}
 \caption{Possible configurations of $e_1$ and $e_2$}
 \label{fig:longside}
 \end{center}
 \end{figure}

\begin{figure}[htb]
 \begin{center}
 \ \psfig{file=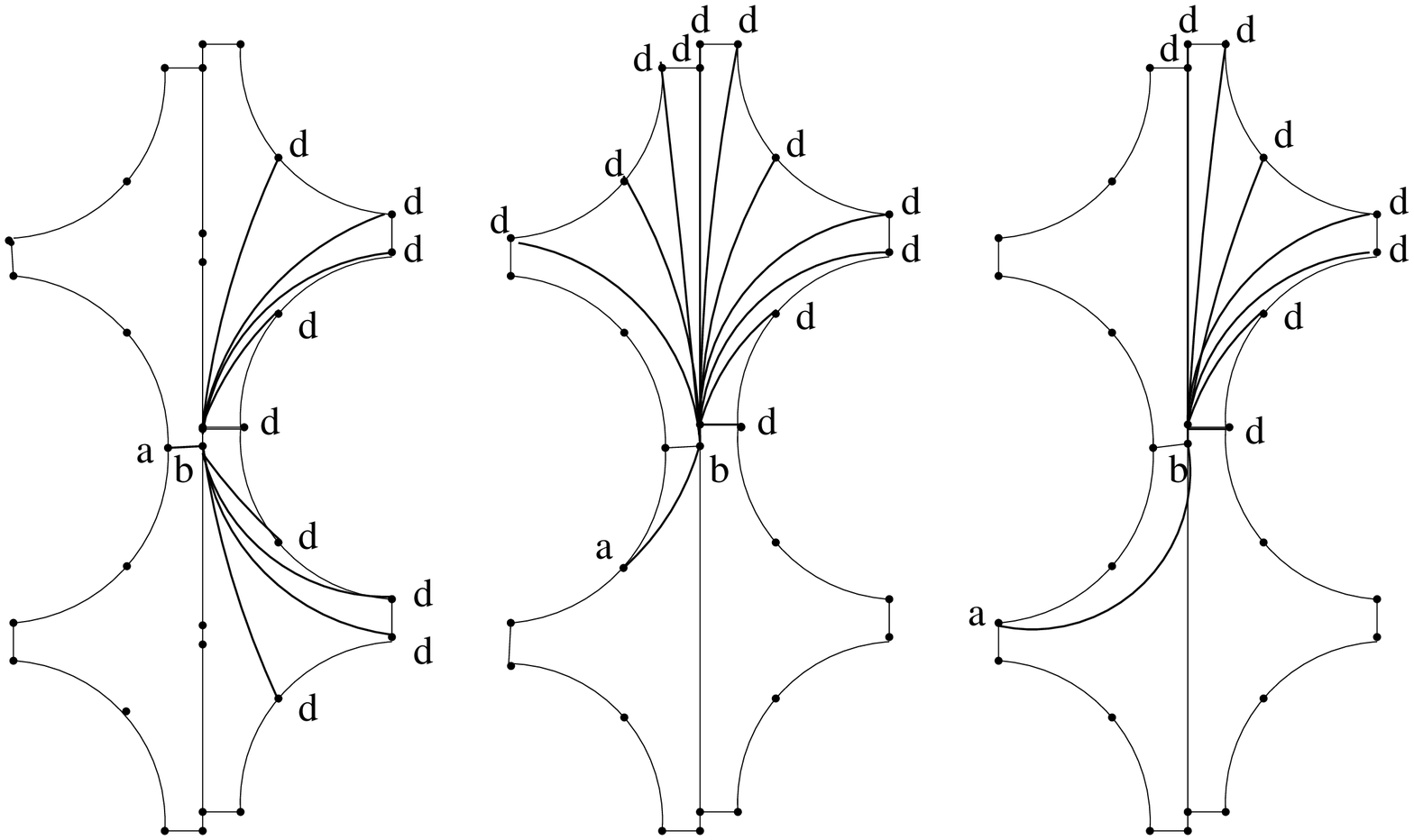,height=2.7in,width=4in}
 \caption{Possible configurations of $e_1$ and $e_2$}
 \label{fig:catty-1}
 \end{center}
 \end{figure}

\begin{figure}[htb]
 \begin{center}
 \ \psfig{file=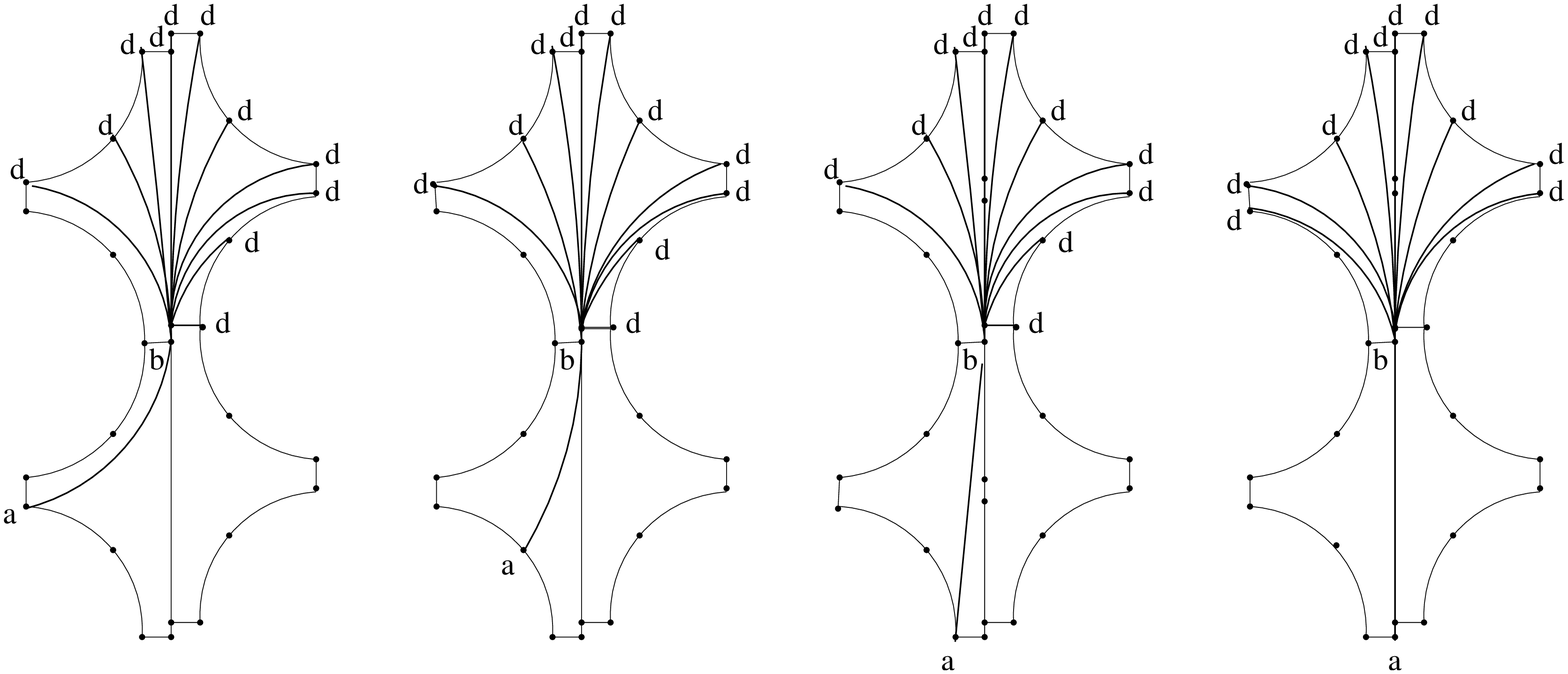,height=3in,width=5.6in}
 \caption{Possible configurations of $e_1$ and $e_2$}
 \label{fig:catty}
 \end{center}
 \end{figure}

The next step is to determine how the least weight property of ${\hat \gamma}$ affects the local geometry of $\gamma$. Let $e_1,e_2$ be consecutive positive-weight edges of $\gamma$. Let $(a,b)$ be the endpoints of $e_1$ and $(c,d)$ be the endpoints of $e_2$ so that $b$ is closest to $e_2$ and $c$ is closest to $e_1$. It may be that $b=c$. The possibilities for $a,b,c,d$ are described in figures \ref{fig:longside}-\ref{fig:catty}. In each figure there is only point $a$ and only one $b$. But there may be many different $c$'s and $d$'s (although $c$ is not always labeled). The interpretation is that if $a$ and $b$ are as in the figure than $c$ and $d$ could be any of the possibilities shown but there are no other possibilities for $c$ and $d$. For example, the leftmost example in figure \ref{fig:longside} shows that if $e_1$ is Type 2 and $b$ is an endpoint of an altitude then there is only one possibility for $e_2$. In particular $e_2$ most be of Type 2. In the middle example of figure \ref{fig:longside}, $e_1$ is of type 7 and there are four different possibilities for $e_2$. Using lemmas \ref{lem:L/2 hexagon}, \ref{lem:altitudes} and \ref{lem:midpoint distances} it can be checked that figures \ref{fig:longside}-\ref{fig:catty} show all the possibilities for $e_1$ and $e_2$ (up to some obvious symmetries).

\begin{proof}(of lemma \ref{lem:gamma})
The proof now follows from lemma \ref{lem:altitude planes} above and the list of possibilities in figures \ref{fig:longside}-\ref{fig:catty}. 
\end{proof}

\begin{proof}(of theorem \ref{thm:incompressible})
From lemma \ref{lem:gamma} above and the figures \ref{fig:longside}-\ref{fig:catty} it follows that if $i<j<k$, $\Pi(e_i),\Pi(e_j),\Pi(e_k)$ are nonempty and there does not exist $l$ with $i<l<k$, $l\ne j$ and $\Pi(e_l)$ nonempty, then $\Pi(e_j)$ separates $\Pi(e_i)$ from $\Pi(e_k)$. From this it follows more generally that if $i<j<k$ and $\Pi(e_i),\Pi(e_j),\Pi(e_k)$ are nonempty then $\Pi(e_j)$ separates $\Pi(e_i)$ from $\Pi(e_k)$. But this implies that $\gamma$ is not a closed curve. Therefore $j({\hat \gamma})$ is homotopically nontrivial. Since ${\hat \gamma}$ is an arbitrary homotopically nontrivial curve in $S$, $j$ is $\pi_1$-injective.  

\end{proof}


\begin{thebibliography}{99}






\bibitem[Buser]{Bus} P. Buser, \textit{Geometry and Spectra of
  Riemann Surfaces}. Progress in Mathematics, 106. Birkhauser Boston
  Inc., Boston, MA, 1992. xiv + 454pp.

\bibitem[Ehrenpreis]{Ehr} L. Ehrenpreis, \textit{Cohomology with
  Bounds}. 1970. Symposia Mathematica. Vol IV (INDAM Rome 1968-69). pp389-395. 

\bibitem[Gendron]{Gen} T. Gendron, \textit{The Ehrenpreis Conjecture
  and the Moduli-Rigidity Gap}, Complex Manifolds and Hyperbolic
  Geometry (Guanojuato 2001), 207-229.


\bibitem[Hedlund1]{Hed} G.A. Hedlund, \textit{Fuchsian groups and Transitive Horocycles}, Duke Math J., (2), 530-543, 1936.

\bibitem[Hedlund2]{Hed2} G.A. Hedlund, \textit{Dynamics of Geodesic Flows}, Bull. Am. Math. Soc. 40, (1939), 241-260. 

\bibitem[Fenchel]{Fen} W. Fenchel, \textit{Elementary Geometry in Hyperbolic Space}, Walter de Gruyter, Berlin, 1989.

\bibitem[Ratcliffe]{Rat} J. Ratcliffe, \textit{Foundations of Hyperbolic
 Manifolds}, Springer-Verlag, New York, 1994.

\bibitem[Ratner]{Ra5} M. Ratner, \textit{Raghunathan's topological conjecture and distributions of unipotent flows}, Duke J. Math. (63), 235-280, 1991.


\bibitem[Waldhausen]{Wal} F. Waldhausen, \textit{The Word Problem in
  Fundamental Groups of Sufficiently Large $3$-Manifolds}. Ann. of
  Math. 88 (1968), pp272-280.


\end{thebibliography}
\end{document}